
\def\misajour{January 24, 2004}

\def\title{Open Diophantine Problems}

\magnification=\magstephalf
 

\catcode`\Œ=\active\defŒ{{\aa}}       
\catcode`\=\active\def{\c c}        
\catcode`\'=\active\def'{\c C}        


\catcode`\Š=\active\defŠ{\"a}        
\catcode`\'=\active\def'{\"e}        
\catcode`\•=\active\def•{\"{\i}}     
\catcode`\š=\active\defš{\"o}        
\catcode`\Ÿ=\active\defŸ{\"u}        
\catcode`\…=\active\def…{\"O}        
\catcode`\†=\active\def†{\"U}        
\catcode`\‡=\active\def‡{\'a}        
\catcode`\Ž=\active\defŽ{\'e}        
\catcode`\'=\active\def'{\'{\i}}     
\catcode`\—=\active\def—{\'o}        
\catcode`\œ=\active\defœ{\'u}        
\catcode`\ƒ=\active\defƒ{\'E}        
\catcode`\ˆ=\active\defˆ{\`a}        
\catcode`\=\active\def{\`e}        
\catcode`\"=\active\def"{\`{\i}}     
\catcode`\˜=\active\def˜{\`o}        
\catcode`\=\active\def{\`u}        
\catcode`\Ë=\active\defË{\`A}        
\catcode`\‹=\active\def‹{\~a}        
\catcode`\–=\active\def–{\~n}        
\catcode`\›=\active\def›{\~o}        
\catcode`\Ì=\active\defÌ{\~A}        
\catcode`\"=\active\def"{\~N}        
\catcode`\Í=\active\defÍ{\~O}        
\catcode`\‰=\active\def‰{\^a}        
\catcode`\=\active\def{\^e}        
\catcode`\"=\active\def"{\^{\i}}     
\catcode`\™=\active\def™{\^o}        
\catcode`\ž=\active\defž{\^u}        

\font\tenmsb=msbm10
\font\sevenmsb=msbm7
\font\fivemsb=msbm5
\newfam\msbfam
\textfont\msbfam=\tenmsb
\scriptfont\msbfam=\sevenmsb
\scriptscriptfont\msbfam=\fivemsb
\def\Bbb#1{{\fam\msbfam\relax#1}}

\font\seventt=cmtt10 scaled 700
 
\font\eightt=cmtt10 scaled 800
\font\ninett=cmtt10 scaled 900 
\font\nineit=cmsl10 at 9pt

    
\font\twelvecmssdc=cmssdc10 at 12 pt
\font\ninecmssi=cmssi10 at 9 pt
\font\twelvecmssi=cmssi10 at 12 pt

\font\niness=cmss10 at 9pt
\font\sevenss=cmss10 at 7pt
\font\fivess=cmss10 at 5pt

\font\tenss=cmss10    

\font\goth=eufm10

\newfam\ssfam
\scriptscriptfont\ssfam=\fivess
\textfont\ssfam=\tenss \scriptfont\ssfam=\sevenss
\def\sanserif{\fam\ssfam\tenss}

\def\textesanserif{
\textfont0=\tenss
\textfont1=\teni
\font\ninerm=cmss10 at 9pt
\font\it=cmssi10
\let\sl=\it
\sanserif}

\textesanserif

\def\petitsanserif
{
\textfont0=\niness
\textfont1=\ninecmssi
\font\bf=cmssdc10  at 9 pt   
\font\it=cmssi10 at 9 pt 
\let\sl=\it
\def\sanserif{\fam\ssfam\niness}
\sanserif
}

\font\tenbfit=cmbxti10

\def\date{
{\ifcase\month\or January\or February\or March\or April
\or May \or June\or July\or August\or September
\or October\or November\or December\fi}\
{\the\day},\
{\the\year}}
\newcount\hour\newcount\minute	\hour=\time
\divide\hour by 60
\minute=-\hour \multiply\minute by 60 \advance\minute by \time  
\def\heure{\the\hour h\kern -.2mm \ifnum\minute<10%
0\fi\the\minute} 


\def\og{\leavevmode\raise.3ex\hbox{$\scriptscriptstyle
\langle\!\langle$}}

\def\fg{\leavevmode\raise.3ex\hbox{$\scriptscriptstyle
\,\rangle\!\rangle$}}

\def\bfabc{ \hbox{\tenbfit  abc}}

\def\calA{{\cal A}}
\def\bC{{\Bbb C}}
\def\calC{{\cal C}}
\def\calD{{\cal D}}
\def\calE{{\cal E}}
\def\bfE{ \hbox{\tenbfit  E}}
\def\calF{{\cal F}}
\def\bF{{\Bbb F}}
\def\bfG{ \hbox{\tenbfit  G}}
\def\bG{{\Bbb G}}
\def\rmh{{\rm h}}
\def\rmH{{\rm H}}
\def\gI{\hbox{\goth I}\relax}
\def\calL{{\cal L}}
\def\calLtilde{\widetilde{\calL}}
\def\ttM{{\tt M}}
\def\rmM{{\rm M}}
\def\bN{{\Bbb N}}
\def\bP{{\Bbb P}}
\def\gP{\hbox{\goth P}\relax}
\def\bQ{{\Bbb Q}}
\def\Qbar{{\overline{\Bbb Q}}}
\def\Qbaretoile{{\overline{\Bbb Q}}^\times}
\def\bR{{\Bbb R}}
\def\calK{{\cal K}}
\def\bZ{{\Bbb Z}}
\def\gZ{\hbox{\goth Z}}
\def\del{\partial}
\def\bGm{\bG_{\rm m}}
\def\ualpha{{\underline{\alpha}}}
\def\ugamma{{\underline{\gamma}}}
\def\utheta{{\underline{\theta}}}
\def\kappatilde{\widetilde{\kappa}}
\def\pps{\,\colon\,}   
\def\and{\quad\hbox{ and }\quad}
\def\for{\quad\hbox{ for }\quad}
\def\with{\quad\hbox{ with }\quad}
\def\where{\quad\hbox{ where }\quad}

\def\lcm{{\rm lcm}}

\def\HPoincare{\hbox{\goth H}\relax}

\def\cqfd{\unskip\kern 6pt\penalty 500
\raise -2pt\hbox{\vrule\vbox to 10pt{\hrule width 4pt
\vfill\hrule}\vrule}\par}
\def\virgule{\raise2pt\hbox{\rm ,}}
\def\pointvirgule{\raise2pt\hbox{;}}

\def\adots{\mathinner{\mkern2mu\raise1pt\hbox{.}
\mkern3mu\raise4pt\hbox{.}\mkern1mu\raise7pt\hbox{.}}}

\def\trdeg{{\rm trdeg}\relax}

\def\dist{ {\rm dist }\relax}

\def\deuxFun{\hbox{$_2F_1$}}

\def\deuxfun#1#2#3#4{\deuxFun\left(
\matrix{#1\;,\; #2\cr#3\cr}
\Bigr|#4\right)}

\def\codim{ {\rm codim }\relax}

\font\xrmhuit=wncyr10 at 8 pt

\def\Beta{\hbox{\xrmhuit V} }
\def\ZH{\hbox{\xrmhuit ZH} }

\def\house#1{\setbox1=\hbox{$\,#1\,$}%
\dimen1=\ht1 \advance\dimen1 by 2pt \dimen2=\dp1 \advance\dimen2 by 2pt
\setbox1=\hbox{\vrule height\dimen1 depth\dimen2\box1\vrule}%
\setbox1=\vbox{\hrule\box1}%
\advance\dimen1 by .4pt \ht1=\dimen1
\advance\dimen2 by .4pt \dp1=\dimen2 \box1~ \relax}

\def\th{\noalign{\hrule}\cr}  
\def\tvi{\vrule height 13pt depth 8pt width 0pt}
\def\tv{\tvi \vrule} 
\def\cc#1{\ #1 \ }
 
\def\grandt{
\hskip -1.15  cm   
{^{^|}}^{\strut} \hskip 1.5  cm   
\hskip -2.95  cm   
{\raise -.3 true cm\hbox{$\vdots$} }\!
 ^{\vrule height .1 pt depth .1pt width 2.6  cm } 
\! {\raise -.3 true cm\hbox{$\vdots$} }
\hskip -2.45  cm} 

\def\frac#1#2{{#1\over #2}}

\def\Numero#1{ N$^\circ$~{#1}}

\def\us{\underline{s}}
\def\ux{\underline{x}}

\def\uz{\underline{z}}
\def\uzeta{\underline{\zeta}}

\def\and{\quad\hbox{and}\quad}
\def\Li{{\rm Li}}
\def\Gbar{\overline{G}}

\let\sl=\cmssi
\let\bf=\cmssdc  


\headline={
\hfill}

\output={\ifnum\pageno>1\headline={\sevenbf \hfil
\folio}\fi\plainoutput}

\footline={\ifnum\pageno>1 \hfill
\eightt
{}\misajour
\else 
\hfill{$
\hbox{ \ninett
$<$http$:$//www.math.jussieu.fr/${ \sim}$miw/%
articles/pdf/odp.pdf$>$ \strut}$}
\fi}

\def\section#1{
\vskip 1 true cm plus .5 true cm minus .5  true cm
\goodbreak 
\noindent{\twelvecmssdc #1}
\vskip .5 true cm plus .25 true cm minus .25  true cm}

\def\subsection#1{
\vskip .75 true cm plus .3 true cm minus .3  true cm
\goodbreak 
\noindent{\bf #1}
\vskip .375 true cm plus .15 true cm minus .15  true cm}


\def\sectionzero{ Abstract}

\def\sectionun{ Diophantine Equations}

\def\subsectionunpointun{Points on Curves}
\def\subsectionunpointdeux{Exponential Diophantine Equations}
\def\subsectionunpointtrois{Markoff Spectrum}

\def\sectiondeux{Diophantine Approximation}

\def\subsectiondeuxpointun{The $\bfabc$ Conjecture}
\def\subsectiondeuxpointdeux{Thue-Siegel-Roth-Schmidt}
\def\subsectiondeuxpointtrois{Irrationality and Linear Independence Measures}

\def\sectiontrois{Transcendence}

\def\subsectiontroispointun{Schanuel's Conjecture}
\def\subsectiontroispointdeux{Multiple Zeta Values}
\def\subsectiontroispointtrois{Gamma, Elliptic, Modular, $\bfG$ and $\bfE$-Functions}
\def\subsectiontroispointquatre{Fibonacci and Miscellanea}

\def\sectionquatre{Heights}

\def\subsectionquatrepointun{Lehmer's Problem}
\def\subsectionquatrepointdeux{Wirsing-Schmidt Conjecture}
\def\subsectionquatrepointtrois{Logarithms of Algebraic Numbers}
\def\subsectionquatrepointquatre{Density: Mazur's Problem}

\def\sectioncinq{Further Topics}

\def\subsectioncinqpointun{Metric Problems}
\def\subsectioncinqpointdeux{Function Fields}
 


\def\numeroProblemeun{1.1}
\def\numeroConjectureCatalan{1.2}
\def\numeroConjecturePillai{1.3}
\def\numeroConjectureShorey{1.4}
\def\numeroConjectureSchlickeweiViola{1.5}
\def\numeroConjectureGrimm{1.6}
\def\numeroConjectureGrimmFaible{1.7}
\def\numeroConjectureGrimmLangevin{1.8}

\def\numeroconjectureMarkoff{1.9}
\def\numeroapproximationMarkoff{1.10}
\def\numeroconjectureabc{2.1}
\def\numeroconjectureErdosWoods{2.2}
\def\numeroConjectureErdosDressler{2.3}
\def\numeroConjecturePhilipppon{2.4}
\def\numeroConjectureLangW{2.5}
\def\ConjecturePillaiQuantitative{2.6}
\def\numeroConjectureHall{2.7}
\def\numeroThueSiegelRoth{2.8}
\def\numeroQuestionBPQ{2.9}
\def\numeroQuestionUPQ{2.10}
\def\numeroinegalitepourRothRaffine{2.11}
\def\numeroRothRaffine{2.12}
\def\numeroCatalan{2.13}
\def\numeroProblemeMahler{2.14}
\def\numeroconjectureMahler{2.15}

\def\numeroConjectureSchanuel{3.1}
\def\numeroConjectureRoyRome{3.2}
\def\numeroConjectureialogs{3.3}
\def\numeroConjectureQuatreExponentiellesForte{3.4}
\def\numeroConjectureCinqExponentiellesForte{3.5}
\def\numeroConjectureDRoyA{3.6}
\def\numeroConjectureQuatreExponentielles{3.7}
\def\numeroProblemeGelfond{3.8}
\def\numeroProblemeSchneider{3.9}
\def\numeroProblemeGelfondSchneider{3.10}
\def\numeroLindemannWeierstrasspadic{3.11}
\def\numeroGelfondpadic{3.12}
\def\numeroBlumCuckerShubSmale{3.13}
\def\numeroConjectureLemmeSchwarz{3.14}
\def\numeroConjectureGoncharov{3.15}
\def\numeroConjectureZagier{3.16}
\def \numeroConjectureZetaRiemann{3.17}
\def\numeroquestionGammaUnCinquieme{3.18}
\def\numeroConjectureRohrlich{3.19}
\def\numeroConjectureNesterenko{3.20}
\def\numeroConjectureBertolin{3.21}
\def\numeroConjectureAndreOort{3.22}
\def\numeroIndependancePeriodesElliptiques{3.23}
\def\numeroConjectureBertrandA{3.24}
\def\numeroConjectureBertrandB{3.25}
\def\numeroConjectureHilbertIrreducibility{3.26}

\def\numeroconjectureMahlerLiouville{3.27}

\def\numeroconjectureLehmer{4.1}
\def\numeroconjectureSchinzelZassenhaus{4.2}
\def\numeroconjectureLehmerSimultanne{4.3}
\def\numeroconjectureAmorosoDavid{4.4}
\def\numeroconjectureAmorosoDavidB{4.5}
\def\numeroConjectureLangevin{4.6}
\def\numeroConjectureLehmerVA{4.7}
\def\numeroFormuleDirichlet{4.8}
\def\numeroFormuleTiroirs{4.9}
\def\numeroTheoremeWirsing{4.10}
\def\numeroConjectureWirsingSchmidt{4.11}
\def\numeroPremiereConjectureGraz{4.12}
\def\numeroConjectureLaurentRoy{4.13}
\def\numeroDeuxiemeConjectureGraz{4.14}

\def\numeroMesureIndLinLo{4.15}
\def\numeroConjectureSchanuelEffectif{4.16}
\def\numeroQuestionMazur{4.17}
\def\numeroQuestionMazurVA{4.18}
\def\numeroconjectureNagoya{4.19}
\def\numeroConjectureDensiteTores{4.20}
\def\numeroConjectureDensiteToresQuantitative{4.21}

\def\numeroConjectureDuffinSchaeffer{5.1}
\def\numeroquestionBugeaud{5.2}

\def\numeroQuatriemeConjectureGraz{5.3}


 \def\numeroConjectureLoxtonVdP{5.4}


\def\refAlloucheShallit{AS-2003}
\def\refAmbros-et-al{AmEMW~2000}
\def\refAmorosoDavidCrelle{AmoD~1999}
\def\refAmorosoDavidAA{AmoD~2000}
\def\refAndersonThakur{AnT~1990}
\def\refAndreLivre{And~1989}
\def\refAndre{And~1997}
\def\refBakerLivre{B~1990}
\def\refBakerEger{B~1998}
\def\refBLSW{BaLSW~1996}
\def\refBallRivoal{BalR~2001}
\def\refBernikDodson{BeD~1999}
\def\refBertolin{Ber~2002}
\def\refBertrandAustralie{Bert~1997a}
\def\refBertrandRamanujanJ{Bert~1997b}
\def\refBilu{Bi~2003}
\def\refBombieriLang{BoL~1970}
\def\refBoydExperimental{Boy~1998}
\def\refBoydZakopane{Boy~1999}
\def\refBriggs{Br~2002}
\def\refBrowkin{Bro~1999}
\def\refBrownawell{Brow~1998}
\def\refBugeaudGraz{Bu~2000a}
\def\refBugeaudJNT{Bu~2000b}
\def\refBugeaudLivre{Bu~2004}
\def\refBugeaudMignotte{BuM~1999}
\def\refBugeaudShorey{BuSh~2001}
\def\refBundschuh{Bun~1979}
\def\refCartierBourbaki{C~2001}
\def\refCassels{Ca~1957}
\def\refCatalan{Cat~1844}
\def\refChambertLoir{Ch~2001}
\def\refChudnovsky{Chu~1980}
\def\refPaula{Co~2003}
\def\refCohn{Coh~1993}
\def\refColliot{ColSS~1997}
\def\refCusickFlahive{CuF~1989}
\def\refCusickPomerance{CuP~1984}
\def\refDavenportSchmidt{DS~1969}
\def\refDavidSinnouLille{Da~2001a}
\def\refDavidSinnou{Da~2001b}
\def\refDavidHindry{DaH~2000}
\def\refDavidPhilipppon{DaP~1999}
\def\refDavisMatiyasevichRobinson{DavMR~1976}
\def\refDeMathanLasjaunias{dML~1999}
\def\refDiazJNTBx{Di~1997}
\def\refDiazRAMA{Di~2000}
\def\refDobrowolskiAA{Do~1979}
\def\refDujella{Du~2004}
\def\refErdosMonthly{E~1980}
\def\refErdosDurham{E~1988}
\def\refFaltings{F~2000}
\def\refFeldmanNesterenko{FeN~1998}
\def\refFeldmanShidlovskii{FeS~1967}
\def\refFinch{Fi~2003}
\def\refFischler{Fis~2002}
\def\refFlattoLagariasPollington{FlLP~1995}
\def\refGalochkin{G~1983}
\def\refGelfondCRAS{Ge~1934}
\def\refGelfondConjecture{Ge~1949}
\def\refGelfondLivre{Ge~1952}
\def\refGramain{Gr~1981}
\def\refGramainWeber{GrW~1985}
\def\refGras{Gra~2002}
\def\refGrimm{Gri~1969}
\def\refGrinspan{Grin~2002}
\def\refGuinness{Gu~2000}
\def\refGuy{Guy~1994}
\def\refGyarmati{Gy~2001}
\def\refHall{H~1971}
\def\refHardyWright{HaW~1979}
\def\refHarman{Har~1998}
\def\refHershonskyPaulin{HeP~2001}
\def\refHilbert{Hi~1900}
\def\refClayMathematicalInstitute{J~2000}
\def\refKontsevichZagier{KZ~2000}
\def\refKraus{Kr~1999}
\def\refLagarias{L~2001}
\def\refLangITN{La~1966}
\def\refLangTNDA{La~1971}
\def\refLangHDDP{La~1974}
\def\refLangDistributions{La~1978a}
\def\refLangECDA{La~1978b}
\def\refLangCyclotomic{La~1978c}
\def\refLangONCDI{La~1990}
\def\refLangEncyclopaedie{La~1991}
\def\refLangAlgebra{La~1993}
\def\refLangGazette{La~1996}
\def\refLangevin{Lan~1977}
\def\refLangevinLuminy{Lan~1992}
\def\refLangevinRocky{Lan~1996}
\def\refLangevinJNTBx{Lan~2001}
\def\refLaurentEger{Lau~1998}
\def\refLaurentSalonique{Lau~1999a}
\def\refLaurentTokyo{Lau~1999b}
\def\refLeopoldt{Le~1962}
\def\refMahlerZnumbers{M~1968}
\def\refMahlerLivre{M~1976}
\def\refMahlerSuggestions{M~1984}
\def\refManin{Ma~1977}
\def\refMasserabc{Mas~1990}
\def\refMatiyasevich{Mat~1999}
\def\refMatveev{Matv~2000}
\def\refBeal{Mau~1997}
\def\refMazurTRP{Maz~1992}
\def\refMazurQDUNT{Maz~1994}
\def\refMazurSTRP{Maz~1995}
\def\refMazurPowersofNumbers{Maz~2000}
\def\refMetsankyla{Me~2004}
\def\refMihailescu{Mi~2004}
\def\refMullerTisserand{M\"uT~1996}
\def\refNarkiewiczCPNT{N~1986}
\def\refNarkiewiczEATAN{N~1990}
\def\refNesterenkoPhilippon{NeP~2001}
\def\refNitaj{Ni}
\def\refOesterle{\OE~1988}
\def\refOort{Oo~1997}
\def\refPhilipponLW{P~1987}
\def\refPhilipponJNT{P~1997}
\def\refPhilipponAustralie{P~1999a}
\def\refPhilipponGammaunquart{P~1999b}
\def\refPillai{Pi~1945}
\def\refPollingtonVelani{PoV~2000}
\def\refDipendra{Pr~2004}
\def\refRamachandraAA{R~1968}
\def\refRauzy{Ra~1976}
\def\refRemond{Re~2000}
\def\refRhinSmyth{RhS~1995}
\def\refRibenboimCatalan{Ri~1994}
\def\refRibenboimNumbersFriends{Ri~2000}
\def\refRivoalCRAS{Riv~2000}
\def\refRoySchanuel{Ro~1989}
\def\refRoyMatrices{Ro~1990}
\def\refRoyInventiones{Ro~1992}
\def\refRoyActaMath{Ro~1995}
\def\refRoyACVEF{Ro~2001a}
\def\refRoyAASNL{Ro~2001b}
\def\refRoyIFAF{Ro~2002}
\def\refRoyRMTNT{Ro~2001c}
\def\refRoyNote{Ro~2003}
\def\refSchinzelBilkent{S~1999}
\def\refSchinzelZassenhaus{SZ~1965}
\def\refSchlickeweiViola{ScV~2000}
\def\refSchmidtLNDA{Sch~1980}
\def\refSchmidtBilkent{Sch~1999}
\def\refSchmidtSurveyAA{Sch~2000}
\def\refSchneiderLivre{Schn~1957}
\def\refSerre{Se~1989}
\def\refShoreySurvey{Sh~1999}
\def\refShoreyConjectures{Sh~2000}
\def\refShoreyTijdeman{ShT~1986}
\def\refSiegelPreuss{Si~1929}
\def\refSiegelLivre{Si~1949}
\def\refSierpinskiA{Sie~1964}
\def\refSierpinskiB{Sie~1970}
\def\refSilverman{Sil~1986}
\def\refSmale{Sm~1998}
\def\refSondow{So~2003}
\def\refSprindzuk{Sp~1979}
\def\refStewartYuKunruiI{StY~1991}
\def\refStewartYuKunruiII{StY~2001}
\def\refTerasoma{T~2002}
\def\refThakur{Th~1998}
\def\refTijdemanHilbert{Ti~1976a}
\def\refTijdemanCatalan{Ti~1976b}
\def\refTijdemanSurveyA{Ti~1989}
\def\refTijdemanSurveyB{Ti~1998}
\def\refTijdemanSADA{Ti~2000}
\def\refVojta{V~2000}
\def\refmiwSemLelong{W~1976}
\def\refmiwGAGDT{W~1986}
\def\refmiwNagoya{W~1999a}
\def\refmiwSurveyIA{W~1999b}
\def\refmiwGraz{W~2000a}
\def\refmiwGL{W~2000b}
\def\refmiwCinquante{W~2000c}
\def\refmiwMZV{W~2000d}
\def\refmiwCetraro{W~2001}
\def\refWaring{Wa~1770}
\def\refWhittakerWatson{WhW~1927}
\def\refWirsing{Wi~1961}
\def\refWong{Wo~1999}
\def\refWoods{Woo~1981}
\def\refZagier{Z~1994}
\def\refZudilin{Zu~2003}


\centerline{\twelvecmssdc \title}
\medskip

\centerline{\sanserif
 by}
\medskip

\centerline{\twelvecmssi  Michel WALDSCHMIDT
\footnote{}{\ninett  
\hskip -1 true cm
2000 Mathematics Subject Classification: 
Primary: 11Jxx \ \
Secondary: 11Dxx, 11Gxx, 14Gxx}
}

\vskip 1 true cm plus .5 true cm minus .5  true cm

\hfill
\hbox{\it DŽdiŽ ˆ un jeune septuagŽnaire, Pierre Cartier, qui m'a beaucoup appris.}


\goodbreak

\section{Contents}


\item{\S 1.} \sectionun

\itemitem{\S 1.1} \subsectionunpointun

\itemitem{\S 1.2} \subsectionunpointdeux

\itemitem{\S 1.3} \subsectionunpointtrois

\item{\S 2.} \sectiondeux

\itemitem{\S 2.1} The $abc$ Conjecture 

\itemitem{\S 2.2} \subsectiondeuxpointdeux

\itemitem{\S 2.3} \subsectiondeuxpointtrois

\item{\S 3.} \sectiontrois

\itemitem{\S 3.1} \subsectiontroispointun

\itemitem{\S 3.2} \subsectiontroispointdeux

\itemitem{\S 3.3} Gamma, Elliptic, Modular, $G$ and $E$-Functions

\itemitem{\S 3.4} \subsectiontroispointquatre

\item{\S 4.} \sectionquatre

\itemitem{\S 4.1} \subsectionquatrepointun

\itemitem{\S 4.2} \subsectionquatrepointdeux

\itemitem{\S 4.3} \subsectionquatrepointtrois

\itemitem{\S 4.4} \subsectionquatrepointquatre

\item{\S 5.} \sectioncinq

\itemitem{\S 5.1} \subsectioncinqpointun

\itemitem{\S 5.2} \subsectioncinqpointdeux

\item{ } \hskip -.7 true cm References


\goodbreak

\section{\sectionzero}


Diophantine Analysis is a very active domain of mathematical research 
where one finds more  conjectures than results.

We collect here a number of open questions concerning Diophantine
equations (including 
Pillai's Conjectures), Diophantine
approximation (featuring the $abc$ Conjecture) and
transcendental number theory (with, for instance, Schanuel's
Conjecture). Some questions related to Mahler's measure
and Weil absolute logarithmic height are then considered
(e.g. Lehmer's Problem). We also discuss Mazur's question on
the density of rational points on a variety, especially in
the particular case of algebraic groups, in connexion with
transcendence problems in several variables. We say only a
few words on metric problems,  equidistribution questions,
 Diophantine approximation on manifolds and Diophantine
analysis on function fields.

{

\narrower
\nineit
This survey grew out of lectures given in several places including India
(October 2000 and September 2002), Ivory Coast (February 2001), Italy (April
2001), Canada (May 2001), Lebanon (November 2002) and France. The author is grateful to all colleagues
who gave  him the opportunity to speak on this topic, and also to all those who
contributed by relevant remarks and suggestions.

}


\goodbreak

\section{1. \sectionun}
 


\subsection{1.1. \subsectionunpointun}

Among the 23  problems of
Hilbert [\refHilbert], [\refGuinness] the tenth one has the shortest statement:

\smallskip
\item{}{\it  Given a Diophantine equation with any number of unknown
quantities and with integral numerical coefficients: To devise a process
according to which it can be determined by a finite number of operations 
whether the equation is solvable in rational integers.}

\smallskip\noindent
An equation of the 
form
$f(\ux)=0$, where $f\in\bQ[X_1,\ldots,X_n]$ is a given polynomial, while
the unknowns $\ux=(x_1,\ldots,x_n)$ take  
rational integer values, is a  {\it  Diophantine equation}. To solve this
equation amounts to determine the integer points on the corresponding
hypersurface of the affine space. Hilbert's tenth Problem is to give  an
algorithm which tells us whether such a given Diophantine equation has a solution
or not.
 
There are other types of Diophantine equations. First of all one may consider
rational solutions instead of integer ones: in this case, one considers
rational points on a hypersurface. Next one may consider integer or rational 
points over a number field. There is a situation which is intermediate between
integer and rational points, where the unknowns take {\it 
$S$-integral points values}. This means that $S$ is a fixed finite set of prime
numbers (rational primes, or prime ideals in the number field), and
the denominators of the solutions are restricted to belong to
$S$. Examples are the  Thue-Mahler equation 
$$
F(x,y)=p_1^{z_1}\cdots
p_k^{z_k}
$$
where $F$ is a homogeneous polynomial with integer coefficients and $p_1,\ldots,p_k$ are fixed primes (the
unknowns are
$x,y,z_1,\ldots,z_k$ and take rational integer values with $z_i\ge 0$) and the  generalized Ramanujan-Nagell equation
$x^2+D=p^n$ where $D$ is a fixed integer, $p$ a fixed prime, and the
unknowns are
$x$,
$n$ and take rational integer values with $n\ge 0$  (see for instance [\refShoreyTijdeman], [\refTijdemanSurveyB] 
[\refShoreySurvey] and [\refBugeaudShorey] for these and  other similar
questions).

Also, it is interesting to deal with simultaneous Diophantine 
equations, i.e., to investigate rational or integer points on algebraic
varieties.

 The final answer to Hilbert original 10th Problem has been given in 1970 
by Yu. Matiyasevich, after the works of  M.~Davis, H.~Putnam
and J.~Robinson.  This was the culminating stage of a rich and
beautiful theory (see  [\refDavisMatiyasevichRobinson],  [\refManin]  and
[\refMatiyasevich]).
The solution is negative: there is no hope nowadays to achieve
a complete theory of the subject. But one may still hope that
there is a positive answer if one restricts Hilbert's initial
question to equations in few variables, say $n=2$, which
amounts to considering integer points on a plane curve. In
this case, deep results have been achieved during the 20th
century and many results are known, but much more remains unveiled.

The most basic results are Siegel's (1929) and Faltings's 
(1983) ones. Siegel's Theorem deals with integer points and
produces an algorithm to decide whether the set of solutions
forms a finite or an infinite set. Faltings's result, solving
Mordell's Conjecture, does the same for rational solutions,
i.e., rational points on curves. To these two outstanding
achievements of the 20th century, one may add Wiles's
contribution which not only settles Fermat's last Theorem, but
also provides a quantity of similar results for other curves
[\refKraus].

\smallskip

Some natural questions  occur:

\smallskip

\item{(a)} To answer Hilbert's tenth Problem for this special case of plane 
curves, which means to give an algorithm to decide whether a given Diophantine
equation
$f(x,y)=0$ has a solution in $\bZ$  (and the same problem  in $\bQ$).

\smallskip

\item{(b)}
To give an upper bound for the number of either rational or integral  points
on a curve.

\smallskip

\item{(c)}
To give an algorithm for solving explicitly a given Diophantine equation in
two unknowns.

\smallskip
\noindent
Further questions may be asked. For instance in question b) one might ask for 
the  exact number of solutions; it may be more relevant to consider more
generally the number of points on any number field, or the number of points of
bounded degree and to investigate related generating series\dots\ The number of
open problems is endless!

Our goal here is not to describe in detail the state of the art on these
questions (see for instance [\refLangEncyclopaedie]). It will suffice to say 

\smallskip
\item{-} that a complete answer to question (a) is not yet available: there is no
algorithm so far (even conjectural ones) to decide whether a curve has a
rational point or not, 
\smallskip
\item{-} that a number of results are known on question (b), the
latest work on this topic being due to G. R\'emond [\refRemond] who produces
an effective upper bound for the number of rational points on a curve of
genus $\ge 2$, 
\smallskip
\item{-} and that
question (c) is unanswered even for integer points, and even for the special case
of curves of genus $2$. 
\smallskip
\noindent
We do not request a practical algorithm, but only (to
start with) a theoretical one. So our first open problem will be an
effective refinement to Siegel's Theorem.

\proclaim Problem {\numeroProblemeun}. 
 Let $f\in\bZ[X,Y]$ be a polynomial such that the equation $f(x,y)=0$ has only
finitely many solutions $(x,y)\in \bZ\times\bZ$. Give an upper bound for
$\max\{|x|,|y|\}$ when $(x,y)$ is such a solution, in terms of the degree of $f$
and of the maximal absolute value of the coefficients of $f$.

 That such a bound exists is part of the hypothesis, but the
problem is to give it explicitly (and, if possible, in a closed form).

Further similar questions might also be asked for more variables (rational
points on varieties), for instance Schmidt's norm form equations. We refer
the reader to  [\refLangHDDP] and 
[\refLangEncyclopaedie] for such questions, including the 
Lang-Vojta Conjectures.

\medskip
Even the simplest case of 
quadratic forms gives rise to open problems: the determination of  all
positive integers which are represented by a given binary form is far from 
being solved. Also it is expected that infinitely many real
quadratic fields have class number one, but it is not even known that
there are infinitely many number fields (without restriction on the degree) with class
number one. Recall that the first complete solution of Gauss's  class
number
$1$ and $2$ Problems (for imaginary quadratic fields) has been reached
by transcendence methods (A.~Baker and H.M.~Stark) so it  may be considered as a Diophantine
problem. Nowadays more efficient methods (Goldfeld, Gross-Zagier,\dots -- 
see [\refLangEncyclopaedie] Chap.~V, \S 5) are available.

A related open problem is the determination of Euler's  {\it  numeri idonei}
[\refRibenboimNumbersFriends]. Fix a positive integer $n$. If $p$ is an odd
prime for which there exist integers $x\ge 0$ and $y\ge 0$ with $p=x^2+ny^2$,
then

\smallskip
\item{(i)} $\gcd(x,ny)=1,$ 

\item{(ii)} the equation $p=X^2+nY^2$ in integers $X\ge 0$
and
$Y\ge 0$  has the only solution $X=x$ and $Y=y$. 

\smallskip\noindent
Now let $p$ be an odd integer such that there exist integers $x\ge 0$ and $y\ge
0$ with $p=x^2+ny^2$ and such that the conditions (i) and (ii) above are
satisfied. If these properties imply that $p$ is prime, then the number $n$ belongs to the set of so-called {\it numeri idonei}. Euler found 65 such integers $n$:
\par
{\narrower
  $1$, $  2$, $3$, $4$, $5$, $6$, $7$, $8$, $9$, $10$, $12$, $13$, $15$, $16$, $18$,
$21$, $22$, $24$, $25$, $28$, $30$, $33$, $37$, $40$, $42$, $45$, $48$, $57$,
$58$,
$60$,
$70$, $72$, $78$, $85$, $88$, $93$, $102$, $105$, $112$, $120$, $130$, $133$, $165$,
$168$, $177$, $190$, $210$, $232$, $240$, $253$, $273$, $280$, $312$, $330$, $345$,
$357$, $385$, $408$, $462$, $520$, $760$, $840$, $1320$, $1365$, $1848$.

\par}

\par

\noindent
It is known that there is at most one more number in the list, but one expects there
is no other one.

\medskip
Here is just one example ([\refSierpinskiA] problem 58, p.~112; [\refGuy] 
D18) of an open problem dealing with simultaneous Diophantine quadratic
equations: {\it  Is there a perfect integer cuboid?} The existence of a box
with integer edges $x_1,x_2,x_3$, integer face diagonals
$y_1,y_2,y_3$ and integer body diagonal $z$, amounts to solving the system of four
simultaneous Diophantine equations in seven unknowns
$$
\left\{
\matrix{\displaystyle
&x_1^2+x_2^2&=&y_3^2\cr
\displaystyle
&x_2^2+x_3^2&=&y_1^2\cr
\displaystyle
&x_3^2+x_1^2&=&y_2^2\cr
\displaystyle
&x_1^2+x_2^2+x_3^2&=&z^2\cr
}
\right.
$$
in $\bZ$. We don't know whether there is a solution, but it is known that
there is no perfect integer cuboid with the smallest edge $\le 2^{31}$.
 
\subsection{1.2. \subsectionunpointdeux}

In a Diophantine equation, the unknowns occur as the variables of
polynomials, while in an {\it  exponential Diophantine
equation} (see [\refShoreyTijdeman]), some exponents also
are variables. One may consider the above-mentioned Ramanujan-Nagell equation
$x^2+D=p^n$ as an exponential Diophantine equation.

 A famous problem which was open until 2002 is Catalan's one which dates back
to 1844 [\refCatalan], the same  year where Liouville constructed the first 
examples of transcendental numbers (see also [\refSierpinskiA] problem 77,
p.~116; [\refSierpinskiB] n$^\circ$~60, p.~42;
[\refShoreyTijdeman] Chap.~12; 
[\refNarkiewiczCPNT] Chap.~11;
[\refRibenboimCatalan]; [\refGuy] D9;
[\refRibenboimNumbersFriends] Chap.~7). The \og Note
extraite d'une lettre adress\'ee
\`a l'\'Editeur par Monsieur E.~Catalan, R\'ep\'etiteur \`a l'\'ecole
polytechnique de Paris\fg, published in Crelle Journal [\refCatalan], reads:

{\narrower\it 
{\og Je vous prie, Monsieur, de bien vouloir \'enoncer, dans votre recueil, 
le th\'eor\`eme suivant, que je crois vrai, bien que je n'aie pas encore
r\'eussi
\`a le d\'emontrer compl\`etement: d'autres seront peut-\^etre plus
heureux: 

Deux nombres entiers cons\'ecutifs, autres que $8$ et $9$, ne peuvent
\^etre  des puissances exactes; autrement dit: l'\'equation $x^m-y^n=1$,
dans laquelle les inconnues sont enti\`eres et positives, n'admet qu'une
seule solution.\fg
\par
 
}
\par
}

This means that  the only example of consecutive numbers which are perfect powers
$x^p$ with $p\ge 2$ should be $(8,9)$. Further information on the history of this question
is available in Ribenboim's book [\refRibenboimCatalan]. Tijdeman's
result [\refTijdemanCatalan] in 1976 shows that there are only
finitely many solutions. More precisely, for any
solution  $x,y,p,q$, the number $\max\{p,q\}$ can be bounded by an effectively
computable absolute constant. Once $\max\{p,q\}$ is bounded, only finitely
many exponential Diophantine equations remain to be considered, and there are
algorithms to complete the solution (based on Baker's method).
Such a bound has been computed, but
it is somewhat large:  M.~Mignotte proved that
any solution  $x,y,p,q$ to Catalan's equation should satisfy
$$
\max\{p,q\}<8\cdot 10^{16}.
$$

Catalan's claim was finally substantiated by P.~Mih\u{a}ilescu  [\refMihailescu] (see also [\refBilu] and [\refMetsankyla]).

\noindent

\proclaim Theorem {\numeroConjectureCatalan} { \rm (Catalan's Conjecture)}. 
The equation 
$$
x^p-y^q=1,
$$
where the unknowns $x$, $y$, $p$ and $q$ take  integer values  all $\ge 2$, has only one
solution, namely $(x,y,p,q)=(3,2,2,3)$.

The final solution by Mih\u{a}ilescu involves  deep
tools from the theory of cyclotomic fields. Initially  sharp measures of linear independence of logarithms of algebraic numbers were required
(namely a specific estimate for two logarithms due to M.~Laurent, M.~Mignotte and Yu.V.~Nesterenko was required), but then a solution using neither results from transcendental number theory nor the help of a computer was derived.

Catalan asked for integral solutions, like in Siegel's Theorem,
while Faltings's Theorem deals with rational points.  D.~Prasad suggested that the set of tuples  
$(x,y,p,q)$ in $\bQ^2\times\bN^2$  satisfying the conditions 
$$
\hbox{\it    $x^p-y^q=1$, and the
curve 
$X^p-Y^q=1$ has  genus $\ge 1$}
$$ 
should be finite -- evidence for this is provided by the $abc$  Conjecture
(see \S 2.1).

\smallskip

The fact that the right hand side in Catalan's equation is $1$ is 
crucial:  nothing is known if one replaces it by another positive
integer. The next conjecture was proposed by S.S.~Pillai [\refPillai] at a
conference of the Indian Mathematical Society held in Aligarh (see also
[\refSierpinskiA] problem 78, p.~117;
[\refShoreyTijdeman]; [\refTijdemanSurveyB]; [\refShoreySurvey]).

\proclaim Conjecture {\numeroConjecturePillai} { \rm  (Pillai)}.
Let $k$ be a positive integer.
The equation 
$$
x^p-y^q=k,
$$
where the unknowns $x$, $y$, $p$ and $q$ take integer values, all $\ge 2$, has only 
finitely many solutions $(x,y,p,q)$.

This means that in the increasing sequence of perfect powers $x^p$, with   $x\ge 2$
and $p\ge 2$:
$$
4,\, 8,\, 9,\, 16,\, 25,\, 27,\, 32,\, 36,\, 49,\, 64,\, 81,\, 100,\, 121,\, 
125,\, 128,\, 144,\, 169,\, \dots,
$$ 
the difference between two consecutive terms tends to infinity. It is not even 
known that for, say, $k=2$, Pillai's equation has only
finitely many solutions. A related open question is whether the number $6$ occurs as
a   difference between two perfect powers: {\it Is there a solution to the
Diophantine equation
$x^p-y^q=6$? } (see 
[\refSierpinskiB] problem 238a, p.~116).

A conjecture which implies Pillai's one has been suggested by T.N.~Shorey in
[\refShoreyConjectures]: this is the very problem which motivated C.L.~Siegel
in [\refSiegelPreuss]. Let
$f\in\bZ[X]$ be a polynomial of degree $n$ with at least two distinct roots and
$f(0)\not=0$. Let $L$ be the number of nonzero coefficients of $f$: write
$$
f(X)=b_1X^{n_1}+\cdots+b_{L-1}X^{n_{L-1}}+b_L
$$
with $n=n_1>n_2>\cdots>n_{L-1}>0$ and $b_i\not=0$ ($1\le i\le L$).
Set $H=\rmH(f)=\max_{1\le i\le L}|b_i|$.

\proclaim Conjecture {\numeroConjectureShorey} { \rm  (Shorey)}.
There exists a positive number $C$ which depends only on $L$ and $H$ with the
following property. Let $m$, $x$ and $y$ be rational integers with $m\ge 2$ and
$|y|>1$ satisfying
$$
y^m=f(x).
$$
Then either $m\le C$, or else there is a proper subsum in
$$
y^m-b_1x^{n_1}-\cdots-b_{L-1}x^{n_{L-1}}-b_L
$$
which vanishes.

An example with a vanishing proper subsum is
$$
y^m=x^{n_1}+x-2
$$
with $H=2$, $L=3$ and a solution $(m,x,y)=(n_1,2,2)$.

Consider now the positive integers which are  perfect powers $y^q$, with $q\ge
2$, and such that all digits in some basis $x\ge 2$ are $1$'s. Examples are
$121$ in basis $3$, $400$ in basis $7$ and $343$ in basis $18$. To find all
solutions  amounts to solving the exponential Diophantine equation  
$$
{x^n-1\over x-1}=y^q,
$$
where the unknown $x,y,n,q$ take positive rational integer values with $x\ge 2$, 
$y\ge 1$,
$n\ge 3$ and $q\ge 2$. Only $3$ solutions are known:
$$
(x,y,n,q)=\; 
(3,11,5,2),
\;
(7,20,4,2),\;
(18,7,3,3),
$$
corresponding to the above-mentioned three examples.
One does not know whether these are the only solutions  (see
[\refShoreyTijdeman]; [\refGuy] D10; 
 [\refTijdemanSurveyB]; [\refShoreySurvey],
 [\refBugeaudMignotte] and [\refShoreyConjectures]), but  it is expected
that there is no other one. 

The next question is to determine all  perfect powers with identical digits in some
basis, which amounts to solving the equation
$$
z{x^n-1\over x-1}=y^q,
$$
where the unknown $x,y,n,q,z$ take positive rational integer values with $x\ge 2$, 
$y\ge 1$,
$n\ge 3$, $1\le z<x$ and $q\ge 2$.

\medskip
Another type of exponential Diophantine equation has been studied in a joint
paper by H.P.~Schlicke\-wei and W.M.~Schmidt  [\refSchlickeweiViola]  where
they state the following conjecture. 

\proclaim Conjecture {\numeroConjectureSchlickeweiViola}.
 Let $k\ge 2$
be an integer and $\alpha_1,\ldots,\alpha_n$ be non-zero elements in a field
$K$ of zero characteristic, such that no quotient
$\alpha_i/\alpha_j$ with $j\not=i$ is a root of unity. Consider the
function
$$
F(X_1,\ldots,X_k)=
\det\pmatrix{
\alpha_1^{X_1}&\cdots&\alpha_k^{X_1}\cr
\vdots&\ddots&\vdots\cr
\alpha_1^{X_k}&\cdots&\alpha_k^{X_k}\cr
}.
$$
Then the equation
$$
F(0,x_2,\ldots,x_k)=0
$$
has only finitely many solutions $(x_2,\ldots,x_k)\in\bZ^{k-1}$ such that
in the corresponding determinant, all
$(k-1)\times k$ and all $k\times (k-1)$ submatrices have rank $k-1$.

Further exponential Diophantine equations are worth of study. See for instance
[\refNarkiewiczCPNT] Chap.~III,  [\refShoreyTijdeman],  
[\refTijdemanSurveyB] and [\refShoreySurvey]. 

Among the numerous applications of Baker's transcendence method are several
questions related to the greatest prime factors of certain numbers.  In
this connexion we mention Grimm's Conjecture ([\refGrimm], 
[\refNarkiewiczCPNT] Chap.~III
\S 3,  [\refGuy] B32):

\proclaim Conjecture {\numeroConjectureGrimm} { \rm  (Grimm)}.  Given $k$
consecutive composite integers $n+1,\ldots,n+k$, there exist
$k$ distinct primes
$p_1,\ldots,p_k$ such that $n+j$ is divisible by $p_j$, $1\le j\le k$.

This conjecture maybe rephrased as follows: {\it given an
increasing sequence of positive integers $n_1<\cdots<n_k$ for which the
product
$n_1\cdots n_k$ has fewer than $k$ distinct prime factors, there is a
prime $p$ in the range $n_1\le p\le n_k$. } The equivalence of this with the original formulation  follows from the ``marriage theorem''.

According to  P.~Erd\H{o}s and  J.L.~Selfridge, a consequence of
Conjecture {\numeroConjectureGrimm} is that between two consecutive squares there is always a prime number.


A weaker form of Conjecture {\numeroConjectureGrimm}, which is also an open problem, is:

\proclaim Conjecture {\numeroConjectureGrimmFaible}. If there is no prime in
the interval
$[n+1,n+k]$, then the product $(n+1)\cdots(n+k)$ has at least $k$ distinct prime
divisors.

M.~Langevin (personal communication) pointed out that Grimm's Conjecture cannot be extended 
to arithmetical progressions without a proviso: the numbers  $12$, $25$, $38$, $51$, $64$, $77$, $90$ belong to an arithmetic progression of ratio $13$, but the number of distinct prime factors of $12\cdot 25\cdot 64\cdot 90$ is only $3$. However he suggested in  [\refLangevin]  a
statement stronger than  Conjecture {\numeroConjectureGrimm}:
 
\proclaim Conjecture {\numeroConjectureGrimmLangevin}
{ \rm  (Langevin)}. Given an increasing sequence $n_1<n_2<\cdots<n_k$ of
positive integers such that $n_1,n_2,\ldots,n_k$
are multiplicatively dependent, there exists
a prime number in the interval $[n_1,n_k]$.

\medskip

Even if they may not be classified as Diophantine questions, the following 
open problems (see [\refLangGazette]) are related to this topic: the twin 
prime conjecture,  Goldbach Problem ({\it  is  every even integer
$\ge 4$ the sum of two primes?}), Bouniakovsky's conjecture, Schinzel's
hypothesis (H) (see also [\refSierpinskiA] \S 29) and the Bateman-Horn
Conjecture.
 
The Diophantine equation
$$
x^p+y^q=z^r
$$
has also a long history in relation with Fermat's last Theorem
([\refKraus], [\refRibenboimNumbersFriends] \S 9.2.D). If we look at the 
solutions in positive integers
$(x,y,z,p,q,r)$   for which
$$
{1\over p}+{1\over q}+{1\over r}<1
$$
and such that $x$, $y$, $z$ are relatively prime, then only $10$
solutions \footnote{($^1$)}{\petitsanserif
Up to obvious symmetries; in particular $1+2^3=3^2$
counts only for one solution.} 
are known:
$$ 
1+2^3=3^2,
\qquad
2^5+7^2=3^4,
\qquad 7^3+13^2=2^9,
\qquad 
2^7+17^3=71^2,
$$
$$
3^5+11^4=122^2,
\qquad
17^7+76271^3=21063928^2, \qquad 
1414^3+2213459^2=65^7,
$$
$$
9262^3+15312283^2=113^7,
\qquad 
43^8+96222^3=30042907^2,
\qquad 
33^8+1549034^2=15613^3.
$$
Since the condition
$$
{1\over p}+{1\over q}+{1\over r}<1
\quad\hbox{implies}\quad
{1\over p}+{1\over q}+{1\over r}\le{41\over 42}\virgule
$$
the $abc$ Conjecture (see
\S 2.1) predicts anyway that the set of such solutions is finite (the
``Fermat-Catalan'' Conjecture formulated by Darmon and Granville -- see [\refBeal]). For all
known solutions, one of
$p$,
$q$,
$r$ is
$2$; this led R.~Tijdeman and D.~Zagier to
conjecture\footnote{($^{2}$)}{\petitsanserif
This conjecture is also known as 
Beal's Conjecture -- see [\refBeal]} that there is no solution with the further
restriction that each of
$p$, $q$ and
$r$ is
$\ge 3$.  

A {\it Diophantine tuple } is a tuple
$(a_1,\ldots,a_n)$ of distinct positive integers such that $a_ia_j+1$ is a
square for $1\le i<j\le n$ (see [\refGuy] and [\refGyarmati]). Fermat  gave the example $(1,3,8,120)$, and Euler
showed that any Diophantine pair
$(a_1,a_2)$ can be extended to a Diophantine quadruple $(a_1,a_2,a_3,a_4)$. It
is not known whether there exists a Diophantine quintuple 
$(a_1,a_2,a_3,a_4,a_5)$, but A. Dujella [\refDujella] proved that each 
Diophantine quintuple has $\max\{a_1,a_2,a_3,a_4,a_5\} \le 10^{10^{26}}$. He
also proved that there is no Diophantine sextuple.
 
\subsection{1.3. \subsectionunpointtrois}

The original Markoff\footnote{($^{3}$)}{\petitsanserif  His name is
spelled {\it Markov} in probability theory.}
 equation (1879) is
$x^2 + y^2 + z^2 = 3 xyz$ (see [\refCassels] Chap.~II; 
[\refCusickFlahive] Chap.~2; [\refGuy] D12
 and [\refRibenboimNumbersFriends]
\S 10.5.B). Here is an algorithm which produces all solutions in positive
integers. Given any solution
$(x,y,z)=(m,m_1,m_2)$, we fix two of the three coordinates; then we get a
quadratic equation in the third coordinate, for which we already know a
solution. By the usual process of cutting with a rational line we deduce
another solution. From one solution $(m,m_1,m_2)$, this produces  three other
solutions 
$$
(m',m_1,m_2),\quad (m,m'_1,m_2),\quad (m,m_1,m'_2),
$$
where
$$
m'=3m_1m_2-m,\quad m'_1=3mm_2-m_1,\quad m'_2=3mm_1-m_2.
$$
These three solutions are called {\it  neighbors} of the original one.
Apart from the two so-called {\it  singular} solutions $(1,1,1)$ and
$(2,1,1)$, the three components of $(m,m_1,m_2)$ are pairwise distinct,
and the three neighbors of $(m,m_1,m_2)$ are pairwise distinct. Assuming $m>m_1>m_2$,
then one checks
$$
m'_2>m'_1>m>m'.
$$
Hence there is one neighbor of $(m,m_1,m_2)$ with maximum component less than
$m$, and two neighbors, namely $(m'_1,m,m_2)$ and
$(m'_2,m,m_1)$, with maximum component greater than $m$. It easily
follows that one produces all solutions, starting from $(1,1,1)$, by taking
successively the neighbors of the known solutions. Here is the Markoff tree, with the
notation of H.~Cohn [\refCohn], where   $(m'_1,m,m_2)$ is written on the right and
$(m'_2,m,m_1)$ on the left:

$$
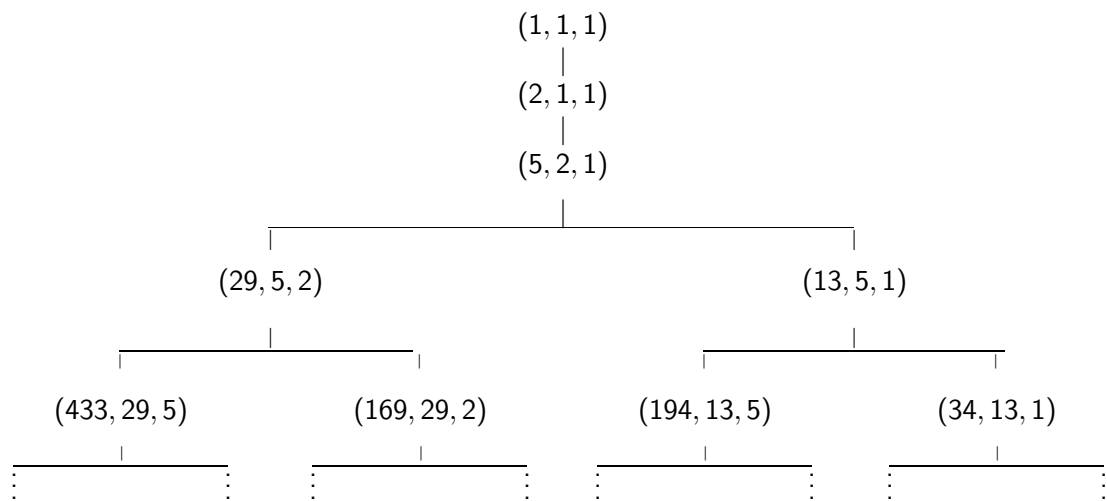
\qquad
\matrix{
        &        &     &(1,1,1)   &          &          &  \cr
        &        &     &|   &          &          &  \cr
        &        &     &(2,1,1)   &          &          &  \cr
        &        &     &|   &          &          &  \cr
        &        &     &(5,2,1)  &          &          &  \cr
 &    _{\vrule height .1 pt depth .1pt width 7.1  cm }^{\strut}\hskip -7.05  cm
    &     &   |      &          &          &  \cr
        &   ^|    &     &          &          &    ^|     &  \cr
        &(29,5,2)&     &          &          &(13,5,1)  & \cr 
  _{_{_{\vrule height .1 pt depth .1pt width 3.55  cm }}}^{\strut}\hskip -3.55  cm
&  
           _|    &     &  
   &_{_{_{\vrule height .1 pt depth .1pt width 3.65  cm }}}^{\strut}\hskip -3.65  cm    
&    _|     &  \cr ^{^|}        &  &      ^{^|}   &          & ^{^|} & & ^{^|}       
\cr 
(433,29,5)&  &(169,29,2)&          &(194,13,5)& & (34,13,1)\cr
\grandt
& 
 & \grandt
  &         
&    \grandt     &         & \grandt && \cr
}
$$

\bigskip\noindent
The main open
problem on this topic ([\refCassels] p.~33, 
[\refCusickFlahive] p.~11 and [\refGuy] D12) is to prove that each largest
component occurs only once in a triple of this tree: 

\proclaim Conjecture {\numeroconjectureMarkoff}.
Fix a positive integer $m$ for which 
the equation 
$$
m^2 + m_1^2 + m_2^2 = 3 mm_1m_2
$$ 
has a solution in positive integers $(m_1,m_2)$
with $0<m_1\le m_2\le m$. Then such a pair $(m_1,m_2)$ is unique.

This conjecture has been checked for $m\le 10^{105}$.

The sequence 
$$
1,\; 2,\;5,\;13,\;29,\;34,\;89,\;169,\;194,\;233,\;
433,\;
610,\;
985,\;
1325,\;1597,\dots
$$
of integers $m$ satisfying the hypotheses of 
Conjecture {\numeroconjectureMarkoff} is closely related to the question of
best   rational approximation to quadratic irrational real numbers: for each
$m$ in this  sequence, there is an explicit quadratic form $f_m(x,y)$ such that
$f_m(x,1)=0$ has a root $\alpha_m$ for which
$$
\limsup_{q\in\bZ,\; q\rightarrow\infty}
(q\Vert q\alpha_m\Vert)={m\over
\sqrt{9m^2-4}}\cdotp
\leqno{(\numeroapproximationMarkoff)}
$$
The sequence of $(m,f_m,\alpha_m,\mu_m)$ with $\mu_m=\sqrt{9m^2-4}/m$ starts 
as follows:
$$
\vbox{
\halign{
\sanserif
# &   \cc{ # }  & # &    \cc{ # } & # & \cc{ # }  & # &  \cc{ # }   &  # &   \cc{ # }
& # \cr
\th
\tv&  $m$ &\tv\tv& $1$&\tv&$2$  &\tv&  $5$&\tv&$13$&\tv\cr
\th
\th\tv&  $f_m(x,1)$ &\tv\tv& $x^2+x-1$&\tv&$x^2+2x-1$  &\tv&  $5x^2+11x-5$&\tv&$
13x^2+29x-13$&\tv\cr
\th
\tv&  $\alpha_m$ &\tv\tv& $\overline 1$&\tv&$\overline 2$  &\tv&  $\overline{2211}$&\tv&$ 
\overline{221111}$&\tv\cr
\th
\tv&  $\mu_m$ &\tv\tv& $\sqrt{5}$&\tv&$\sqrt{8}$  &\tv&  $\sqrt{221}/5$&\tv&
$\sqrt{1517}/13$&\tv\cr
\th
\cr
}}
$$
The third row gives the continued fraction expansion for $\alpha_m$, where
$\overline{2211}$, for instance, stands for $[2,2,1,1,2,2,1,1,2,2,1,1,\dots]$.
Conjecture {\numeroconjectureMarkoff} amounts to claiming that there is no
ambiguity in the notation
$f_m$: given $m$, two quadratic numbers $\alpha_m$ satisfying
(\numeroapproximationMarkoff) should be roots of equivalent quadratic forms.

Hence Markoff spectrum is closely related to rational approximation to a single
real  number. A generalization  to simultaneous approximation is considered in
 \S 2.2 below.

\goodbreak

\section{ 2. \sectiondeux}

In this section we restrict ourselves to problems in Diophantine approximation 
which do not require introducing a notion of height for algebraic
numbers: these will be discussed in
\S 4 only.

\subsection{2.1. \subsectiondeuxpointun}

For a positive integer $n$, we
denote by
$$
R(n)=\prod_{p|n}p
$$
the {\it  radical\/} or {\it  squarefree part\/} of $n$. 

The
$abc$ Conjecture was born from a discussion between D.W.~Masser
and J.~{\OE}sterl\'e  ([\refOesterle] p.~169; see also  [\refMasserabc], 
as well as [\refLangONCDI], [\refLangEncyclopaedie] Chap.~II
\S 1; 
[\refLangAlgebra] Ch.~IV  \S 7; 
[\refGuy] B19; [\refBrowkin]; 
[\refRibenboimNumbersFriends] \S 9.4.E; [\refVojta], [\refMazurPowersofNumbers] and [\refNitaj].

\proclaim Conjecture {\numeroconjectureabc} { \rm   ($abc$ Conjecture)}.
For each $\varepsilon>0$
there exists a positive number $\kappa(\varepsilon)$ which has the
following property: if $a$, $b$ and $c$ are three positive
rational integers which are relatively prime and satisfy $a+b=c$, 
 then
$$
c<\kappa(\varepsilon)R(abc)^{1+\varepsilon}.
$$

Conjecture {\numeroconjectureabc} implies a
previous conjecture by L.~Szpiro on the conductor of elliptic curves:
{\it   Given any $\varepsilon> 0$, there exists a constant $C(\varepsilon)>
0$ such that, for every elliptic curve with minimal discriminant
$\Delta$ and conductor $N$, we have $|\Delta| < C(\varepsilon) N^{6+\varepsilon}$. }

When $a$, $b$ and $c$ are three positive relatively prime integers
satisfying $a+b=c$, define 
$$
\lambda(a,b,c)= {\log c\over\log R(abc)} 
$$
and
$$
\varrho(a,b,c)= {\log abc\over\log R(abc)}\cdotp
$$

Here are the six largest known values for $\lambda(abc)$ (in
[\refBrowkin] p.~102--105  as well as in [\refNitaj], one can find all the $140$
known values of  $\lambda(a,b,c)$ which are $\ge 1.4$).
$$
\vbox{
\halign{
\sanserif
#&  \ #  & # \hfill& \hfill\  # &#& #\hfill &#&\ \ \hfill   # \hfill   &  #
& \quad  # \hfill& #\cr
\th
\tv&   &\tv\tv& $a+b$&=&$c$  &\tv&  $\lambda(a,b,c)$&\tv&author(s)&\tv\cr
\th\th
\tv& $1$ &\tv\tv& $2+3^{10}\cdot 109 $&$=$&$23^5$ &\tv&  $1.629912\dots$&\tv
&\'E. Reyssat&\tv\cr
\th
\tv& $2$ &\tv\tv& $11^2+3^2 5^6 7^3$&$=$&$2^{21}\cdot 23$ &\tv&  $1.625991\dots$&\tv
&B.~de Weger&\tv\cr
\th\tv& $3$ &\tv\tv& $19\cdot1307 + 7\cdot29^2\cdot31^8$&$=$&$2^8\cdot3^{22}\cdot5^4
$ &\tv&  $1.623490\dots$&\tv
&J.~Browkin -- J.~Brzezinski&\tv\cr
\th\tv& $4$ &\tv\tv& $283+5^{11}\cdot13^2$&$=$&$2^8\cdot3^8\cdot17^3
$ &\tv&  $1.580756\dots$&\tv
&
$\hbox{J.~Browkin -- J.~Brzezinski,}\atop\hbox{A.~Nitaj\hfill}$&\tv\cr
\th\tv& $5$ &\tv\tv& $1+2\cdot 3^7$&$=$&$5^4\cdot 7
$ &\tv&  $1.567887\dots$&\tv
&B.~de Weger&\tv\cr
\th\tv& $6$ &\tv\tv& $7^3+3^{10}$&$=$&$2^{11}\cdot 29
$ &\tv&  $1.547075\dots$&\tv
&B.~de Weger&\tv\cr
\th
\cr
}}
$$
Here are the six largest known values for $\varrho(abc)$, according to    
 [\refNitaj], where one can find the
complete list of $46$ known triples $(a,b,c)$ with $0<a<b<c$, $a+b=c$ and
$\gcd(a,b)=1$ satisfying 
$\varrho(a,b,c)> 4$.
$$
\vbox{
\halign{
\sanserif
#&  \ #  & # \hfill& \hfill\  # &#& #\hfill &#&\ \ \hfill   # \hfill   &  #
& \quad  # \hfill& #\cr
\th
\tv&   &\tv\tv& $a+b$&=&$c$  &\tv&  $\varrho(a,b,c)$&\tv&author(s)&\tv\cr
\th\th
\tv& $1$ &\tv\tv& $13 \cdot 19^6+ 	2^{30} \cdot 5$&$=$&$ 	3^{13} \cdot 11^2\cdot 31 $ 
&\tv&  
$4. 41901\dots$	&\tv
&A.~Nitaj&\tv\cr
\th
\tv& $2$ &\tv\tv& $
2^5 \cdot 11^2 \cdot 19^9 	+  
5^{15} \cdot 37^2 \cdot 47 
$&$=$&$
	3^7 \cdot 7^{11} \cdot
743 	$ &\tv&  $4.26801 \dots$&\tv
&A.~Nitaj&\tv\cr
\th\tv& $3$ &\tv\tv& $
2^{19} \cdot 13 \cdot 103	+ 7^{11}
$&$=$&$
	3^{11} \cdot 5^3 \cdot 11^2	
$ &\tv&  $4. 24789	\dots$&\tv
&B.~de Weger&\tv\cr
\th\tv& $4$ &\tv\tv& $ 
2^{35} \cdot 7^2 \cdot 17^2 \cdot 19 +	3^{27} \cdot 107^2 
$&$=$&$
5^{15} \cdot 37^2
\cdot 2311   
$ &\tv&  $	4.23069\dots$&\tv
&
A.~Nitaj&\tv\cr
\th\tv& $5$ &\tv\tv& $
	3^{18} \cdot 23 \cdot 2269 +	17^3 \cdot 29 \cdot 31^8 
$&$=$&$
	2^{10}
\cdot 5^2
\cdot 7^{15}$ &\tv&  $4.22979 \dots$&\tv
&A.~Nitaj&\tv\cr
\th\tv& $6$ &\tv\tv& $
17^4 \cdot 79^3 \cdot 211 +	2^{29} \cdot 23 \cdot
29^2
$&$=$&$
  5^{19}
$ &\tv&  $4.22960\dots$&\tv
&A.~Nitaj&\tv\cr
\th
\cr
}}
$$


As noticed by M.~Langevin [\refLangevinLuminy], a consequence of the $abc$
Conjecture 
is the solution of the following open problem
[\refErdosMonthly]:

\proclaim 
Conjecture {\numeroconjectureErdosWoods} { \rm   (Erd\H{o}s-Woods)}. 
There exists a positive integer $k$ such that, for $m$ and
$n$  positive integers, the conditions
$$
R(m+i)=R(n+i)\quad (i=0,\ldots,k-1)
$$
imply $m=n$.

Conjecture {\numeroconjectureErdosWoods} is motivated by the following
question raised by J.~Robinson: {\it Is   first
order arithmetic  definable using only the successor function $S:x\mapsto
x+1$   and   the coprimarity $x\perp y \Leftrightarrow (x,y)=1$?} It would
suffice to decide whether the function
$x\mapsto 5^x$ can be defined in the language $(S,\perp)$; see
[\refWoods], [\refGuy] B29 and B35, [\refBLSW].

\par From the $abc$ Conjecture  (or even the weaker version with ``some $\varepsilon<1$''  rather than ``all $\varepsilon>0$''), it follows that, apart from possibly
finitely many exceptions $(m,n)$, $k=3$ is an admissible value.
Indeed, assume $m>n$. Using the $abc$ Conjecture  with $a=m(m+2)$,
$b=1$, $c=(m+1)^2$, we get
$$
m^2\le \kappa(\varepsilon) R\bigl(m(m+1)(m+2)\bigr)^{1+\varepsilon}.
$$
Now if $R(m+i)=R(n+i)$ for $i=0,1,2$ then  $R(m+i)$ divides $m-n$, hence the number
$$ 
R\bigl(m(m+1)(m+2)\bigr)=\lcm\bigl(R(m), R(m+1),R(m+2)\bigr)
$$ 
divides $m-n$ and therefore $m^2\le
\kappa(\varepsilon)m^{1+\varepsilon}$. This shows that $m$ is bounded.

One 
suspects that there is no exception at all with $k=3$:
this would mean that if
$m$ and
$n$ have the same prime divisors, 
$m+1$, $n+1$ have the same prime divisors and $m+2$, $n+2$ have the same prime 
divisors, then $m=n$.

That $k=2$ is not an admissible value is easily seen:
 $75$ and $1215$ have the same prime divisors, and this
is true also for  $76$ and $1216$:
$$
R(75)=15=R(1215),\qquad
R(76)=2\cdot 19=R(1216).
$$
Apart from this sporadic example, there is also a sequence of examples: for
$m=2^h-2$ and 
$n=m(m+2)=2^hm$ we have
$$
R(m)=R(n)\and
R(m+1)=R(n+1)
$$
because $n+1=(m+1)^2$.

A generalization of the Erd\H{o}s-Woods Problem to arithmetic progressions has
been suggested by T.N.~Shorey: 

\smallskip

\item{(?)}{\it Does there exist a positive integer $k$ such
that, for any non zero integers
$m$, $n$, $d$ and $d'$  satisfying $\gcd(m,d)=\gcd(n,d')=1$,
the conditions
$$
R(m+id)=R(n+id')\quad (i=0,\ldots,k-1)
$$
imply $m=n$ and $d=d'$?}

\smallskip\noindent
On the one hand, if the answer is positive the integer $k$ is at least $4$, as shown by
several examples of quadruples
$(m,n,d,d')$, like  $(2,2,1,7)$, $(2,8,79,1)$ or $(4,8,23,1)$:
$$
R(2)=R(2),\quad R(3)=R(2+7),\quad R(4)=R(2+2\cdot 7),
$$
$$
R(2)=R(4)=R(8),\quad
R(2+79)=R(4+23)=R(9),\quad
R(2+2\cdot 79)=R(4+2\cdot 23)=R(10).
$$
On the other hand, under the $abc$ Conjecture, Shorey's question has a positive answer with $k= 5$ (see [\refLangevinRocky]).

\smallskip

Another related problem of T.S.~Motzkin and E.G.~Straus ([\refGuy] B19) is to 
determine the pairs of integers
$m,n$ such that $m$ and $n+1$ have the same prime divisors, and also $n$ and $m+1$ 
have the same set of prime divisors. The known examples are
$$
m=2^k+1,\qquad n=m^2-1\qquad (k\ge 0)
$$
and the sporadic example $m=35=5\cdot 7$, $n=4374=2\cdot 3^7$,  which gives
$m+1=2^2\cdot 3^2$ and $n+1=5^4\cdot 7$.
\medskip

 We also quote another related
conjecture attributed to P.~Erd\H{o}s in [\refLangevinLuminy] and to 
R.E.~Dressler in [\refNitaj].

\proclaim Conjecture {\numeroConjectureErdosDressler}  {\rm  
(Erd\H{o}s--Dressler)} . If $a$ and
$b$ are two positive integers with $a<b$ and
$R(a)=R(b)$  then there is a prime $p$ with $a< p< b$.

\medskip

The first estimates in the direction of
the $abc$ Conjecture  have been achieved by
C.L.~Stewart and R.~Tijdeman, and then refined by C.L.~Stewart and Yu Kunrui (see
[\refStewartYuKunruiI], [\refStewartYuKunruiII]), using ($p$-adic) lower bounds for linear forms in logarithms:
if $a$, $b$, $c$ are relatively prime positive integers with $a+b=c$, then
$$
\log c\le\kappa 
R^{1/3} (\log R )^3
$$
with $R=R(abc)$.

An explicit version has been worked out by Wong Chi Ho in 1999 [\refWong], following an earlier version of [\refStewartYuKunruiII]:
for $c>2$ the estimate
$$
\log c\le  R^{(1/3)+(15/\log\log R)}
$$
holds.

Further connexions between the $abc$ Conjecture 
and measures of linear  independence of logarithms of algebraic numbers
have been pointed out by A.~Baker [\refBakerEger] and P.~Philippon
[\refPhilipponAustralie] (see also [\refmiwGL] exercise 1.11). We
reproduce here the main conjecture of the addendum of 
[\refPhilipponAustralie]. For a rational number $a/b$ with relatively
prime integers $a,b$, we denote by $\rmh(a/b)$ the number
$\log\max\{|a|,|b|\}$.

\proclaim Conjecture {\numeroConjecturePhilipppon} { \rm  
(Philippon)}.
There exist real numbers $\varepsilon$, $\alpha$ and $\beta$ with
$0<\varepsilon<1/2$, $\alpha\ge 1$ and $\beta\ge 0$, and a positive integer
$B$, such that for any nonzero rational numbers $x$, $y$ satisfying
$xy^B\not=1$, the following is true: if $S$ denotes the set of prime numbers for which 
$|xy^B+1|_p<1$, then
$$
-\sum_{p\in S}\log|xy^B+1|_p\le
B\Bigl(\alpha \rmh(x)+\varepsilon\rmh(y)+
\bigl(\alpha B+\varepsilon\bigr)\bigl(\beta+\sum_{p\in S}\log p\bigr)\Bigr).
$$

The conclusion is a lower bound for the $p$-adic distance between $-xy^B$
and $1$; the main point is that several $p$'s are involved. The conjecture  {\numeroConjecturePhilipppon} is telling us something about the prime decomposition of all numbers $xy^B+1$ for some fixed but unspecified value of $B$ -- and it implies the $abc$ Conjecture.

Examples of optimistic archimedean estimates related to measures of
linear independence of logarithms of algebraic numbers are the
Lang-Waldschmidt Conjectures in [\refLangECDA] (introduction to Chap.~X and
XI, p.~212--217). Here is a simple example:

\proclaim Conjecture {\numeroConjectureLangW} { \rm  
(Lang-Waldschmidt)}.
For any $\varepsilon>0$, there 
exists a constant $C(\varepsilon)>0$ such that, for any 
nonzero rational integers $a_1,\ldots,a_m$, $b_1,\ldots,b_m$
with
$a_1^{b_1}\cdots a_m^{b_m}\not=1$,
$$
\left|a_1^{b_1}\cdots a_m^{b_m}-1\right|\ge
{C(\varepsilon)^m B\over (|b_1|\cdots |b_m|\cdot |a_1|\cdots
|a_m|)^{1+\varepsilon}}\virgule
$$
where 
$B=\max_{1\le i\le m}|b_i|$.

Similar questions related to Diophantine approximation on tori are
discussed in [\refLangEncyclopaedie] Chap.~IX, \S 7.

 Conjecture {\numeroConjectureLangW} deals with rational
integers; we shall consider more generally algebraic numbers in \S 4, once
we have defined a notion of height in \S 3.

 
A very sharp conjectured lower bound for infinitely many elements in a specific  sequence
$$
\left|e^{b_0}a_1^{b_1}\cdots a_m^{b_m}-1\right|\
$$
with $b_0$ arbitrary, and where all the exponents $b_i$  have the same sign (compare with Conjecture \numeroProblemeMahler)  is shown by
J.~Sondow in [\refSondow] to yield the irrationality of Euler's
constant.


\par From either the $abc$ Conjecture  or Conjecture
{\numeroConjectureLangW} one deduces a quantitative refinement to Pillai's 
Conjecture
\numeroConjecturePillai:

\proclaim Conjecture {\ConjecturePillaiQuantitative}. For any
$\varepsilon>0$,  there is a constant $C(\varepsilon)>0$ such that, for any positive
integers  $x$,
$y$,
$p$, $q$ satisfying 
$x^p\not= y^q$, the inequality
$$
|x^p-y^q| \ge C(\varepsilon)\max\{x^p, y^q\}^{1-(1/p)-(1/q)-\varepsilon}
$$
holds. 

\smallskip

We consider two special cases of Conjecture \ConjecturePillaiQuantitative:
first $(p,q)=(2,3)$, which gives rise to  Hall's Conjecture [\refHall]
(also [\refLangEncyclopaedie]  Chap.~II, \S 1): 

\proclaim Conjecture {\numeroConjectureHall} { \rm  (Hall)}. If $x$
and $y$ are positive integers with $y^2\not=x^3$, then
$$
|y^2-x^3|\ge C\max\{y^2,x^3\}^{1/6}.
$$

In this statement there is no $\varepsilon$. On the one hand, maybe
Conjecture {\numeroConjectureHall} is true by a sort of accident, but one
may also expect that the estimate is too strong to be true. On the other hand, with the exponent $(1/6)-\varepsilon$, the $abc$ Conjecture provides a lower bound not only for $|y^2-x^3|$, but also for its radical [\refLangevinRocky]: {\it for $x$ and $y$ relatively prime positive integers with $y^2\not=x^3$, 
$$
R\bigl(|y^2-x^3|\bigr)\ge C(\varepsilon)\max\{y^2,x^3\}^{(1/6)-\varepsilon}.
$$
}

The exponent $1/6$ in  Conjecture {\numeroConjectureHall} is optimal, as shown by L.V.~Danilov and A.~Schinzel: using the polynomial identity
$$
(X^2-6X+4)^3-(X^2+1)(X^2-9X+19)^2=27(2X-11)
$$
(which is related to Klein's identity for the icosahedron -- cf. [\refLangevinJNTBx], Th.~6), they show that there exist infinitely many  pairs of positive integers $(x,y)$ such that
$$
0<|y^2-x^3| < {54\over 25\sqrt{5}}\cdot \sqrt{x}.
$$

The smallest known value for
$|y^2-x^3|/\sqrt{x}$ (N.~Elkies, 1998) is
$0.0214\dots$, with
$$
x=3\cdot 7\,211\cdot
38\,791\cdot
6\,975\,841,\quad
y=2\cdot
3^2\cdot
15\,228\,748\,819\cdot
1\,633\,915\,978\,229,
$$
$$
x^3-y^2=3^3\cdot 7^2\cdot 17\cdot73.
$$

\medskip
The second special case of Conjecture \ConjecturePillaiQuantitative\ we  consider is $(x,y)=(3,2)$. The question of how small 
$3^n-2^m$ can be in comparison with $2^m$ has been raised by
J.E.~Littlewood [\refGuy] F23. The example   
$$
{3^{12}\over 2^{19}}=1+{7153\over 524288}=1.013\dots
$$
is related to music scales. 

\smallskip

For further questions dealing with exponential Diophantine equations, we
refer to Chap.~12 of the book of T.N.~Shorey and R.~Tijdeman [\refShoreyTijdeman],
as well as to the more recent surveys [\refTijdemanSurveyB] and
[\refShoreySurvey].

\subsection{2.2. \subsectiondeuxpointdeux}

One of the main open problems in Diophantine approximation is to
produce an effective  version of the Thue-Siegel-Roth  
Theorem: {\it  For any $\varepsilon>0$ and any irrational algebraic number
$\alpha$, there is a positive constant $C(\alpha,\varepsilon)>0$ such that, for
any rational number
$p/q$,
$$
\left|\alpha-{p\over q}\right|>{C(\alpha,\varepsilon) \over q^{2+\varepsilon}}
\cdotp
\leqno{(\numeroThueSiegelRoth)}
$$
}
In connexion with the negative answer to Hilbert's 10th Problem by
Yu.~Matiyasevich, it has been suggested by M.~Mignotte that an effective
version of Schmidt's Subspace Theorem (which extends the Thue-Siegel-Roth
Theorem to simultaneous Diophantine approximation) may be impossible. If this
turns out to be the case also for the special case of the Thue-Siegel-Roth
Theorem itself, then, according to E.~Bombieri (see [\refNitaj]), an effective version
of the
$abc$ Conjecture would also be out of reach. M.~Langevin noticed that the
$abc$ Conjecture yields a stronger inequality than Roth's one:
$$
\left|\alpha-{p\over q}\right|> {C(\varepsilon)\over R(pq)q^{\varepsilon}}\cdotp
$$

So far, effective improvements are known only for Liouville's  bound, and
to improve them is already a big challenge.

\smallskip

Another goal would be to improve the estimate in  Roth's Theorem: in the
lower bound
 (\numeroThueSiegelRoth)
one would like to  replace $q^{-2-\varepsilon}$ by, say, $q^{-2}(\log
q)^{-1-\varepsilon}$. It is expected that for  any irrational  real algebraic 
number $\alpha$ of degree $\ge 3$, in the inequality 
 (\numeroThueSiegelRoth) the term
$q^{-2-\varepsilon}$  cannot be
replaced by
$q^{-2}$, but the set of $\alpha$ for which the answer is known is empty! This
question is often asked for the special case of the number
$\root {3} \of 2$, but another interesting example (due to  Stanislaw Ulam -- see for
instance [\refGuy] F22) is the real algebraic number $\xi$ defined by
$$
\xi={1\over\xi+y}
\with
y={1\over 1+y}\cdotp
$$
Essentially nothing is known on the continued fraction expansion of a real
algebraic number of degree
$\ge 3$; one does not know the answer to any of the following two questions:

\smallskip\noindent
{(\numeroQuestionBPQ~?)} {\it  Does there exist a real algebraic number of degree $\ge 3$ with bounded partial
quotients? }

\smallskip\noindent
{(\numeroQuestionUPQ~?)} {\it  Does there exist a real algebraic number of degree $\ge 3$ with unbounded
partial quotients? }

\smallskip\noindent
It is usually expected is that the continued fraction expansion of a 
real algebraic number of degree at least $3$ has always unbounded partial
quotients.
More precisely one expects
that real algebraic numbers of degree
$\ge 3$ behave like ``almost all'' real numbers (see \S 5.1). 

Let
$\psi(q)$ be a continuous positive real valued function. Assume that the
function
$q\psi(q)$ is non-increasing. Consider the inequality
$$
\left|\theta-{p\over q}\right|> {\psi(q)\over q}\cdotp
\leqno{(\numeroinegalitepourRothRaffine)}
$$

\proclaim Conjecture {\numeroRothRaffine}. 
Let $\theta$ be real algebraic number of degree at least $3$.
Then inequality (\numeroinegalitepourRothRaffine) has infinitely
many solutions in integers $p$ and $q$ with $q>0$
if and only if the integral
$$
\int_1^\infty \psi(x) dx
$$
diverges.  

\medskip

A far reaching generalization of the Thue-Siegel-Roth Theorem to simultaneous
approximation is the Schmidt Subspace Theorem.  Here are two special cases:

\smallskip
\item{$\bullet$}{\it  Given real algebraic numbers $\alpha_1,\ldots,\alpha_n$
such that
$1,\alpha_1,\ldots,\alpha_n$ are linearly independent over $\bQ$, for any
$\varepsilon>0$  the inequality
$$
\max_{1\le i\le n} 
\left|\alpha_i-{p_i\over q}\right|<{1\over
q^{1+(1/n)+\varepsilon}}
$$
has only finitely many solutions $(p_1,\ldots,p_n,q)$  in $\bZ^{n+1}$ with $q>0$.
}
\smallskip
\item{$\bullet$}{\it  Given real algebraic numbers $\alpha_1,\ldots,\alpha_n$
such that
$1,\alpha_1,\ldots,\alpha_n$ are linearly independent over $\bQ$, for any
$\varepsilon>0$   the inequality
$$
\left|q_1\alpha_1+\cdots+q_n\alpha_n-p\right|<{1\over
q^{n +\varepsilon}}
$$
has only finitely many solutions $(q_1,\ldots,q_n,p)$ in $\bZ^{n+1}$ with
$q=\max\{|q_1|,\ldots,|q_n|\}>0$.}

\medskip\noindent
These  two types of Diophantine statements are parallel to the two types 
of Pad\'e Approximants. It would be interesting to consider the analog of 
Schmidt's Subspace Theorem in case of Pad\'e Approximants, and also to
investigate a corresponding analogue of  Khinchine's transference principle
[\refCassels].

\medskip

One of the most important consequences of  Schmidt's
Subspace Theorem is the finiteness of nondegenerate solutions of the equation
$$
x_1+\cdots+x_n=1,
$$
where the unknowns take integer values (or $S$-integer values) in a number field. Here,
non-degenerate means that no proper subsum vanishes. One main open question
is to  prove an effective version of this result. Schmidt's Theorem, which
is a generalization of Roth's Theorem, is not effective. Only for $n=2$ does one
know bounds for the solutions of the $S$-unit equation $x_1+x_2=1$, thanks
to Baker's method (see [\refBakerLivre] Chap.~5;
[\refLangECDA] Chap.~VI;
 [\refShoreyTijdeman]
Chap~1; [\refSerre]
  and [\refLangEncyclopaedie]). One would like to extend Baker's method (or
any other effective method) to the higher dimensional case.

\medskip
A generalization of the Markoff spectrum to simultaneous approximation is not yet
available: even the first step is missing. Given a positive integer $n$ and real
numbers
$(\xi_1,\ldots,\xi_n)$, not all of which are rational,     define 
$c_n=c_n(\xi_1,\ldots,\xi_n)$ to be the infimum of all
$c$  in the range $0<c\le 1$ for which the inequality
$$
q|q\xi_i-p_i|^n<c
$$
has infinitely many solutions. Then define the {\it  $n$-dimensional
simultaneous Diophantine approximation constant} $\gamma_n$   to be the supremum
of $c_n$ over tuples  $(\xi_1,\ldots,\xi_n)$ as above. 
Following [\refFinch], here is a
summary of what is known about the first values of the approximation constants: 
$$\matrix{
&&&&\gamma_1&=&\displaystyle{
1\over\sqrt{5}}&=&0.447\dots& \hbox{(Hurwitz)}\cr\cr
0.285\dots&=&\displaystyle{2\over 7}&\le& 
\gamma_2 &\le&\displaystyle{64\over169}
&=&0.378 \dots&\hbox{(Cassels and
Nowak)}\cr\cr
0.120\dots&=&\displaystyle{2\over5\sqrt{11}}&\le& 
\gamma_3 &\le&\displaystyle{1\over 2(\pi-2)}
&=&0.437 \dots&\hbox{(Cusick and
Spohn)}\cr
}
$$
The question remains open as to whether there are pairs with approximation constant larger than $2/7$ (see [\refBriggs]).

 \bigskip
We illustrate now with Waring's Problem the importance of proving
effective Roth-type inequalities for irrational algebraic numbers.

 In 1770, a few months before 
J.L.~Lagrange proved that every positive integer is the sum of at
most four squares of integers, E.~Waring ([\refWaring] Chap.~5,  Theorem
47 (9))  wrote:

{\narrower

{\it 
``Every  integer is a cube or the sum of two, three, \dots nine cubes; every
integer is also the square of a square, or the sum of up to nineteen such; and
so forth. Similar laws may be affirmed for the correspondingly defined numbers
of quantities of any like degree.''}

\par
}

\par
\noindent
See also Note 15 of the translator in [\refWaring].

\medskip
For $k\ge 2$
define $g(k)$ as the smallest positive integer $g$ such that any integer is the
sum of $g$ elements of the form
$x^k$ with $x\ge 0$. In other terms, for each positive integer $n$ the equation
$$
n=x_1^k+\cdots+x_m^k
$$
has  a solution if $m=g(k)$, while there is a $n$ which is not  the
sum of 
$g(k)-1$ such $k$-th powers. 

Lagrange's Theorem, which solved a conjecture of Bachet and Fermat, is
$g(2)=4$. Following   Chap.~IV of [\refNarkiewiczCPNT],  here are the values of 
$g(k)$ for the first integers
$k$, with the  name(s) of the author(s) and the date:
$$
\vbox{
\halign{
 #&  \hfill   #\hfill   &#& \hfill  # \hfill  &#&\hfill     # \hfill     &  #
& \hfill   # \hfill  & # & \hfill   # \hfill & # & \hfill   # \hfill  & #\cr
\th
\tv&\ \ \hskip -.5 true cm \hfill $g(2)=4$ \hfill $\!$&\tv& \hfill $g(3)=9$   \hfill&\tv&   \hfill $g(4)=19$  \hfill &\tv& \hfill \ \ $g(5)=37$ \ \  \hfill &\tv& \hfill\  \ $ g(6)=73$\ \  \hfill&\tv&
 \hfill\ $ g(7)=143$\   \hfill&\tv\cr
\th
\tv&     
\ J.L.~Lagrange &\tv&\ A.~Wieferich  &\tv&  
$\matrix{
\hbox{\ R.~Balasubramanian}^{\strut}\cr
\hbox{J-M.~Deshouillers}\cr
\hbox{F.~Dress}_{\strut}\cr}$&\tv&\ J.~Chen&\tv&\
S.S.~Pillai &\tv &\ L.E.~Dickson &\tv\cr\th
\tv& $1770$ &\tv& $1909$  &\tv&  $1986$ &\tv& $1964$&\tv& $1940$ &\tv
& $1936$ &\tv\cr 
\th
\cr
}}
$$
For each integer $k\ge 2$, define 
$$
I(k)=2^k+[(3/2)^k]-2.
$$
It is easy to check $g(k)\ge I(k)$: write 
$$
3^k=2^kq+r \with 
0<r<2^k,\quad q=[(3/2)^k],
$$
and  consider the integer
$$
N=2^kq-1=(q-1)2^k+ (2^k-1)1^k .
$$
Since $N<3^k$, writing $N$ as a sum of $k$-th powers can involve no term
$3^k$, and since $N<2^kq$, it involves at most $(q-1)$ terms $2^k$, all
others being $1^k$; hence it requires a total number of at least
$(q-1)+ (2^k-1)=I(k)$ terms. 

 L.E.~Dickson and S.S.~Pillai (see for instance [\refHardyWright] or [\refNarkiewiczCPNT] 
 Chap.~IV) proved independently in 1936 that $g(k)=I(k)$ provided that the number $r=3^k-2^kq$ satisfies
$$
r\le 2^k-q-2.
$$
Otherwise there is another formula for $g(k)$.

It has been checked that the condition $r\le 2^k-q-2$ is satisfied for $3\le k\le 471~600~000$, and
K.~Mahler proved that it is also true for any sufficiently large $k$.  Hence $g(k)=I(k)$  for these values of $k$. The
problem is that Mahler's proof relies on a $p$-adic version of the
Thue-Siegel-Roth Theorem, and therefore is not effective. So there is a
gap, of which we don't even know the size. The conjecture, dating back to
1853, is $g(k)=I(k)$ for any $k\ge 2$, and this is true as soon as
$$
\left\Vert \left( {
3\over 2}\right)^k\right\Vert
\ge 
 \left({3\over 4}\right)^k,
$$
where $\Vert\cdot\Vert$ denote the distance to the nearest integer.
As remarked by S.~David, such an estimate (for sufficiently
large $k$) follows not only from Mahler's estimate, but also from the $abc$ Conjecture!

\medskip
In [\refMahlerZnumbers] K.~Mahler defined a {\it  $Z$-number} as a nonzero real
number
$\alpha$ such that the fractional part $r_n$ of $\alpha (3/2)^n$ satisfies
$0\le r_n<1/2$ for any positive integer $n$. It is not known whether
$Z$-numbers exist  (see [\refFlattoLagariasPollington]). A related remark
by J.E.~Littlewood ([\refGuy] E18) is that we are not yet able to prove
that the fractional part of $e^n$ does not tend to
$0$ as $n$ tends to infinity (see also Conjecture \numeroProblemeMahler\ below).

\medskip
A well known conjecture of Littlewood ([\refBakerLivre] Chap.~10, \S 1 and
[\refPollingtonVelani])
asserts that {\it for any pair $(x,y)$ of real numbers and any
$\varepsilon>0$, there exists a positive integer $q$ such that
$$
q\Vert qx\Vert\cdot \Vert qy\Vert<\varepsilon.
$$}
According to G.~Margulis (communication of G.~Lachaud), the
proofs
in a 1988 paper by B.F.~Skubenko (see { \tt 
M.R{.}$\,$94d$:$11047}) are not correct and cannot be fixed.

\medskip

There are several open questions known as ``view obstruction Problems''. 
One of them is the following. {\it Given $n$ positive integers
$k_1,\ldots,k_n$, there exists a real number $x$ such that
$$
\Vert k_ix\Vert\ge {1\over n+1}\for 1\le i\le n.
$$}
It is known that $1/(n+1)$ cannot be replaced by a larger number
[\refCusickPomerance].

\subsection{2.3. \subsectiondeuxpointtrois}

Given a real number $\theta$, the first Diophantine question is to decide
whether $\theta$ is rational or not. This is a qualitative question, and it
is remarkable that an answer is provided by a   quantitative property  of
$\theta$: it depends ultimately on the quality of rational Diophantine
approximations to $\theta$. Indeed, on the one hand,  if $\theta$ is rational,
then there exists a positive constant $c=c(\theta)$ such that
$$
\left|\theta -{p\over q}\right|>{c\over q}
$$
for any $p/q\in\bQ$.
An admissible value for $c$ is $1/b$ when $\theta=a/b$. On the other hand, if
$\theta$ is irrational, then there are infinitely many rational numbers
$p/q$ such that
$$
0<\left|\theta -{p\over q}\right|<{1\over q^2}\cdotp
$$
Hence, in order to prove that $\theta$ is irrational, it suffices to prove
that
 for any $\varepsilon>0$
there is a rational number
$p/q$ such that
$$
0<\left|\theta -{p\over q}\right|<{\varepsilon\over q}\cdotp
$$
This is a rather weak requirement: there are rational approximations in $1/q^2$,
and we need only to produce rational approximations better than the trivial
ones in $c/q$. Accordingly one should expect that it is rather easy to
prove the irrationality of a given real number. In spite of that, the class of
``interesting'' real numbers which are known to be irrational is not as large as
one would expect 
[\refKontsevichZagier]. For instance no proof of irrationality has been given so far
for Euler's constant 
$$
\gamma=\lim_{n\rightarrow\infty} \left(1+{1\over 2}+{1\over 3}+\cdots
+{1\over n}-\log n\right)=0.577215\dots,
$$ 
nor for Catalan's constant
$$ 
G=\sum_{n\ge 0}{(-1)^n\over (2n+1)^2}
= 
0.915965\dots
\virgule
 \leqno{(\numeroCatalan)}
$$
nor for 
$$
\Gamma(1/5)=\int_0^\infty e^{-t}t^{-4/5}dt=4.590843\dots
$$
or for numbers like
$$
e+\pi= 5.859874\dots,\quad
e^\gamma=1.781072\dots,\quad \zeta(5)=1.036927\dots ,\quad 
\zeta(3)/\pi^3=
0.038768\dots
$$
and
$$
\sum_{n\ge 1}{\sigma_k(n)\over n!} \quad (k=1,2)
\where
\sigma_k(n)=\sum_{d\mid n}d^k
$$
(see [\refGuy] B14). 

Here is another  irrationality  question raised by P.~Erd\H{o}s and
E.~Straus in 1975 (see [\refErdosDurham] and [\refGuy] E24). Define an
{\it  irrationality sequence} as an increasing sequence
$(n_k)_{k\ge 1}$ of positive integers such that, for any sequence
$(t_k)_{k\ge 1}$ of positive integers, the real number
$$
\sum_{k\ge 1}{1\over n_kt_k}
$$
is irrational. On the one hand, it has been proved by Erd\H{o}s that $(2^{2^k})_{k\ge 1}$
is an irrationality sequence. On the other hand, the sequence $(k!)_{k\ge 1}$
is not, since
$$
\sum_{k\ge 1}{1\over k!(k+2)}={1\over 2}\cdotp
$$
An open question is whether an irrationality sequence must increase
very rapidly. No irrationality sequence $(n_k)_{k\ge 1}$ is known for which 
$n_k^{1/2^k}$ tends to $1$ as $k$ tends to infinity.

Many further open irrationality
questions are raised in [\refErdosDurham]. Another related example is
Conjecture {\numeroConjectureLoxtonVdP} below.

\medskip
Assume now that the first step has been completed and that we know
our number $\theta$ is irrational. Then there are (at least) two directions
for further investigation: 

\smallskip
\item{(1)} Considering several real numbers $\theta_1,\ldots,\theta_n$, a
fundamental question is to decide whether or not they are linearly independent
over
$\bQ$.  One main example is to start with the successive powers of one number:
$1,\theta,\theta^2,\ldots,\theta^{n-1}$; the goal is to decide
whether $\theta$ is algebraic of degree $<n$. If $n$ is not fixed,  the
question is whether
$\theta$ is transcendental. This question, which is relevant also for
complex numbers, will be considered in the next section. Notice also that the
problem of algebraic independence is included here: it amounts to the linear
independence of monomials.

\smallskip
\item{(2)} Another direction of research is to consider a
quantitative refinement of the irrationality statement, namely an {\it  
irrationality measure}: we wish to bound from below the nonzero number
$|\theta-(p/q)|$ when
$p/q$ is any rational number; this lower bound will depend on $\theta$ as well
as the denominator $q$ of the rational approximation.
In case where  a statement weaker  than an irrationality result is known,
namely if one can prove only that at least one of
$n$ numbers
$\theta_1,\ldots,\theta_n$ is irrational, then a quantitative refinement will be
a lower bound (in terms of
$q$) for 
$$
\max\left\{\left|\theta_1-{p_1\over q}\right|,\ldots,
\left|\theta_n-{p_n\over q}\right|\right\},
$$
when $p_1/q, \ldots,p_n/q$ are $n$ rational numbers and $q>0$ a common
denominator. 

\smallskip\noindent
On the one hand, the study of rational approximation of real numbers is achieved in a
satisfactory way for numbers whose ``regular'' \footnote{($^4$)}{
{\petitsanserif
A ``regular'' continued fraction expansion 
$$
\left[a_0+{1\over a_1+}\; {1\over a_2+}\;   \cdots\;
\right]
$$
is written}  $[a_0 , a_1 , a_2 , \ldots]$.  {\petitsanserif A   continued fraction expansion of the form
$$
\left[a_0+{b_1\over a_1+}\; {b_2\over a_2+}\;   \cdots\;
\right]
$$
is called  ``irregular''.}} continued fraction expansion
is known. This is the case for rational numbers (!), for quadratic
numbers, as well as for a small set of transcendental numbers, like 
$$
e=[2,1,2,1,1,4,1,1,6,1,1,\dots]=[2,\{\overline{1,2m,1}\}_{m\ge 1}]
$$
$$
e^2=[7,2,1,1,3,18,5,1,1,6,30,8,1,1,9,42,11,\dots]=
[7,\{\overline{3m-1,1,1,3m,12m+6}\}_{m\ge
1}]
$$
and
$$
e^{1/n}=[1,n-1,1,1,3n-1,1,1,5n-1,1,1,\dots]=
[\{\overline{1,(2m-1)n-1,1}\}_{m\ge
1}]
$$
for $n>1$.  On
the other hand, even for a real number $x$ for which an ``irregular''
continued fraction expansion is known, like
$$
\log 2=\left[{1\over 1+}\; {1\over 1+}\; {4\over 1+}\;
{9\over 1+}\; \cdots\; {n^2\over 1+} \cdots\;
\right]
$$
or
$$
{\pi\over 4}=
\left[{1\over 1+}\; {9\over 2+}\; {25\over 2+}\;
{49\over 2+}\; \cdots \; {(2n+1)^2\over 2+}\;\cdots\;
\right]
$$
one does not know how well $x$ can be approximated by rational numbers. No regular
pattern has been observed or is expected  from the regular continued fraction of 
$\pi$:
$$
\pi=[3,7,15,1,292,1,1,1,2,1,3,1,14,2,1,1,
2,2,2,2,1,84,2,1,1,15,3,13,1,4,2,6,6,6,1,\dots],
$$
nor from any number ''easily'' related to $\pi$.

One expects
that for any
$\varepsilon>0$ there are constants
$C(\varepsilon)>0$ and $C'(\varepsilon)>0$ such that
$$
\left|\log 2-{p\over q}\right|>{C(\varepsilon)\over q^{2+\varepsilon}}
\and
\left|\pi-{p\over q}\right|>{C'(\varepsilon)\over q^{2+\varepsilon}}
$$
hold for any $p/q\in\bQ$, but this is known only with larger exponents,
namely $3.8913\dots$ and $8,0161\dots$ respectively
(Rukhadze and Hata). The sharpest known exponent for an irrationality measure of 
$$
\zeta(3)=\sum_{n\ge 1}{1\over n^3}=1.202056\dots
$$
is $5.513891\dots$, while for $\pi^2$ (or for $\zeta(2)=\pi^2/6$) it is
$5.441243\dots $ (both results due to Rhin and Viola). For a number like
$\Gamma(1/4)$, the existence of absolute positive constants
$C$ and $\kappa$ for which
$$
\left|\Gamma(1/4)-{p\over q}\right|>{C \over q^{\kappa}}
$$
has been proved only recently [\refPhilipponGammaunquart]. The similar problem
for
$e^\pi$ is not yet solved: there is no proof so far that $e^\pi$ is not a
Liouville number.

Earlier we distinguished two directions for research once we know the
irrationality of some given numbers: either, on the qualitative side,  one
studies the linear dependence relations, or else, on the quantitative side, one
investigates the quality of rational approximation. One can combine both:  a
quantitative version of a result of $\bQ$-linear independence of $n$ real numbers
$\theta_1,\ldots,\theta_n$,  is a lower bound, in
terms of $\max\{|p_1|,\ldots,|p_n|\}$, for 
$$
\bigl|p_1\theta_1+\cdots+p_n\theta_n\bigr|
$$
when $(p_1,\ldots,p_n)$ is in $\bZ^n\setminus\{0\}$. 

For some specific classes of transcendental numbers,
A.I. Galochkin [\refGalochkin], A.N.~Korobov (Th.~1.22
of [\refFeldmanNesterenko] Chap.~1
\S 7) and more recently P.~Ivankov  proved extremely sharp
measures of linear independence (see [\refFeldmanNesterenko] Chap.~2
\S 6.2 and \S 6.3). 

A  general and important problem is to improve the known
measures of linear independence for logarithms of algebraic numbers, as
well as elliptic logarithms, abelian logarithms, and more generally  
logarithms of algebraic points on commutative algebraic groups.
For instance the conjecture that $e^\pi$ is not a Liouville number should
follow from improvements of known linear independence measures for logarithms of
algebraic numbers.

The next step, which is to get sharp measures of algebraic
independence for transcendental numbers, will be considered later
(see \S 4.3). 

\medskip

The so-called Mahler Problem (see
[\refmiwCetraro] 
\S 4.1)  is related to linear combination of logarithms $|b-\log a|$:  

\proclaim Conjecture {\numeroProblemeMahler} { \rm  (Mahler)}. There exists an absolute constant $c>0$ such that
$$
\Vert \log a\Vert>a^{-c}
$$
for all integers $a\ge 2$.

Equivalently, one has
$$
|a-e^b|>a^{-c}
$$
for some absolute constant $c>0$ for all integers $a$, $b$ $>1$.

\smallskip
A stronger  conjecture is suggested in  [\refmiwCetraro]
(4.1):
$$
\Vert \log a\Vert>(\log a)^{-c}
$$
for some absolute constant $c>0$ for all integers $a\ge 3$, or equivalently
$$
|a-e^b|>b^{-c}
$$
for some absolute constant $c>0$ for all integers $a$, $b$ $>1$.
So far the best known estimate is 
$$
|a-e^b|>e^{-c(\log a)( \log b)},
$$
so the problem is to replace in the exponent the product $(\log a)( \log b)$ by
the sum
$\log a+\log b$.

Explicit such lower bounds have interest in theoretical computer science [\refMullerTisserand].

\smallskip
Another topic which belongs to Diophantine approximation is the theory of  {\it
equidistributed  sequences}. For a positive integer $r\ge 2$, a {\it normal
number}  in base $r$ is a real number such that the sequence $(xr^n)_{n\ge 1}$ is
equidistributed modulo $1$. Almost all real numbers for Lebesgue measure are normal
(i.e., normal in basis $r$ for any $r>1$), but it is not known whether any
irrational  real algebraic number is normal to any integer basis, and it is also not known whether there is an
integer $r$ for which any number like
$e$, $\pi$, $\zeta(3)$, $\Gamma(1/4)$, $\gamma$, $G$, $e+\pi$, $e^\gamma$
is normal in basis
$r$ (see  [\refRauzy]). Further investigation by D.H.~Bailey and M.E.~Crandall
have been recently developed by J.C.~Lagarias in [\refLagarias].

The
digits of the  expansion (in any basis $\ge 2$) of an irrational real algebraic
 number should be equidistributed -- in particular any digit should
appear infinitely often. But even the following special case is unknown.

\proclaim Conjecture {\numeroconjectureMahler} { \rm  (Mahler)}.
Let $(\varepsilon_n)_{n\ge 0}$ be a sequence of elements in $\{0,1\}$. Assume 
that the   real number 
$$
\sum_{n\ge 0}\varepsilon_n 3^{-n}  
$$
is irrational, then it is transcendental.

\goodbreak

\section{3. \sectiontrois}

When $K$ is a field and $k$ a subfield,  we denote by $\trdeg_k K$
the transcendence degree of the extension $K/k$. In the case $k=\bQ$ we
write simply $\trdeg K$  (see [\refLangAlgebra] Chap.~VIII, \S 1).

\subsection{3.1. \subsectiontroispointun}

We concentrate here on problems related to transcendental number theory. 
To start with, we consider the classical exponential function
$e^z=\exp(z)$. A recent reference on this topic is [\refmiwGL].

Schanuel's Conjecture is a simple but far reaching statement
-- see the historical note to Chap.~III of [\refLangITN].

\proclaim Conjecture {\numeroConjectureSchanuel} { \rm   (Schanuel)}. 
Let 
$x_1,\ldots,x_n$ be
$\bQ$-linearly independent complex numbers. Then the 
transcendence degree over $\bQ$ of the field
$\bQ\bigl(x_1,\ldots,x_n,e^{x_1},\ldots,e^{x_n}
\bigr)$
is at least $n$.

According to S.~Lang
([\refLangITN] p.~31): ``From this statement, one would get most statements
about algebraic independence of values of
$e^t$ and
$\log t$ which one feels to be true''. See also [\refLangTNDA] p.~638--639 
and [\refRibenboimNumbersFriends] \S 10.7.G. For instance the following
statements [\refGelfondCRAS] are consequences of Conjecture
\numeroConjectureSchanuel.

\smallskip
\item{(?)}{\it  Let $\beta_1,\ldots,\beta_n$ be $\bQ$-linearly independent
algebraic numbers and let $\log\alpha_1,\ldots,\log\alpha_m$ be $\bQ$-linearly 
independent logarithms of algebraic numbers. Then the numbers
$$
e^{\beta_1},\ldots,e^{\beta_n},
\; \log\alpha_1,\ldots,\log\alpha_m
$$
are algebraically independent  over $\bQ$.}

\smallskip
\item{(?)}{\it  
Let $\beta_1,\ldots,\beta_n$ be algebraic
numbers with $\beta_1\not=0$ and let $\log\alpha_1,\ldots,\log\alpha_m$ be
 logarithms of algebraic numbers with $\log\alpha_1\not=0$ and 
$\log\alpha_2\not=0$. Then the numbers
$$
e^{\beta_1 e^{\beta_2 e^{{\adots}^{{\mathstrut}^{\beta_{n-1} e^{\beta_n}}}}}}
\and
\alpha_1^{\alpha_2^{{\adots}^{{\mathstrut}^{\alpha_m}}}}
$$
are transcendental, and there is no nontrivial algebraic relation between
such numbers.}

\smallskip\noindent
A quantitative refinement of  Conjecture {\numeroConjectureSchanuel}
is suggested in [\refmiwSurveyIA] Conjecture 1.4.

A quite interesting approach to Schanuel's Conjecture is given in 
[\refRoyACVEF] where D.~Roy states the next conjecture which he shows to be
equivalent to Schanuel's one. Let $\calD$ denote the derivation
$$
\calD = {\del\over\del X_0} + X_1 {\del\over\del X_1}
$$
over the ring $\bC[X_0,X_1]$. The {\it  height} of a 
polynomial $P\in\bC[X_0,X_1]$ is defined as the maximum of the absolute
values of its coefficients.  

\proclaim Conjecture {\numeroConjectureRoyRome} { \rm   (Roy)}.
Let $k$ be a positive integer, $y_1,\ldots,y_k$ complex numbers which are
linearly independent over $\bQ$,  $\alpha_1,\ldots,\alpha_k$ nonzero
complex numbers and $s_0,s_1,t_0,t_1,u$  positive real numbers
satisfying
$$
\max\{1,t_0,2t_1\} < \min\{s_0,2s_1\} 
  \and 
\max\{s_0,s_1+t_1\} < u < {1\over 2}(1+t_0+t_1).
$$
Assume that, for any sufficiently large positive integer $N$, there
exists a nonzero polynomial $P_N\in\bZ[X_0,X_1]$ with partial degree 
$\le N^{t_0}$ in $X_0$, partial degree $\le N^{t_1}$ in $X_1$ and 
height $\le e^N$ which satisfies
$$
\Bigg| \big(\calD^kP_N\big) 
       \Big(\sum_{j=1}^k m_jy_j, 
             \prod_{j=1}^k\alpha_j^{m_j} \Big) \Bigg|
  \le \exp(-N^u)
$$
for any nonnegative integers $k$, $m_1,\ldots,m_k$ with $k\le N^{s_0}$ 
and $\max\{m_1,\ldots,m_k\}\le N^{s_1}$.  Then, we have
$$
\trdeg \bQ(y_1,\ldots,y_k, \alpha_1,\ldots,\alpha_k) \ge k.
$$
\par

This work of Roy also provides an interesting connexion with other open
problems related to Schwarz Lemma for complex functions of several variables
(see [\refRoyIFAF] Conjectures 6.1 and 6.3).

The most important special case of Schanuel's Conjecture is the {\it 
Conjecture  of algebraic independence of logarithms of algebraic
numbers:}

\proclaim Conjecture {\numeroConjectureialogs}  { \rm  
(Algebraic Independence of Logarithms of Algebraic Numbers)}.
Let $\lambda_1,\ldots, \lambda_n$ be
$\bQ$-linearly independent complex numbers. Assume that the
numbers 
$e^{\lambda_1},\ldots,e^{\lambda_n}$ are algebraic.
Then the
numbers $\lambda_1,\ldots, \lambda_n$  are algebraically
independent  over $\bQ$.

An interesting reformulation of Conjecture {\numeroConjectureialogs} is due to
D.~Roy [\refRoyActaMath]. Denote by  $\calL$ the set of complex numbers $\lambda$
for which $e^\lambda$ is algebraic. Hence $\calL$ is a
$\bQ$-vector subspace of $\bC$. Roy's statement is: 

\smallskip
\item{(?)}{\it For any algebraic
subvariety $V$ of $\bC^n$ defined over the field $\Qbar$ of algebraic 
numbers, the set $V\cap \calL^n$ is the union of the sets $E\cap \calL^n$,
where $E$ ranges over the set of vector subspaces of $\bC^n$ which are
contained in $V$.}

\smallskip\noindent
Such a statement is reminiscent of several conjectures of Lang in 
Diophantine geometry (e.g. [\refLangEncyclopaedie]
Chap.~I, \S 6, Conjectures 6.1 and 6.3).

Not much is known on the algebraic independence of
logarithms of algebraic numbers, apart from the work of D.~Roy on the
rank of matrices whose entries are either logarithms of algebraic
numbers, or more generally linear combinations of logarithms of algebraic
numbers. We refer to [\refmiwGL] for a detailed study of this question as well as
related ones.

Conjecture {\numeroConjectureialogs} has many consequences.  
The next three ones are suggested by the work of D.~Roy ([\refRoySchanuel]
and [\refRoyMatrices]) on matrices whose entries are linear combinations of
logarithms of algebraic numbers (see also  [\refmiwGL] Conjecture 11.17,
\S 12.4.3 and Exercise 12.12).
 
Consider  the $\Qbar$-vector space $\calLtilde$
spanned by $1$ and $\calL$. In other words $\calLtilde$ is the set of complex
numbers which can be written
$$
\beta_0+\beta_1\log\alpha_1+\cdots+\beta_n\log\alpha_n,
$$
where $\beta_0,\beta_1,\ldots,\beta_n$ are algebraic numbers,
$\alpha_1,\ldots,\alpha_n$ are nonzero algebraic numbers, and 
finally
$\log\alpha_1,\ldots,\log\alpha_n$ are logarithms of 
$\alpha_1,\ldots,\alpha_n$ respectively.

\proclaim Conjecture 
{\numeroConjectureQuatreExponentiellesForte} { \rm   (Strong Four
Exponentials  Conjecture)}.
Let
$x_1,x_2$ be two 
$\Qbar$-linearly independent complex numbers and 
$y_1,y_2$ be also two 
$\Qbar$-linearly independent complex numbers. Then at  least
one of the four numbers  
$x_1y_1$, $x_1y_2$, $x_2y_1$, $x_2y_2$  does not belong to  
$\calLtilde$.

The following special case is also open:

\proclaim
Conjecture {\numeroConjectureCinqExponentiellesForte} { \rm  
(Strong Five Exponentials
Conjecture)}.
Let $x_1, x_2$ be two $\bQ$-linearly independent complex 
numbers and 
$y_1, y_2$ be also two $\bQ$-linearly independent complex 
numbers.  Further let $\beta_{ij}$ \ ($ i=1,2$, $j=1,2$),
$\gamma_1$ and
$\gamma_2$ be six algebraic numbers with $\gamma_1\not=0$. 
Assume that  the five numbers
$$
e^{x_1y_1-\beta_{11}},\; e^{x_1y_2-\beta_{12}},\; 
e^{x_2y_1-\beta_{21}},\; e^{x_2y_2-\beta_{22}},\;  
e^{(\gamma_1 x_1/x_2)-\gamma_2} 
$$
are algebraic.
Then all five exponents vanish:
$$
x_iy_j=\beta_{ij} \quad  (i=1,2, \quad    j=1,2)
\quad  \hbox{ and }\quad \gamma_1 x_1=\gamma_2 x_2.
$$.

A consequence of Conjecture {\numeroConjectureCinqExponentiellesForte} is the solution of the open problem of the transcendence of the number $e^{\pi^2}$, and more generally of $\alpha^{\log\alpha}=e^{\lambda^2}$ when $\alpha$ is a non zero algebraic number and $\lambda=\log\alpha$ a non zero logarithm of $\alpha$.

The next conjecture is proposed in [\refRoyActaMath].

\proclaim Conjecture {\numeroConjectureDRoyA} { \rm   (Roy)}.
For any $4\times 4$ skew-symmetric 
matrix  $\ttM$ with 
entries in 
$\calL$ and rank $\le 2$, either the rows of $\ttM$ are 
linearly dependent over $\bQ$, or the column space of
$\ttM$ contains a nonzero element of $\bQ^4$.

Finally a special case of  Conjecture {\numeroConjectureDRoyA} is the well
known  Four Exponentials
Conjecture due to Schneider [\refSchneiderLivre] Chap.~V, end of  \S 4,
Problem 1;  S.~Lang [\refLangITN] Chap.~II, \S 1;  [\refLangTNDA] p.~638 and
K.~Ramachandra [\refRamachandraAA\ II] \S 4.

\proclaim
Conjecture {\numeroConjectureQuatreExponentielles}
{ \rm  (Four Exponentials
Conjecture)}.   Let
$x_1,x_2$ be two $\bQ$-linearly independent complex numbers
and $y_1,y_2$ also two $\bQ$-linearly independent complex
numbers. Then at least one of the four numbers
$$
\exp(x_iy_j)\qquad(i=1,2,\; j=1,2)
$$
is transcendental.

The four exponentials  Conjecture can be stated as follows: 
{\it  consider a $2\times2$
matrix whose entries are logarithms of algebraic numbers:
$$
M=\pmatrix{\log\alpha_{11}& \log\alpha_{12}\cr
\log\alpha_{21}& \log\alpha_{22}\cr};
$$
assume that the two rows of this matrix are linearly independent
over $\bQ$ (in $\bC^2$), and also that the two columns are
linearly independent over $\bQ$; then the rank of this matrix is
$2$.}

We refer to [\refmiwGL] for a detailed discussion of this topic, including
the notion of {\it structural rank of a matrix} and the result, due to
D.~Roy, that Conjecture {\numeroConjectureialogs} is equivalent to a
conjecture on the rank of matrices whose entries are logarithms of 
algebraic numbers.

A classical problem on algebraic independence of algebraic powers of
algebraic numbers has been raised by A.O.~Gel'fond [\refGelfondConjecture] 
and Th.~Schneider [\refSchneiderLivre]  Chap.~V, end of \S 4, Problem 7. 
The data are an irrational algebraic number $\beta$ of degree $d$ and a
nonzero algebraic number
$\alpha$ with a nonzero logarithm $\log\alpha$. We write
$\alpha^z$ in place of $\exp\{z\log\alpha\}$. Gel'fond's
problem is:

\proclaim Conjecture {\numeroProblemeGelfond}
{ \rm  (Gel'fond)}.  The two numbers 
$$
\log\alpha\and  \alpha^\beta 
$$
are algebraically independent over $\bQ$.

Schneider's question is 

\proclaim Conjecture {\numeroProblemeSchneider}
{ \rm  (Schneider)}. The $d-1$
numbers 
$$
\alpha^\beta,\; \alpha^{\beta^2},\ldots,
\alpha^{\beta^{d-1}}
$$
are algebraically independent  over $\bQ$.

The first partial results in the direction of  Conjecture
{\numeroProblemeSchneider} are due to A.O.~Gel'fond [\refGelfondLivre].
For the more recent ones, see [\refNesterenkoPhilippon], Chap.~13 and 14.

Combining both questions {\numeroProblemeGelfond} and
{\numeroProblemeSchneider} yields a stronger conjecture:

\proclaim Conjecture {\numeroProblemeGelfondSchneider}
{ \rm  (Gel'fond-Schneider)}.  The $d$ numbers 
$$
\log\alpha,\; \alpha^\beta,\; \alpha^{\beta^2},\ldots,
\alpha^{\beta^{d-1}}
$$
are algebraically independent  over $\bQ$.

Partial results are known. They deal, more generally, with the values of
the usual exponential function at products $x_iy_j$, when
$x_1,\ldots,x_d$ and $y_1,\ldots,y_\ell$ are $\bQ$-linearly independent complex
(or $p$-adic) numbers.  The six exponentials Theorem states that, in
these circumstances, the $d\ell$ numbers $e^{x_iy_j}$  \ ($1\le i\le d$,
$1\le j\le \ell$) cannot all be algebraic, provided that $d\ell>d+\ell$.
Assuming stronger conditions on $d$ and $\ell$, namely $d\ell\ge
2(d+\ell)$, one deduces that two at least of these $d\ell$ numbers
$e^{x_iy_j}$ are algebraically independent  over $\bQ$. Other results are available
involving also either the numbers 
$x_1,\ldots,x_d$ themselves, or $y_1,\ldots,y_\ell$, or both. But an
interesting point is that, if we wish to get higher transcendence degree,
say to obtain that three at least of the numbers $e^{x_iy_j}$ are
algebraically independent  over $\bQ$, one needs so far a further assumption, which
is a measure of linear independence  over $\bQ$ for the tuple 
$x_1,\ldots,x_d$ as well as for the tuple $y_1,\ldots,y_\ell$. To remove
this so-called {\it  technical hypothesis} does not seem to be an easy 
challenge (see
[\refNesterenkoPhilippon] Chap.~14, \S 2.2 and \S 2.3).

The need for such a technical hypothesis seems to be connected with
the fact that the actual transcendence methods produce not only a
qualititative statement (lower bound for the transcendence degree), but
also quantitative statements (see \S 4). 

Several complex results have not yet been
established in the ultrametric situation. Two noticeable instances are:

\proclaim Conjecture {\numeroLindemannWeierstrasspadic} { \rm  ($p$-adic
analog of Lindemann-Weierstrass's Theorem)}. Let $\beta_1,\ldots,\beta_n$ 
be
$p$-adic algebraic numbers in the domain of convergence of the $p$-adic
exponential function $\exp_p$. Then the
$n$ numbers $\exp_p\beta_1,\ldots, \exp_p\beta_n$ are algebraically
independent over $\bQ$.

\proclaim Conjecture {\numeroGelfondpadic}  { \rm  ($p$-adic analog of an
algebraic independence result of Gel'fond)}.  Let $\alpha$ be a non-zero
algebraic number in the domain of convergence of the $p$-adic logarithm
$\log_p$, and let
$\beta$ be a
$p$-adic cubic  algebraic number, such that $\beta\log_p\alpha$ is in the
domain of convergence of the $p$-adic exponential function $\exp_p$.  
Then the two numbers 
$$
\alpha^\beta=\exp_p(\beta\log_p\alpha)
\and
\alpha^{\beta^2}=\exp_p(\beta^2\log_p\alpha)
$$ 
are algebraically independent
over
$\bQ$.

The
$p$-adic analog of  Conjecture {\numeroConjectureialogs} would solve 
Leopoldt's Conjecture on the $p$-adic rank
of the units of an algebraic number field
[\refLeopoldt] (see also [\refNarkiewiczEATAN] and [\refGras]), by proving the
nonvanishing of the
$p$-adic regulator. 

Algebraic independence results for the values of the exponential
function (or more generally for analytic subgroups of algebraic groups) in
several variables have already been established, but they are not yet
satisfactory. The conjectures stated  p.~292--293 of [\refmiwGAGDT] as
well as those of [\refNesterenkoPhilippon] Chap.~14, \S 2 are not yet 
proved. One of the main obstacles is the above-mentioned open problem with
the technical hypothesis.

The problem of extending  the Lindemann-Weierstrass Theorem to 
commutative algebraic groups is not yet completely solved (see conjectures
by P.~Philippon in [\refPhilipponLW]).

Algebraic independence proofs use elimination theory. Several methods are
available; one of them, developed by Masser, W\"ustholz and Brownawell,
rests on Hilbert Nulstellensatz. In this context we quote the following
conjecture of Blum, Cucker, Shub and Smale (see
[\refSmale]  
and [\refNesterenkoPhilippon] Chap.~16, \S 6.2), related to the open problem
``$P=NP~?$'' [\refClayMathematicalInstitute]~:

\proclaim Conjecture {\numeroBlumCuckerShubSmale}  { \rm 
(Blum, Cucker, Shub and Smale)}. 
Given an absolute constant $c$ and  polynomials $P_1,\ldots,P_m$ with a
total of $N$ coefficients and no common complex zeros, there is no program
to find, in at most $N^c$ step, the coefficients of polynomials 
$A_i$ satisfying B\'ezout's relation
$$
A_1P_1+\cdots+A_mP_m=1.
$$

In connexion with complexity in theoretical computer science,
W.D.~Brownawell suggests to investigate Diophantine approximation from a
new point of view in  [\refNesterenkoPhilippon] Chap.~16, \S 6.3.

Complexity theory may be related with a question raised by M.~Kontsevich and D.~Zagier in [\refKontsevichZagier]: they define a {\it period} as a   complex number whose real and imaginary part are values of absolutely convergent integrals of rational functions with rational coefficients over domains of $\bR^n$ given by polynomials (in)equalities with rational coefficients. The statement of Problem 3 in [\refKontsevichZagier]
is: {\it
exhibit at least one number which is {\rm not} a period.} This is the analog for {\it periods} of Liouville's Theorem for {\it algebraic numbers}. A more difficult question is to prove that specific numbers like
$$
e,\quad 1/\pi,\quad \gamma
$$
(where $\gamma$ is Euler's constant) are not periods. Since every algebraic number is a period, a number which is not a period is transcendental.

\medskip

Another important tool missing for transcendence proofs in higher
dimension is a Schwarz Lemma in several variables. The following conjecture is
suggested in [\refmiwSemLelong] \S 5. For a finite subset  $\Sigma$ of
$\bC^n$ and   a positive integer $t$, denote by $\omega_t(\Sigma)$ the least
total  degree of a nonzero polynomial $P$ in $\bC[z_1,\ldots,z_n]$ which vanishes
on
$\Sigma$ with multiplicity at least $t$:
$$
\displaylines{\qquad
\left({\partial\over \partial z_1}\right)^{\tau_1}\cdots
\left({\partial\over \partial z_n}\right)^{\tau_n}
P(z)=0\hfill
\cr
\hfill
\hbox{for any}\quad
z\in \Sigma \and
\tau=(\tau_1,\ldots,\tau_n)\in\bN^n
\with
\tau_1+\cdots+\tau_n<t.\quad\cr}
$$
Further, when $f$ is an analytic function in an open neighborhood of a closed 
polydisc
$|z_i|\le r$ \ ($1\le i\le n$) in $\bC^n$, denote by $\Theta_f(r)$ the average
mass of the set of zeroes of $f$ in the polydisc (see [\refBombieriLang]). 

\proclaim Conjecture {\numeroConjectureLemmeSchwarz}. 
Let $\Sigma$ be a finite subset of $\bC^n$ and 
$\varepsilon$ a positive number. There exists a positive number
$r_0(\Sigma,\varepsilon)$ such that, for any positive integer $t$ and any entire
function $f$ in $\bC^n$ which vanishes on $\Sigma$ with multiplicity $\ge t$, 
$$
\Theta_f(r)\ge \omega_t(\Sigma)-t\varepsilon
\for r\ge r_0(\Sigma,\varepsilon).
$$

The next question is to compute $r_0(\Sigma,\varepsilon)$. One may expect that for
$\Sigma$ a chunk of a finitely generated subgroup of $\bC^n$, say
$$
\Sigma=\bigl\{s_1y_1+\cdots+s_\ell y_\ell\; ;\;
(s_1,\ldots,s_\ell)\in\bZ^\ell,\; |s_j|\le S\;
(1\le j\le \ell)\bigr\}\subset\bC^n,
$$
an admissible value for the number $r_0(\Sigma,\varepsilon)$ will depend only on
$\varepsilon$, $y_1,\ldots,y_\ell$, but not on $S$. This would have interesting
applications, especially in the special case $\ell=n+1$.

Finally we refer to [\refChudnovsky] for a connexion between the numbers
$\omega_t(S)$ and Nagata's work on Hilbert's 14-th Problem.

\subsection{3.2. \subsectiontroispointdeux}

Many recent papers (see for instance [\refCartierBourbaki]) are devoted to the
study of algebraic relations among  ``multiple zeta values'' 
$$
\sum_{n_1>\cdots>n_k\ge 1}
n_1^{-s_1}\cdots n_k^{-s_k},
$$
(where $(s_1,\ldots,s_k)$ is a $k$-tuple of positive integers
with 
$s_1\ge 2$). The main Diophantine conjecture, suggested by the work
of D.~Zagier,  A.B.~Goncharov,
M.~Kontsevich, M.~Petitot, Minh Hoang Ngoc, K.~Ihara, M.~Kaneko and others (see [\refZagier], [\refCartierBourbaki] and  [\refZudilin]), is  that all such relations can be deduced from the linear and
quadratic ones arising from the {\it  shuffle} and {\it  stuffle} products (including
the relations occurring from the study of divergent series -- see [\refmiwMZV] for
instance). For $p\ge 2$, let $\gZ_p$ denote the
$\bQ$-vector subspace of
$\bR$ spanned by the real numbers
$\zeta(\us)$ satisfying $\us=(s_1,\ldots,s_k)$  and $s_1+\cdots+s_k=p$. Set
$\gZ_0=\bQ$ and
$\gZ_1=\{0\}$.
Then the
$\bQ$-subspace $\gZ$ spanned by all $\gZ_p$, $p\ge 0$, is a
subalgebra of $\bR$ and part of the Diophantine conjecture states:

\proclaim Conjecture {\numeroConjectureGoncharov} {\rm 
(Goncharov)}.  As a $\bQ$-algebra, $\gZ$ is the direct sum of
$\gZ_p$ for $p\ge 0$.

In other terms all algebraic relations should
be consequences of homogeneous ones, involving values $\zeta(\us)$ with
different
$\us$ but with the same weight $s_1+\cdots+s_k$. 

Assuming this conjecture \numeroConjectureGoncharov, the question of {\it
algebraic independence} of the numbers $\zeta(\us)$ is reduced to the
question of  {\it
linear independence} of the same numbers. The conjectural situation is
described by the next conjecture of Zagier  [\refZagier] on the dimension
$d_p$ of the
$\bQ$-vector space
$\gZ_p$.

\proclaim Conjecture {\numeroConjectureZagier} {\rm (Zagier)}. 
For $p\ge 3$ we have
$$
d_p=d_{p-2}+d_{p-3}
$$
with $d_0=1$, $d_1=0$, $d_2=1$.

That the actual dimensions of the spaces $\gZ_p$ are bounded above by the
integers which are defined inductively in Conjecture \numeroConjectureZagier\
has been proved by T.~Terasoma in [\refTerasoma], who expresses multiple
zeta values as periods of relative cohomologies and uses mixed Tate Hodge
structures (see also the work of A.G.~Goncharov referred to in
[\refTerasoma]). Further work on Conjectures {\numeroConjectureGoncharov}
and {\numeroConjectureZagier} is due to J.~\'Ecalle.
In case $k=1$ (values of the Riemann zeta
function)  the conjecture is: 

\proclaim Conjecture {\numeroConjectureZetaRiemann}. The numbers
$\pi$,
$\zeta(3),\zeta(5),\ldots,\zeta(2n+1),\ldots$  are algebraically
independent  over $\bQ$.

So far the only known results on this topic [\refFischler] are: 

\smallskip
\item{$\bullet$} {\it  $\zeta(2n)$ is transcendental for $n\ge 1$}
(because $\pi$ is transcendental and $\zeta(2n)\pi^{-2n}\in\bQ$),

\smallskip
\item{$\bullet$} {\it  $\zeta(3)$ is irrational} (Ap\'ery, 1978),

\smallskip
\noindent
and

\smallskip
\item{$\bullet$} {\it  For any $\varepsilon>0$   the
$\bQ$-vector space spanned by the $n+1$ numbers
$1,\zeta(3),\zeta(5),\ldots,\zeta(2n+1)$ has dimension
$$
\ge {1-\varepsilon\over 1+\log 2}\log n
$$
for  $n\ge n_0(\varepsilon)$} (see [\refRivoalCRAS] and [\refBallRivoal]).
For instance infinitely many of these numbers $\zeta(2n+1)$ \ ($n\ge 1$) are
irrational. W.~Zudilin proved that at least one of the four numbers 
$\zeta(5)$, $\zeta(7)$, $\zeta(9)$, $\zeta(11)$ is irrational.

\smallskip
\noindent
Further, more recent results are due to T.~Rivoal and W.~Zudilin: for instance
in a joint paper they proved that infinitely many numbers among
$$
\sum_{n\ge 1} {(-1)^n\over (2n+1)^{2s}}\qquad (s\ge 1)
$$
are irrational, 
but, as pointed out earlier, the irrationality of   Catalan's
constant
$G$ -- see (\numeroCatalan) -- is still an open problem.

\smallskip
It may turn out to be more efficient to work with a larger set of
numbers, including special values of multiple polylogarithms
$$
\sum_{n_1>\cdots>n_k\ge 1}
{z_1^{n_1}\cdots z_k^{n_k}\over
n_1^{s_1}\cdots n_k^{s_k}}\cdotp
$$
An interesting set of points $\uz=(z_1,\ldots,z_k)$ to consider is
the set of $k$-tuples consisting of roots of unity.
The function of a single variable
$$
\Li_{\us}(z)=\sum_{n_1>\cdots>n_k\ge 1}
{z^{n_1} \over
n_1^{s_1}\cdots n_k^{s_k}} 
$$
is worth of study from a Diophantine point of view. For instance, Catalan
constant mentioned above is the imaginary part of $\Li_2(i)$:
$$
\Li_2(i) =
\sum_{n\ge 1} {i^n\over n^2}=-{1\over 8}\zeta(2)+iG.
$$
Also no proof for the irrationality of the numbers
$$
\zeta(4,2)=
\sum_{n>k\ge 1} {1\over n^4k^2}=\zeta(3)^2-{4\pi^6\over 2835}\virgule
$$
$$
\Li_2(1/2)=
\sum_{n\ge 1}{1\over n^2 2^n}=
{\pi^2\over 12}-{1\over 2}(\log 2)^2
$$
and
$$
\Li_{2,1}(1/2)=
\sum_{n\ge k\ge 1}{1\over 2^n n^2 k}=\zeta(3)-{1\over 12}\pi^2\log 2,
\leqno{\hbox{(Ramanujan)}}
$$
is known so far.  

According to P.~Bundschuh [\refBundschuh], the transcendence of the numbers
$$
\sum_{n=2}^\infty {1\over n^s-1}
$$
for even $s\ge 4$ is a consequence of Schanuel's Conjecture
{\numeroConjectureSchanuel}. For $s=2$ the sum is $3/4$, and for $s=4$ the
value is $(7/8)-(\pi/4)\coth \pi$, which is a transcendental number since
$\pi$ and $e^\pi$ are algebraically independent  over $\bQ$
(Yu.V.~Nesterenko [\refNesterenkoPhilippon]).

Nothing is known on the arithmetic nature of the values of Riemann
zeta function at rational or algebraic points which are not integers.
\subsection{3.3. \subsectiontroispointtrois}

On the one hand, the transcendence problem of the values of Euler beta function at
rational points has been solved as early as 1940, by Th.~Schneider:
{\it for any rational numbers
$a$ and $b$ which are not integers and such that $a+b$ is not an integer, the
number
$$
\Beta(a,b)={\Gamma(a)\Gamma(b)\over \Gamma(a+b)}
$$
is transcendental.}  On the other hand, transcendence results for the values of
the gamma function itself are not so sharp: apart from G.V.~Chudnovsky's results
which imply the transcendence of $\Gamma(1/3)$ and $\Gamma(1/4)$ (and
Lindemann's result on the transcendence of $\pi$ which implies that
$\Gamma(1/2)=\sqrt{\pi}$ is also transcendental), not so much is known. For
instance, as we said earlier, there is no proof so far that $\Gamma(1/5)$ is
transcendental. This is because the Fermat curve of exponent $5$, viz.
$x^5+y^5=1$, has genus $2$. Its  Jacobian is an  abelian surface, and the
algebraic independence results known for elliptic curves like 
$x^3+y^3=1$ and $x^4+y^4=1$, which were sufficient for dealing with
$\Gamma(1/3)$ and $\Gamma(1/4)$, are not yet known for abelian varieties (see
[\refGrinspan]).

One might expect that
Nesterenko's results (see [\refNesterenkoPhilippon], Chap.~3) on the
algebraic independence of $\pi,\;
\Gamma(1/4),\;  e^\pi$ and of $\pi,\; \Gamma(1/3), \; e^{\pi\sqrt 3}$ should be
extended as follows:

\proclaim Conjecture {\numeroquestionGammaUnCinquieme}. At least three of the four numbers
$$
\pi,\; \Gamma(1/5),\;  \Gamma(2/5), \; e^{\pi\sqrt 5}
$$
are algebraically independent  over $\bQ$.

So the  challenge is to extend Nesterenko's results on modular
functions in one variable (and elliptic curves) to several variables
(and abelian varieties).

This may be one of the easiest questions to answer on this topic (but
it is still open). At the opposite, one may ask for a general statement
which would provide all algebraic relations between gamma values at
rational points. Here is a conjecture of Rohrlich
[\refLangDistributions]. Define
$$
G(z)={1\over\sqrt{2\pi}}\Gamma(z).
$$
According to the multiplication theorem of Gauss and
Legendre, [\refWhittakerWatson] \S 12.15, for each positive integer
$N$ and each complex number $x$ such that
$Nx\not\equiv 0 \pmod \bZ$,
$$
\prod_{i=0}^{N-1}G\left(x+{i\over N}\right)
=N^{(1/2)-Nx}G(Nx).
$$
The gamma function has no zero and defines a map from $\bC\setminus\bZ$ to $\bC^\times$. We restrict to $\bQ\setminus\bZ$ and we compose  with the canonical map
$\bC^\times\rightarrow \bC^\times/\Qbar^\times$ - which amounts to consider its
values modulo the algebraic numbers. The composite map has period $1$,
and the resulting mapping 
$$
\Gbar:{\bQ \over\bZ}\setminus\{0\}
\rightarrow {\bC^\times\over \Qbar^\times}
$$ 
is an odd {\it  distribution} on $(\bQ/\bZ)\setminus\{0\}$:
$$
\prod_{i=0}^{N-1}\Gbar\left(x+{i\over N}\right)
= \Gbar(Nx)\for x\in{\bQ \over\bZ}\setminus\{0\}
\and
\Gbar(-x)=\Gbar(x)^{-1}.
$$
Rohrlich's Conjecture
([\refLangDistributions],  
[\refLangCyclotomic] Chap.~II, Appendix, p.~66) asserts :

\proclaim Conjecture {\numeroConjectureRohrlich} { \rm  (Rohrlich)}.
$\Gbar$ is a universal odd distribution with values in groups where
multiplication by $2$ is invertible. 

In other terms, any multiplicative relation between gamma values at rational points
$$
 \pi^{b/2}\prod_{a\in\bQ}\Gamma(a)^{m_a}\in \Qbar
 $$
 with $b$ and $m_a$ in $\bZ$ can be derived for the standard relations satisfied by the gamma function.
This leads to the question whether the
distribution relations, the oddness relation and the functional equations 
of the gamma function generate the ideal over $\Qbar$ of all algebraic
relations among the values of $G(x)$ for $x\in\bQ$.

In [\refNesterenkoPhilippon] (Chap.~3,
\S 1, Conjecture 1.11) Yu.V.~Nesterenko proposed another conjectural extension of his algebraic independence result 
on Eisenstein series of weight $2$, $4$ and $6$:
$$
\eqalign{
P(q)&=1-24\sum_{n=1}^\infty{nq^n\over 1-q^n}
=1-24\sum_{n=1}^\infty\sigma_1(n)q^n
,\cr
Q(q)&=1+240\sum_{n=1}^\infty{n^3q^n\over 1-q^n}
=1+240\sum_{n=1}^\infty\sigma_3(n)q^n,\cr
R(q)&=1-504\sum_{n=1}^\infty{n^5q^n\over 1-q^n}
=1-504\sum_{n=1}^\infty\sigma_5(n)q^n.\cr}
$$

\proclaim Conjecture {\numeroConjectureNesterenko} { \rm  (Nesterenko)}.
Let $\tau\in\bC$ have positive imaginary part.
Assume that $\tau$ is not quadratic.
Set $q=e^{2i\pi\tau}$.
Then at least $4$ of the $5$ numbers
$$
\tau,\; q,\; P(q),\; Q(q),\; R(q)
$$
are algebraically independent.

Finally we remark that essentially nothing is known on the arithmetic nature of
the values of either the beta or the gamma function at algebraic irrational points.

A  wide range of open problems in transcendental 
number theory, including not only 
Schanuel's Conjecture
{\numeroConjectureSchanuel} and Rohrlich's Conjecture {\numeroConjectureRohrlich}  on the values of
the gamma function,  but also  a conjecture of Grothendieck
on the periods of an algebraic variety (see
[\refLangITN] Chap.~IV, Historical Note;
[\refLangTNDA] p.~650;
[\refAndreLivre] p.~6
 and 
[\refChambertLoir] \S 3), are
special cases of  very general conjectures due to Y.~Andr\'e [\refAndre],
which deal  with periods of mixed motives.  A discussion of Andr\'e's
conjectures for certain
$1$-motives related to products of elliptic curves and their connexions 
with elliptic and modular functions is given in [\refBertolin]. 
 Here is 
 a special case of the {\it elliptico-toric Conjecture} in  [\refBertolin]. 

\proclaim Conjecture {\numeroConjectureBertolin} { \rm  (Bertolin)}. 
Let $\calE_1,\ldots,\calE_n$ be pairwise non isogeneous elliptic curves with modular invariants $j(\calE_h)$. For $h=1,\ldots,n$, let $\omega_{1h},\omega_{2h}$ be a pair of fundamental periods of $\wp_h$ with  $\eta_{1h},\eta_{2h}$ the associated quasi-periods, $P_{ih}$ points on $\calE_h(\bC)$ and  $p_{ih}$  (resp. $d_{ih}$)  elliptic integrals of the first (resp. second) kind associated to  $P_{ih}$. Define $\kappa_h=[k_h:\bQ]$ and let $d_h$ be the dimension of the $k_h$-subspace of $\bC/(k_h\omega_{1h}+k_h\omega_{2h})$ spanned by $p_{1h},\ldots,p_{r_hh}$. Then the transcendence degree of the field
$$
\bQ\Bigl(\bigl\{
j(\calE_h), \omega_{1h},\omega_{2h}, \eta_{1h},\eta_{2h}, 
P_{ih}, p_{ih}, d_{ih}
\bigr\}_{1\le i\le r_h \atop 1\le h\le n}
\Bigr)
$$
is at least
$$
2\sum_{h=1}^n d_h
+4\sum_{h=1}^n \kappa_h^{-1}-n+1.
$$

A new approach to Grothendieck's Conjecture via  Siegel's   $G$-functions was
initiated in [\refAndreLivre] Chap.~IX. A development of this method led
Y.~Andr\'e to his conjecture on the special points on Shimura varieties 
[\refAndreLivre] Chap.~X, \S 4, which gave rise to the Andr\'e--Oort
Conjecture [\refOort] (for a discussion of this topic, including a precise
definition of ``Hodge type'', together with relevant references, see
[\refPaula]).

\proclaim Conjecture {\numeroConjectureAndreOort} { \rm  (Andr\'e--Oort)}. 
Let $\calA_g(\bC)$ denote the moduli space of principally polarized
complex  abelian varieties of dimension $g$. Let $Z$ be an irreducible
algebraic subvariety of  $\calA_g(\bC)$  such that the complex 
multiplication points on $Z$ are dense for the Zariski topology. Then $Z$
is a subvariety of 
$\calA_g(\bC)$  of Hodge type.

Conjecture {\numeroConjectureAndreOort} is a far reaching generalization 
of Schneider's Theorem on the transcendence of $j(\tau)$, where $j$ is
the modular invariant and $\tau$ an algebraic point in the Poincar\'e
upper half plane $\HPoincare$, which is not imaginary quadratic
([\refSchneiderLivre] Chap.~II, \S 4, Th.~17).  We also mention a
related conjecture of D.~Bertrand  (see [\refNesterenkoPhilippon] Chap.~1, \S 4
Conjecture 4.3) which may be viewed as a nonholomorphic
analog of Schneider's  result and which would answer the following question raised by N.~Katz:

\smallskip\item{(?)}{\it
Assume a lattice $L=\bZ\omega_1+\bZ\omega_2$ in $\bC$ has algebraic
invariants $g_2(L)$ and $g_3(L)$ and no complex multiplication. Does
this implies that the number
$$
G_2^*(L)=\lim_{s\rightarrow 0}
\sum_{\omega\in L\setminus\{0\}}
\omega^{-2} |\omega|^{-s}
$$
is transcendental?
}

\smallskip\noindent 
Many   open transcendence problems dealing with elliptic
functions   are consequences of  Andr\'e's
conjectures (see [\refBertolin]), most of which are likely to be very hard.
The next one, which is still open,  may be easier, since a number of
partial results are already known,  after the work of G.V.~Chudnovsky  and others
(see [\refGrinspan]).
 
\proclaim Conjecture {\numeroIndependancePeriodesElliptiques}. Given an
elliptic curve with Weierstrass equation
$y^2=4x^3-g_2x-g_3$, a nonzero period $\omega$,  the associated quasi-period
$\eta$ of the zeta function and a complex number $u$ which is not a pole of
$\wp$, we have
$$
\trdeg\bQ
\bigl(g_2,g_3,\pi/\omega,\wp(u),\zeta(u)-(\eta/\omega)u\bigr)\ge 2.
$$

Given
a lattice $L=\bZ\omega_1+\bZ\omega_2$ in $\bC$ with invariants
$g_2(L)$ and $g_3(L)$,  denote by
$\eta_i=\zeta_L(z+\omega_i)-\zeta_L(z)$  \ ($i=1,2$)   the
corresponding fundamental  quasi-periods of the Weierstrass zeta function.
  Conjecture
{\numeroIndependancePeriodesElliptiques} implies that  the transcendence degree
over
$\bQ$ of the field $\bQ\bigl(g_2(L),g_3(L),\omega_1,\omega_2,\eta_1,\eta_2\bigr)$   is
at least $2$.This would be optimal in the CM case, while in the non CM case,   we expect   it  to be
$\ge 4$. These lower bounds are given by   the  period conjecture of Grothendieck applied to an elliptic curve.

\smallskip 
According to [\refDiazRAMA] conjectures 1 and 2, p.~187, the following 
special case of  Conjecture {\numeroIndependancePeriodesElliptiques} can be stated in two equivalent ways: either in terms of
values of elliptic functions, or in terms of values of Eisenstein
series
$E_2$, $E_4$ and $E_6$ (which are $P$, $Q$ and $R$ in Ramanujan's
notation).

\smallskip
\item{(?)}{\it 
For any lattice $L$ in $\bC$ without complex multiplication  and for any nonzero 
period
$\omega$ of $L$,   
$$
\trdeg  \bQ\bigl(g_2(L),g_3(L),\omega/\sqrt{\pi},\eta/\sqrt{\pi}\bigr)\ge 2.
$$
 }

\smallskip
\item{(?)}{\it 
For any $\tau\in \HPoincare$ which is not imaginary quadratic, 
$$
\trdeg  \bQ\bigl(\pi E_2(\tau),
\pi^2 E_4(\tau),\pi^3 E_6(\tau)\bigr)\ge 2.
$$
}

\smallskip
\noindent
Moreover, each of these two statements  implies the following
one, which is stronger than a conjecture of Lang ([\refLangTNDA]
p.~652):

\smallskip
\item{(?)}{\it 
For any $\tau\in \HPoincare$ which is not imaginary quadratic, 
$$
\trdeg  \bQ\bigl(j(\tau),
j'(\tau),j''(\tau)\bigr)\ge 2.
$$
}

\smallskip
\noindent
Further related open problems are proposed by G.~Diaz in
[\refDiazJNTBx] and [\refDiazRAMA], in connexion with conjectures due to 
D.~Bertrand on the values of the modular function $J(q)$, where
$j(\tau)=J(e^{2i\pi \tau})$ 
(see [\refBertrandRamanujanJ] as well as  [\refNesterenkoPhilippon] Chap.~1,
\S 4 and Chap.~2, \S 4).

\proclaim Conjecture {\numeroConjectureBertrandA} { \rm  (Bertrand)}. 
Let $q_1,\ldots,q_n$ be nonzero algebraic numbers in the unit open disc such
that the $3n$ numbers
$$
J(q_i), \; DJ(q_i),\; D^2J(q_i)\qquad (i=1,\ldots,n)
$$
are algebraically dependent over $\bQ$. Then there exist two indices $i\not =
j$ \ ($1\le i\le n$, $1\le j\le n$) such that $q_i$ and $q_j$ are
multiplicatively dependent.

\proclaim Conjecture {\numeroConjectureBertrandB} { \rm  (Bertrand)}. 
Let $q_1$ and $q_2$ be two nonzero algebraic numbers in the unit open disc.
Suppose that there is an irreducible element $P\in\bQ[X,Y]$ such that
$$
P\bigl(J(q_1),J(q_2)\bigr)=0.
$$
Then there exist a constant $c$ and a positive integer $s$ such that
$P=c\Phi_s$, where $\Phi_s$ is the modular polynomial of level $s$. Moreover
$q_1$ and $q_2$ are
multiplicatively dependent.
 
\bigskip

Among Siegel's   $G$-functions are the algebraic functions. Transcendence
methods produce some information, in particular in connexion with Hilbert's
Irreducibility Theorem.
Let $f\in\bZ[X,Y]$ be a polynomial which is   irreducible in $\bQ(X)[Y]$. 
According to Hilbert's Irreducibility Theorem,  the set of positive integers
$n$ such that
$P(n,Y)$ is irreducible in $\bQ[Y]$ is infinite. Effective upper bounds for an
admissible value for $n$ have been investigated (especially by M.~Fried,
P.~D\`ebes and U.~Zannier), but do not yet answer the next: 

\proclaim Conjecture {\numeroConjectureHilbertIrreducibility}. 
Is there such a bound
depending polynomially on the degree and height of $P$? 

Such questions are also related to the {\it Galois inverse Problem}
[\refSerre].

Also the polylogarithms 
$$
\Li_s(z)=\sum_{n\ge 1}{z^m\over n^s}
$$
(where $s$ is a positive integer) are $G$-functions; unfortunately the 
Siegel-Shidlovskii method cannot be used so far to prove irrationality of the values of
the Riemann zeta function ([\refFeldmanNesterenko] Chap.~5, \S 7, p.~247).

\smallskip
With $G$-functions, the other class of analytic functions introduced by 
C.L.~Siegel in 1929 is the class of
$E$-functions, which includes the hypergeometric ones.  One main open question
is to investigate  the arithmetic nature of the values at algebraic points of
hypergeometric functions with {\it algebraic} parameters:
$$
\deuxfun{\alpha}{\beta}{\gamma}{z}=
\sum_{n\ge 0} {(\alpha)_n(\beta)_n\over (\gamma)_n}\cdot{z^n\over n!}
\virgule
$$ 
defined for $|z|<1$ and $\gamma\not\in\{0,-1,-2,\ldots\}$. 

\smallskip
In 1949,   C.L.~Siegel  ([\refSiegelLivre] Chap.~2, \S 9, p.~54 and 58; see also
[\refFeldmanShidlovskii] p.~62  and
[\refFeldmanNesterenko] Chap.~5,
\S 1.2) asked if {\it any $E$-function satisfying a linear differential equation with
coefficients in $\bC(z)$ can be expressed as a polynomial 
in $z$ and a finite number of hypergeometric $E$-functions or functions 
obtained from them by a change of variables of the form $z\mapsto \gamma z$ with algebraic
$\gamma$'s?}

\medskip

 Finally, we quote from [\refmiwSurveyIA]:
a folklore conjecture is that the zeroes of
the Riemann zeta function (say their imaginary parts, assuming it $>0$) are
algebraically independent. As suggested by J-P.~Serre, one might be tempted to
consider also
\smallskip
\item{$\bullet$}{The eigenvalues of the zeroes of the hyperbolic Laplacian in the
upper half plane modulo ${\rm SL}_2(\bZ)$ (i.e., to study the algebraic independence of
the zeroes of the Selberg zeta function).}
\smallskip
\item{$\bullet$}{The eigenvalues of the Hecke operators acting on the corresponding
eigenfunctions (Maass forms).}

\subsection{3.4. \subsectiontroispointquatre}

Many further open problems arise in transcendental number theory. 
An intriguing question is to investigate the arithmetic nature of real numbers
given in terms of power series involving  the Fibonacci sequence
$$
F_{n+2}=F_{n+1}+F_n,\qquad F_0=0,\quad F_1=1.
$$
Several results are due to  P.~Erd\H{o}s, R.~Andr\'e-Jeannin, C.~Badea,
J.~S\'andor, P.~Bundschuh, A.~Peth\H{o}, P.G.~Becker, T.~T\"opfer,  
D.~Duverney, Ku.~et Ke.~Nishioka, I.~Shiokawa, T.~Tanaka\dots\
It is known that the number
$$
\sum_{n=1}^\infty{1\over F_nF_{n+2}}=1
$$
is rational, while
$$
\sum_{n=0}^\infty{1\over F_{2^n}}={7-\sqrt{5}\over 2}\virgule\qquad
\sum_{n=1}^\infty{(-1)^n\over F_nF_{n+1}}={1-\sqrt{5}\over 2}
$$
and
$$
\sum_{n=1}^\infty{1\over F_{2n-1}+1}={\sqrt{5}\over 2}\virgule
$$
are  irrational algebraic numbers. Each of the numbers
$$
\sum_{n=1}^\infty{1\over F_n}\virgule\quad
\sum_{n=1}^\infty{1\over F_n+F_{n+2}}\and
\sum_{n\ge 1}{1\over F_1F_2\cdots F_n}
$$
is irrational, but it is not known whether they are algebraic or 
transcendental. The numbers 
$$
\sum_{n=1}^\infty{1\over F_{2n-1}}\virgule\quad
\sum_{n=1}^\infty{1\over F_n^2}\virgule\quad
\sum_{n=1}^\infty{(-1)^n\over F_n^2}\virgule\quad
\sum_{n=1}^\infty{n\over F_{2n}}\virgule\quad
\sum_{n=1}^\infty{1\over F_{2^n-1}+F_{2^n+1}} \and
\sum_{n=1}^\infty{1\over F_{2^n+1}}
$$
are all transcendental (further results of algebraic independence are
known).  The first challenge here is to formulate a conjectural statement which
would give a satisfactory  description of the situation.

\medskip

There is a similar situation for infinite sums $\sum_{n} f(n)$ where $f$
is a rational function [\refTijdemanSADA]: while
$$
\sum_{n=1}^\infty {1\over n(n+1)}=1
$$
and
$$
\sum_{n=0}^\infty
\left( {1\over 4n+1}-{3\over 4n+2}+{1\over 4n+3}+{1\over 4n+4}
\right)
=0
$$
are rational numbers, the sums
$$
 \sum_{n=0}^\infty {1\over (2n+1)(2n+2)}=\log2,
\qquad 
\sum_{n=0}^\infty {1\over (n+1)(2n+1)(4n+1)}={\pi\over 3} \virgule
$$
$$
\sum_{n=1}^\infty {1\over n^2}={\pi^2\over 6}\virgule
\qquad
\sum_{n=0}^\infty {1\over n^2+1}={1\over 2}
+{\pi\over 2}\cdot { e^\pi+e^{-\pi}
\over e^\pi-e^{-\pi}}\virgule
\qquad
\sum_{n=0}^\infty {(-1)^n \over n^2+1}= 
{2\pi \over e^\pi-e^{-\pi}} 
$$
and
$$
\sum_{n=0}^\infty {1\over (6n+1)(6n+2)
(6n+3)(6n+4)(6n+5)(6n+6)}={1\over 4320} 
(192\log 2-81\log 3-7\pi\sqrt{3})
$$
are transcendental.  The simplest example of the Euler sums $\sum_{n}n^{-s}$
(see
\S 3.2) illustrates the difficulty of the question: here again, even a
sufficiently general conjecture is missing. One may remark that there is no known algebraic irrational number of the form
 $$
 \sum_{n\ge 0\atop Q(n)\not=0} {P(n)\over Q(n)}\virgule
 $$
 with $P$ and $Q$ non zero polynomials having rational coefficients and $\deg Q\ge 2+\deg P$.
 
\medskip

The arithmetic study of the values of power series gives rise to many open
problems. We only mention a few of them.

The next question is due to K.~Mahler
[\refMahlerSuggestions]:

\proclaim Question {\numeroconjectureMahlerLiouville} { \rm 
(Mahler)}. Are there entire transcendental functions $f(z)$ such
that if $x$ is a Liouville number then so is $f(x)$?

The study of integral valued entire functions gives rise to several
open problems; we quote only one of them, which arose in the work of
D.W.~Masser and F.~Gramain on entire functions $f$ of one complex variable
which map the ring of Gaussian integers $\bZ[i]$ into itself. The initial
question (namely to derive an analog of P\'olya's Theorem in this setting)
has been solved by F.~Gramain in [\refGramain] (after previous work of
Fukasawa, Gel'fond, Gruman and Masser):  {\it If $f$ is not a polynomial,
then
$$
\limsup_{r\rightarrow\infty}
{1\over r^2}\log|f|_r\ge{\pi\over 2e}\cdot
$$}
Preliminary works on this estimate gave rise to the following 
problem, which is still unsolved. For each integer
$k\ge 2$, let $A_k$ be the minimal area of a closed disk in $\bR^2$ containing at least $k$ points of $\bZ^2$, and for $n\ge 2$ define
 %
$$
\delta_n=-\log n+\sum_{k=2}^n {1\over A_k}\cdotp
$$
The limit $\delta=\lim_{n\rightarrow \infty}\delta_n$ exists (it is an
analog in dimension $2$ of Euler constant), and the best known estimates
for it are [\refGramainWeber]
$$
1.811\dots<\delta<1.897\dots
$$
(see also  [\refFinch]).
F.~Gramain conjectures:
$$
\delta=1+{4\over\pi}\bigl(\gamma L(1)+ L'(1)\bigr),
$$
where $\gamma$ is Euler's constant and 
$$
L(s)=\sum_{n\ge 0}(-1)^n (2n+1)^{-s}
$$
is the $L$ function of the quadratic field $\bQ(i)$ (Dirichlet beta function). Since
$$
L(1)={\pi\over 4}\quad\hbox{and}\quad
L'(1)=\sum_{n\ge 0}(-1)^{n+1}\cdot{\log(2n+1)\over  2n+1}
=
{\pi\over 4}\bigl(3\log
 \pi+2\log 2+\gamma-4\log\Gamma(1/4)\bigr),
$$
Gramain's conjecture is equivalent to
$$
\delta=1+3\log
 \pi+2\log 2+2\gamma-4\log\Gamma(1/4)
=1.822825\dots
$$
Other problems related to the lattice $\bZ[i]$ are described in the
section ``{\it On the borders of geometry and arithmetic}'' of
[\refSierpinskiA].

\goodbreak

\section{4. \sectionquatre}


For  a nonzero polynomial $f\in\bC[X]$ of degree $d$:
$$
f(X)=a_0X^d+a_1X^{d-1}+\cdots+a_{d-1}X+a_d=
a_0\prod_{i=1}^d(X-\alpha_i),
$$
 define its {\it  usual height} by
$$
\rmH(f)=\max\{|a_0|,\ldots,|a_d|\}
$$
and its {\it  Mahler's
measure}  by 
$$
\rmM(f)=|a_0|\prod_{i=1}^d\max\{1,|\alpha_i|\}
=\exp\left(\int_0^1\log|f(e^{2i\pi t})|dt\right).
$$
The equality between these two formulae follows from Jensen's formula (see
[\refMahlerLivre] Chap.~I, \S 7, as well as  [\refmiwGL] Chap.~3 and 
[\refSchinzelBilkent]; the latter includes an extension to several
variables).

When $\alpha$ is an algebraic number with minimal polynomial
$f\in\bZ[X]$, define its {\it  Mahler's measure} by
$\rmM(\alpha)=\rmM(f)$ and its {\it  usual height} by
$\rmH(\alpha)=\rmH(f)$. Further, if $\alpha$ has degree $d$, 
define its {\it  logarithmic height} as 
$$
\rmh(\alpha)={1\over d}\log \rmM(\alpha).
$$
Furthermore, if $\alpha_1,\ldots,\alpha_d$ are the complex roots of $f$ (also called the {\it complex conjugates } of $\alpha$), then
the {\it house} of $\alpha$ is
$$
\house{\alpha}=\max\{|\alpha_1|,\ldots,|\alpha_d|\}.
$$
 
The height of an algebraic number is the prototype of a whole collection of height, like the height of  projective point 
$(1\colon \alpha_1\colon \cdots\colon
\alpha_n)\in\bP_n$ which is denoted by
$\rmh(1\colon \alpha_1\colon \cdots\colon
\alpha_n)$ (see for instance
[\refmiwGL]
\S 3.2) and the height of a subvariety $\hat{h}(V)$ (see for instance [\refDavidSinnouLille]
and
[\refDavidSinnou]).

Further, if
$\ualpha=(\alpha_1,\ldots,\alpha_n)$ is a $n$-tuple of 
multiplicatively independent algebraic numbers,
$\omega(\ualpha)$ denotes the minimum
degree of a nonzero  polynomial in
$\bQ[X_1,\ldots, X_n ]$ which vanishes at $\ualpha$.

A side remark is that Mahler's measure of a polynomial in a single
variable with algebraic coefficients is an algebraic number. The situation
is much more intricate for polynomials in several variables and gives rise to
further very interesting open problems [\refBoydExperimental,
\refBoydZakopane].

\subsection{4.1. \subsectionquatrepointun}

The smallest known value for $d\rmh(\alpha)$, which was found in 1933 by
D.H.~Lehmer, is $\log\alpha_0=0.162357\ldots$, where
$\alpha_0=1.176280\ldots$ is the real root\footnote{($^{5}$)}{\petitsanserif
Further
properties of this smallest known Salem number are described by  D.~Zagier in
his paper: Special values and functional equations of polylogarithms, Appendix
A  of {\it
 Structural properties of polylogarithms},  ed. L. Lewin,
 Mathematical Surveys and Monographs, vol. {\bf 37}, Amer. Math. Soc.
1991, p.~377--400.
}
of the degree
$10$ polynomial 
$$
X^{10}+X^9-X^7-X^6-X^5-X^4-X^3+X+1.
$$
D.H.~Lehmer asked
whether it is true that  for every positive $\varepsilon$ there
exists an algebraic integer $\alpha$ for which
$1<\rmM(\alpha)<1+\varepsilon$. 

\proclaim Conjecture {\numeroconjectureLehmer} { \rm  (Lehmer's Problem)}.
 There exists 
a positive absolute constant $c$ such that,
 for any 
nonzero algebraic number $\alpha$   which
is not a root of unity, 
$$
\rmM(\alpha)\ge  1+c.
$$ 
Equivalently, there exists 
a positive absolute constant $c$ such that,
 for any 
nonzero algebraic number $\alpha$ of degree at most $d$ which
is not a root of unity, 
$$
\rmh(\alpha)\ge  {c\over d}\cdotp
$$

Since $\rmh(\alpha)\le \log \house{\alpha}$,
the following statement [\refSchinzelZassenhaus] is a weaker
assertion than Conjecture
\numeroconjectureLehmer.

\proclaim Conjecture {\numeroconjectureSchinzelZassenhaus}
{ \rm   (Schinzel-Zassenhaus)}. There exists an
absolute  constant
$c>0$ such that, for any nonzero algebraic integer of degree 
$d$ which is not a root of unity,
$$
\house{\alpha}\ge 1+{c\over d}\cdotp
$$

Lehmer's 
Problem  is related to the multiplicative group $\bGm$. Generalizations to
$\bGm^n$ have been considered by many authors (see for instance
[\refBertrandAustralie]  and 
[\refSchmidtBilkent]).
In [\refAmorosoDavidCrelle] Conjecture 1.4, F.~Amoroso and S.~David extend
Lehmer's Problem {\numeroconjectureLehmer} to simultaneous approximation. 

\proclaim Conjecture {\numeroconjectureLehmerSimultanne}
{ \rm  (Amoroso-David)}.
For each positive integer $n\ge 1$ there
exists a positive number $c(n)$ having the
following property. Let $\alpha_1,\ldots,\alpha_n$ be
multiplicatively independent algebraic numbers. Define
$D=[\bQ(\alpha_1,\ldots,\alpha_n)\pps \bQ]$. Then
$$
\prod_{i=1}^n\rmh(\alpha_i)\ge{c(n)\over D}\cdotp
$$

The next statement ([\refAmorosoDavidCrelle] Conjecture 1.3 and 
[\refAmorosoDavidAA] Conjecture 1.3)
is stronger. 

\proclaim Conjecture {\numeroconjectureAmorosoDavid}
{ \rm  
(Amoroso-David)}. 
For each
positive integer
$n\ge 1$ there exists a positive number $c(n)$ such that, if
$\ualpha=(\alpha_1,\ldots,\alpha_n)$ is a $n$-tuple of 
multiplicatively independent algebraic numbers, then
$$
\rmh(1\colon \alpha_1\colon \cdots\colon
\alpha_n)\ge{c(n)\over
\omega(\ualpha)}\cdotp
$$

Many open questions are related to the height of subvarieties
[\refDavidSinnouLille],
[\refDavidSinnou]. 
The next one, dealing with the height of subvarieties of
$\bGm^n$ and proposed by F.~Amoroso and S.~David in [\refAmorosoDavidAA]
Conjecture 1.4 (see also  Conjecture  1.5 of [\refAmorosoDavidAA],
which is due to S.~David and P.~Philippon [\refDavidPhilipppon]), is
more general than  Conjecture {\numeroconjectureAmorosoDavid}:

\proclaim Conjecture {\numeroconjectureAmorosoDavidB}
{ \rm  
(Amoroso-David)}. 
For each integer $n\ge 1$ there exists a positive constant $c(n)$
such that, for any algebraic subvariety $V$ of $\bGm^n$ which is
defined over $\bQ$, which is $\bQ$-irreducible, and which is not a
union of translates of algebraic subgroups by torsion points, we have
$$
\hat{h}(V)\ge c(n)\deg(V)^{(s-\dim V-1)/(s-\dim V)},
$$
where $s$ is the dimension of the smallest algebraic subgroup of
$\bGm^n$ containing $V$.

\smallskip
Let $V$ be an
open subset of $\bC$. The {\it  Lehmer-Langevin constant} of $V$  is
defined as 
$$
L(V)=\inf \rmM(\alpha)^{1/[\bQ(\alpha):\bQ]},
$$
where $\alpha$ ranges over the set of nonzero and non-cyclotomic algebraic
numbers $\alpha$ lying with all their conjugates outside of $V$. It has been
proved by  M.~Langevin in 1985 that $L(V)>1$ as soon as $V$ contains a point
on the unit circle $|z|=1$.

\proclaim Problem {\numeroConjectureLangevin}. For $\theta\in(0,\pi)$, define
$$
V_\theta=\{re^{it}\; ;\;
r>0,\; |t|>\theta\}.
$$
Compute $L(V_\theta)$ in terms of $\theta$.
 
The solution is known only for a very few values of $\theta$: in 1995  G.~Rhin
and C.~Smyth [\refRhinSmyth] computed $L(V_\theta)$  for nine values of
$\theta$, including
$$
L(V_{\pi/2})=1.12\dots
$$

\smallskip
In a different direction, an analog of Lehmer's Problem has been raised for
elliptic curves, and more generally for abelian varieties.  Here is Conjecture
1.4 of [\refDavidHindry]. Let $A$ be an abelian variety defined over a number
field $K$ and equipped with a symmetric ample line bundle $\calL$. For any
$P\in A(\Qbar)$ define
$$
\delta(Q)=\min \deg(V)^{1/\codim(V)},
$$
where $V$ ranges over the proper subvarieties of $A$, defined over
$K$, $K$-irreducible  and
containing $Q$, while $\deg(V)$ is the degree of $V$ with respect to $\calL$.
Also denote by $\hat{h}_{\calL}$ the N\'eron-Tate canonical height on
$A(\Qbar)$ associated to $\calL$.

\proclaim Conjecture {\numeroConjectureLehmerVA} { \rm  (David-Hindry)}.
 There exists a positive constant
$c$, depending only on $A$ and $\calL$, such that for any $P\in A(\Qbar)$
which has infinite order modulo any abelian subvariety, we have
$$
\hat{h}_{\calL}(P)\ge c \delta(P)^{-1}.
$$

An extension of Conjecture {\numeroConjectureLehmerVA} to linearly
independent tuples is also stated in [\refDavidHindry] Conjecture
1.6.

The dependence on $A$ of these ``constants'' also raise interesting questions.
Just take an elliptic curve $E$ and consider the N\'eron-Tate height $\hat{h}(P)$
of a nontorsion rational point on a number field $K$. Several invariants are related
to $E$: the modular invariant   $j_E$, the
discriminant
$\Delta_E$
and Faltings height $h(E)$. A conjecture of Lang asserts
$$
\hat{h}(P)\ge c(K)\max\{1,h(E)\},
$$
while S.~Lang ([\refLangECDA] p.~92) and J.~Silverman ([\refSilverman] Chap.~VIII,
\S 10, Conjecture 9.9) conjecture 
$$
\hat{h}(P)\ge c(K)\max\{\log|N_{K/\bQ}(\Delta_E)|,\;h(j_E)\}.
$$
Partial results are known (J.~Silverman, M.~Hindry and J.~Silverman, S.~David), but
the conjecture is not yet proved.  

There is another abelian question related to
Mahler's measure. According to D.A.~Lind, Lehmer's Problem 
is known to be equivalent to the existence of a continuous endomorphism of
the infinite torus $(\bR/\bZ)^{\bZ}$ with finite entropy. A similar question
has been asked by P.~D'Ambros, G.~Everest, R.~Miles and T.~Ward
[\refAmbros-et-al] for elliptic curves, and it can be extended to abelian
varieties, and more generally to commutative algebraic groups.  

\subsection{4.2. \subsectionquatrepointdeux}

According to Dirichlet's box principle, for any irrational real number
$\theta$ there is an infinite set of rational numbers $p/q$ with $q>0$ such
that 
$$
\left|\theta-{p\over q}\right|\le {1\over q^2}\cdotp
\leqno{(\numeroFormuleDirichlet)}
$$
There are several extensions of this result. For the first one, we write 
(\numeroFormuleDirichlet) as $|q\theta-p|<1/q$ and we
replace
$q\theta-p$ by $P(\theta)$ for some polynomial $P$:

\item{(\numeroFormuleTiroirs)}{\it  Let $\theta$ be a real
number, $d$ and $H$ positive integers. There exists a nonzero
polynomial $P\in\bZ[X]$, of degree $\le d$ and usual height $\le H$, such that
$$ 
|P(\theta)|\le c H^{-d}
$$ 
where $c=1+|\theta|+\cdots+|\theta|^d$.
}

\smallskip\noindent
There is no assumption on $\theta$, but if $\theta$ is algebraic of degree
$\le d$ then there is a trivial solution!

A similar result holds for complex numbers, and more generally for
$\theta$ replaced by a
$m$-tuple
$(\theta_1,\ldots,\theta_m)\in\bC^m$  (see  for instance
[\refmiwGL] Lemma 15.11). For simplicity, we deal here only with the
easiest case. 

Another extension of (\numeroFormuleDirichlet) is interesting to consider,
where
$p/q$ is replaced by an algebraic number of degree $\le d$. If the  polynomial
$P$ given by (\numeroFormuleTiroirs) has a single simple  root  $\gamma$ close
to $\theta$, then
$$ 
|\theta-\gamma|\le c' H^{-d}
$$
where $c'$ depends only on $\theta$ and $d$. However, the root of $P$ which is
nearest $\gamma$ may be a multiple root, and may be not unique: this occurs
precisely when the first derivative $P'$ of $P$ has a small absolute value at
$\theta$. Dirichlet's box principle does not allow us to construct a
polynomial $P$ like in  (\numeroFormuleDirichlet) with a {\it  lower} bound
for
$|P'(\theta)|$. 

However  E.~Wirsing [\refWirsing] succeeded to prove the following theorem:
\smallskip

\item{(\numeroTheoremeWirsing)}{\it  
There exist two positive absolute constants $c>0$ and $\kappa$ such that, for
any trans\-cen\-den\-tal real number
$\theta$ and any positive integer $n$, there are infinitely many algebraic numbers
$\gamma$ of degree $\le n$ for which
$$
|\theta-\gamma|\le c\rmH(\gamma)^{-\kappa n}.
$$
}

\smallskip\noindent
E.~Wirsing himself obtained his estimate in 1960 with $c$ replaced by
$c(n,\varepsilon_n,\theta)$ and $\kappa n$  by
$(n/2)+2-\varepsilon_n$, where $\varepsilon_n\rightarrow 0$ as $n\rightarrow\infty$.
He conjectured that the exponent $\kappa n$  
can be replaced  by $n+1-\varepsilon$ [\refWirsing].
For $n=2$, H.~Davenport and W.M.~Schmidt in 1967 reached the  exponent $3$
without $\varepsilon$: {\it For any transcendental real number $\theta$, there
exists a positive real number $c(\theta)$ such that  the inequality
$$
|\theta-\gamma|\le c(\theta)\rmH(\gamma)^{-3}
$$
has infinitely many solutions $\gamma$ with $[\bQ(\gamma):\bQ]\le 2$.}
A conjecture of  Schmidt  ([\refSchmidtLNDA] Chap.~VIII, \S 3; see also
[\refBugeaudGraz] and [\refBugeaudJNT]) asserts that  (\numeroTheoremeWirsing)
holds with $\kappa n$ replaced by $n+1$:

\proclaim Conjecture {\numeroConjectureWirsingSchmidt} { \rm 
(Wirsing and Schmidt)}. For any positive integer
$n$ and any real number $\theta$ which is
either transcendental or else is algebraic of degree $>n$, there exists a
positive constant $c=c(n,\theta)$ with the following property:
there exist infinitely many algebraic numbers $\gamma$ of degree $\le n$
for which
$$
0<|\theta-\gamma|<c\rmH(\gamma)^{-n-1}.
$$

A third extension of (\numeroFormuleDirichlet) is to investigate the simultaneous rational approximation of successive powers of a real number. Let  $n\ge 2$ be an integer; denote by $\calE_n$ the set of real numbers which are not algebraic of degree $\le n$. For $\xi\in \calE_n$, let $\alpha_n(\xi)$ be the infimum of the set of real numbers $\alpha$ such that, for any sufficiently large real number $X$, there exists $(x_0,x_1,\ldots,x_n)\in\bZ^n $ satisfying  
$$
0<x_0\le X\quad 
\hbox{and}\quad \max_{1\le j\le n}|x_0\xi^j-x_j|\le X^{-1/\alpha}.
$$ From Dirichlet's box principle one deduces $\alpha_n(\xi)\le n$ for any $\xi\in \calE_n$ and any $n\ge 2$. Moreover, for any $n\ge 2$, the set of $\xi\in \calE_n$ for which $\alpha_n(\xi)<n$ has Lebesgue measure zero. H.~Davenport and W.M.~Schmidt proved in [\refDavenportSchmidt] that  $\alpha_2(\xi)\ge \gamma$ for any $\xi\in \calE_2$, where $\gamma=(1+\sqrt{5})/2=1.618\dots$\  It was expected that $\alpha_2(\xi)$ would be equal to $2$ for any $\xi\in \calE_2$, but D.~Roy [\refRoyNote]  has produced a $\xi\in \calE_2$ for which $\alpha_2(\xi)= \gamma$, showing that the result of Davenport and Schmidt is optimal. This raises a number of open problems and suggests to investigate the set 
$$
\calA_n=\{\alpha_n(\xi)\; ; \; \xi\in \calE_n\}.
$$
Recent results concerning the set $\calA_2$, by D.~Roy, Y.~Bugeaud and M.~Laurent, S.~Fischler, indicate a structure like the Markoff spectrum. For further references on this topic, see [\refBugeaudLivre].

Connexions between the problems of algebraic independence
we  considered in \S 1, the question of measures of linear
independence of logarithms of algebraic numbers discussed in \S 2,
and the notion of height which was introduced in
\S 3 (see for instance [\refmiwGL] Chap.~15) arise from the study of simultaneous approximation of complex numbers
by algebraic numbers. For a $m$-tuple $\ugamma=(\gamma_1,\ldots,\gamma_m)$ of
algebraic numbers, we define
$$
\mu(\ugamma)=
[\bQ(\ugamma):\bQ]\max_{1\le j\le m}\rmh(\gamma_j),
$$
so that for $m=1$ and $\gamma\in\Qbar$ we have simply
$\mu(\gamma)=\log\rmM(\gamma)$.

So far, relatiions between simultaneous approximation and algebraic
independence have been established only for small transcendence
degree. The missing link for large transcendence degree is 
given by the next statement (see [\refmiwGL] Conjecture 15.31;  
 [\refLaurentEger] \S 4.2 Conjecture 5;
[\refLaurentTokyo]  Conjecture 1;  [\refmiwGraz] Conjecture 2;  as well
as [\refRoyAASNL] Conjectures 1 and 2).

\proclaim Conjecture {\numeroPremiereConjectureGraz}.
Let $\utheta=(\theta_1,\ldots,\theta_m)$ be a $m$-tuple of
complex numbers. Define
$$
t=\trdeg  \bQ(\utheta)
$$
and assume $t\ge 1$. There
exist positive constants
$c_1$ and
$c_2$ with the following property.
 Let $(D_\nu)_{\nu\ge 0}$ and $(\mu_\nu)_{\nu\ge 0}$
be  sequences of real numbers satisfying
$$
c_1\le D_\nu\le \mu_\nu,\quad
D_\nu\le D_{\nu+1}\le 2 D_\nu,
\quad
\mu_\nu\le  \mu_{\nu+1}\le 2 \mu_\nu
 \qquad (\nu\ge 0).
$$
Assume also that the sequence $(\mu_\nu)_{\nu\ge 0}$ is
unbounded. Then for infinitely  many
$\nu$  there exists a $m$-tuple $(\gamma_1,\ldots,\gamma_m)$
of  algebraic numbers satisfying
$$
[\bQ(\ugamma):\bQ]\le D_\nu,
\quad
\mu(\ugamma)\le \mu_\nu
$$
and
$$
\max_{1\le i\le m}|\theta_i-\gamma_i|\le \exp\{-  c_2
D_\nu^{ 1/t }
\mu_\nu\}.
$$

There are two different (related) quantitative refinements to a
trans\-cen\-dence result: for a trans\-cen\-dental number $\theta$, either one
proves   {\it  a trans\-cen\-dence measure}, which is a lower bound for
$|P(\theta)|$ when $P$ is a nonzero polynomial with integer coefficients,
or else one proves a {\it  measure of algebraic approximation for} $\theta$,
which is a lower bound for $|\theta-\gamma|$ when $\gamma$ is an
algebraic number. In both cases such a lower bound will depend, usually,
on the degree (of the polynomial $P$, or on the algebraic number
$\gamma$), and on the height of the same.

Next, given several transcendental numbers $\theta_1,\ldots,\theta_n$,
one may consider either a measure of simultaneous approximation by
algebraic numbers, namely a lower bound for 
$$
\max\{|\theta_i-\gamma_i|\}
$$
when $\gamma_1,\ldots,\gamma_n$  are algebraic numbers, or a measure of
algebraic independence, which is a lower bound for 
$$
|P(\theta_1,\ldots,\theta_n)|
$$
when $P$ is a non-zero polynomial with integer coefficients. The first
estimate deals with algebraic points (algebraic sets of zero dimension),
the second with hypersurfaces (algebraic sets of codimension $1$). There
is a whole set of intermediate possibilities which have been
investigated by Yu.V.~Nesterenko and P.~Philippon, and are closely
connected. 

For instance Conjecture {\numeroPremiereConjectureGraz} deals with
simultaneous approximation by algebraic points; M.~Laurent and D.~Roy asked 
 general 
questions for the approximation by algebraic subsets of $\bC^m$, defined over
$\bQ$. For instance Conjecture
2 in  [\refLaurentTokyo] as well as the conjecture in \S 9 of
[\refRoyRMTNT] deal with the more general problem of approximation of points
in
$\bC^n$ by points located on
$\bQ$-varieties of a given dimension.

For an algebraic subset $Z$ of $\bC^m$, defined over $\bQ$, denote by $t(Z)$
the size of a Chow form of $Z$.

\proclaim Conjecture {\numeroConjectureLaurentRoy}
{ \rm  (Laurent-Roy)}. Let $\theta\in\bC^m$. There is a positive constant
$c$, depending only on $\theta$ and $m$, with the following property. Let $k$
be an integer with
$0\le k\le m$. For infinitely many integers
$T\ge 1$, there exists an algebraic set $Z\subset\bC^m$, defined over $\bQ$, 
of dimension $k$, and a point $\alpha\in Z$, such that
$$
t(Z)\le T^{m-k}\and
|\theta-\alpha|\le \exp\{-cT^{m+1}\}.
$$

Further far reaching  open
problems in this direction are proposed by P.~Philippon as Probl\`emes
7, 8 and 10 in [\refPhilipponJNT] \S 5.

\subsection{4.3. \subsectionquatrepointtrois}

We already suggested several questions related to linear independence
measures over the field of rational numbers for logarithms of rational
numbers (see Conjectures
{\numeroConjecturePhilipppon}, {\numeroConjectureLangW} and
{\numeroProblemeMahler}). Now that we have a notion of height for algebraic
numbers at our disposal, we extend our study to  linear independence
measures over the field of algebraic numbers for logarithms of  algebraic 
numbers.

The next statement is Conjecture 14.25 of [\refmiwGL].

\proclaim Conjecture {\numeroDeuxiemeConjectureGraz}.
There exist two positive absolute
constants $c_1$ and $c_2$ with the following property. Let
$\lambda_1,\ldots,\lambda_m$ be logarithms of algebraic
numbers with
$\alpha_i=e^{\lambda_i}$ \  ($1\le i\le m$), let $\beta_0,\ldots,
\beta_m$ be algebraic numbers, $D$ the degree of the number field
$$
\bQ(\alpha_1,\ldots,\alpha_m,\beta_0,\ldots,\beta_m)
$$ 
and finally 
let $h\ge 1/D$ satisfy
$$
h\ge \max_{1\le i\le m} \rmh(\alpha_i), \quad 
h\ge {1\over D}\max_{1\le i\le m}|\lambda_i|
\quad\hbox{and}\quad
h\ge \max_{0\le j\le m} \rmh(\beta_j).
$$
(1) Assume that the number
$$
\Lambda=\beta_0+\beta_1\lambda_1+\cdots+\beta_m\lambda_m
$$
is nonzero. Then
$$
|\Lambda|\ge 
\exp\bigl\{-c_1mD^2h \bigr\}.
$$
(2) Assume $\lambda_1,\ldots,\lambda_m$ are linearly independent 
over $\bQ$. Then
$$
\sum_{i=1}^m|\lambda_i-\beta_i|\ge 
\exp\bigl\{-c_2mD^{1+(1/m)}h \bigr\}.
$$

Assuming both Conjecture {\numeroPremiereConjectureGraz} and part 2
of Conjecture \numeroDeuxiemeConjectureGraz,
one deduces not only Conjecture \numeroConjectureialogs,
but also further special cases of  Conjecture {\numeroConjectureSchanuel}
(these
connexions are described in [\refmiwGraz] as well as [\refmiwGL] Chap.~15).

As far as part 1 of 
Conjecture {\numeroDeuxiemeConjectureGraz} is concerned, weaker estimates are
available (see [\refmiwGL] \S 10.4). Here is a much weaker (but still
open) statement than  either Conjecture
{\numeroConjectureLangW} or part 1 of 
Conjecture \numeroDeuxiemeConjectureGraz:

\proclaim Conjecture {\numeroMesureIndLinLo}. There exists a
positive absolute constant
$C$ with the following property.
 Let 
$\alpha_1,\ldots,\alpha_n$ be nonzero  algebraic numbers and
$\log\alpha_1,\ldots,\log\alpha_n$ logarithms of 
$\alpha_1,\ldots,\alpha_n$ respectively. Assume that the numbers
$\log\alpha_1,\ldots,\log\alpha_n$ are $\bQ$-linearly independent. Let
$\beta_0,\beta_1,\ldots,\beta_n$  be algebraic  numbers, not all of which are zero.
Denote by $D$  the degree of the number field
$$
\bQ(\alpha_1,\ldots,\alpha_n,\beta_0,\beta_1,\ldots,\beta_n)
$$
 over
$\bQ$. Further, let $A_1,\ldots,A_n$ and $B$ be positive
real numbers, each $\ge e$, such that
$$
\log A_j\ge\max\left\{\rmh(\alpha_j),\  {|\log\alpha_j|\over
D}\virgule \ {1\over D}\right\}
\quad (1\le j\le n),
$$
$$
B\ge\max_{1\le j\le n-1}
\rmh(\beta_j).
$$
Then the number
$$
\Lambda=\beta_0+\beta_1\log\alpha_1+\cdots+\beta_n\log\alpha_n
$$
satisfies
$$
|\Lambda|>
\exp\{-C^n 
D^{n+2}(\log A_1)\cdots  (\log A_n)(\log B+\log D) (\log D)\}.
$$

One is rather close to such an estimate (see [\refmiwCetraro] \S 5 and \S 6, as
well as [\refMatveev]). The result is proved now in the so-called rational case,
where 
$$
\beta_0=0\and \beta_i\in\bQ\for 1\le i\le n.
$$
In the general case, one
needs a further condition, namely 
$$
B\ge \max_{1\le
i\le n}\log  A_i.
$$
Removing this extra condition would enable one to prove that numbers like $e^\pi$
or
$2^{\sqrt 2}$ are not Liouville numbers.

These questions are the first and simplest ones concerning
transcendence measures,measures of Diophantine approximation,
measures of linear independence and measures of algebraic
independence. One may ask many further questions on this topic,
including an effective version of Schanuel's conjecture. It is
interesting to notice that in this case a ``technical condition''
cannot be omitted ([\refmiwSurveyIA] Conjecture  1.4).

Recall that  the rank of a prime ideal $ 
\gP\subset\bQ[T_1,\ldots,T_m]$ is the largest integer $  r\ge 0$ such
that there exists an increasing chain of prime ideals   
$$ 
(0)=\gP_0\subset\gP_1\subset\cdots \subset\gP_r=\gP.
$$ 
The  rank of an ideal $ \gI\subset\bQ[T_1,\ldots,T_m]$ is the minimum
rank of a prime ideal containing $ \gI$.  

\proclaim Conjecture {\numeroConjectureSchanuelEffectif}
{ \rm  (Quantitative Refinement of Schanuel's Conjecture)}.
Let $x_1,\ldots,x_n$ be $\bQ$-linearly independent
complex numbers.
Assume that for any $\varepsilon>0$, there exists a positive number $H_0$ such
that, for any $H\ge H_0$ and $n$-tuple $(h_1,\ldots,h_n)$ of rational integers
satisfying $0<\max\{|h_1|,\ldots,|h_n|\}\le H$, the inequality
$$
|h_1x_1+\cdots+h_nx_n|\ge 
\exp\bigl\{- H^\varepsilon\bigr\}
$$
holds.
 Let $d$ be a
positive integer. Then there exists a positive number
$C=C(x_1,\ldots,x_n,d)$ with the following property: for any integer $H\ge 2$ and
any $n+1$ tuple
$P_1,\ldots,P_{n+1}$ of polynomials in $\bZ[X_1,\ldots,X_n,Y_1,\ldots,Y_n]$ with
degrees
$\le d$ and usual heights $\le H$, which generate an ideal of
$\bQ[X_1,\ldots,X_n,Y_1,\ldots,Y_n]$ of rank
$n+1$, we have
$$
\sum_{j=1}^{n+1}\bigl|P_j(x_1,\ldots,x_n,e^{x_1},\ldots, e^{x_n})\bigr|\ge  H^{-C}.
$$

A consequence of Conjecture {\numeroConjectureSchanuelEffectif} is a
quantitative refinement to 
Conjecture {\numeroConjectureialogs}  on algebraic independence of logarithms of
algebraic numbers [\refmiwSurveyIA]: 

\smallskip
\item{(?)}{\it If $\log\alpha_1,\ldots,\log\alpha_n$
are $\bQ$-linearly independent logarithms of algebraic numbers and $d$ a
positive integers, there exists a constant $C>0$ such that, for any nonzero
polynomial $P\in\bZ[X_1,\ldots,X_n]$ of degree $\le d$ and height $\le H$ with
$H\ge 2$,
$$
|P(\log\alpha_1,\ldots,\log\alpha_n)|\ge H^{-C}.
$$
}

\subsection{4.4. \subsectionquatrepointquatre}

Let $K$ be a number field with a given real embedding. Let $V$ be a smooth 
variety over $K$. Denote by $Z$ the closure, for the real topology, of  $V(K)$ 
in
$V(\bR)$. In his paper [\refMazurTRP] on the topology of rational points,   
Mazur asks:

\proclaim Question {\numeroQuestionMazur}  { \rm  (Mazur)}. Assume that
$K=\bQ$ and that
$V(\bQ)$ is Zariski dense; is
$Z$  a union of connected components of $V(\bR)$?

An interesting fact is that Mazur asks this question in connexion with the
rational version of Hilbert's tenth Problem  (see [\refMazurQDUNT]
and [\refMazurSTRP]).

 The answer to question
{\numeroQuestionMazur}  is negative: an example is given in
[\refColliot] by J-L.~Colliot-Th\'el\`ene,   A.N.~Skorobogatov and  P.~Swinnerton-Dyer   of a smooth surface
$V$ over
$\bQ$, whose
$\bQ$-rational  points are Zariski-dense, but such that the closure $Z$ in
$V(\bR)$ of the set  of
$\bQ$-points is not a union of connected components.

However for the special case of abelian varieties, there are good reasons to
believe that the answer to question
{\numeroQuestionMazur}  
is positive. Indeed for this special case a reformulation of 
question {\numeroQuestionMazur}   is the following:

\smallskip

\item{(?)}{\it Let $A$ be a simple abelian variety over $\bQ$. Assume the
Mordell-Weil group $A(\bQ)$ has rank $\ge 1$. Then
$A(\bQ)\cap A(\bR)^0$ is dense in the neutral component $A(\bR)^0$
of $A(\bR)$.}

\smallskip
\noindent
This statement is equivalent to the next one:

\proclaim Conjecture {\numeroQuestionMazurVA}.
Let $A$ be a simple abelian variety over $\bQ$, 
$\exp_A:\bR^g\rightarrow A(\bR)^0$ the exponential map of the Lie group
$A(\bR)^0$  and
$\Omega=\bZ\omega_1+\cdots+\bZ\omega_g$ its kernel.  Let
$u=u_1\omega_1+\cdots+u_g\omega_g\in\bR^g$ satisfy $\exp_A(u)\in A(\bQ)$.
Then $1,u_1,\ldots,u_g$ are linearly independent over $\bQ$.

The following quantitative refinement of Conjecture {\numeroQuestionMazurVA} is
suggested in 
[\refmiwNagoya] Conjecture 1.1.

 For $\zeta=(\zeta_0:\cdots:\zeta_N)$ and
$\xi=(\xi_0:\cdots:\xi_N)$ in $\bP_N(\bR)$, write
$$
\dist(\zeta,\xi)={\displaystyle \max_{0\le i,j\le
N}|\zeta_i\xi_j-\zeta_j\xi_i| \over \displaystyle
\max_{0\le i \le N}|\zeta_i|\cdot
\max_{0\le  j\le N}| \xi_j |}\cdotp
$$

\proclaim Conjecture {\numeroconjectureNagoya}. 
Let $A$ be a simple Abelian
variety of  dimension
$g$ over a number field
$K$ embedded in $\bR$. Denote by $\ell$ the rank over $\bZ$ of the
Mordell-Weil group $A(K)$. For any $\varepsilon>0$, there exists $h_0>0$ (which
depends only on the Abelian variety 
$A$, the real number field 
$K$ and
$\varepsilon$) such that, for any $h\ge h_0$ and any $\zeta\in A(\bR)^0$, there is a
point $\gamma\in  A(K)$ with N\'eron-Tate height $\le
h$ such that
$$
\dist(\zeta,\gamma)
\le h^{-(\ell/2g)+\varepsilon}.
$$

Similar problems arise for commutative algebraic groups. Let us consider the
easiest case of a torus $\bGm^n$ over the field of real algebraic numbers. We
replace the simple abelian variety $A$ of dimension $g$ by the torus $\bGm^n$
of dimension $n$, the Mordell-Weil group $A(K)$  by a finitely
generated multiplicative subgroup of $\bigl(\Qbaretoile\bigr)^n$ and the
connected component $A(\bR)^0$ of the origin in $A(\bR)$ by
$\bigl(\bR_+^\times\bigr)^n$. The corresponding problem is then: given positive
algebraic numbers $\gamma_{ij}$
\ ($1\le i\le n$, $1\le j\le m$), consider the approximation of a tuple
$(\zeta_1,\ldots,\zeta_n)\in\bigl(\bR_+^\times\bigr)^n$ by tuples of algebraic
numbers of the form
$$
\bigl(\gamma_{11}^{s_1}\cdots \gamma_{1m}^{s_m},
\ldots,
  \gamma_{n1}^{s_1}\cdots \gamma_{nm}^{s_m}\bigr)
$$
with $\us=(s_1,\ldots,s_m)\in\bZ^m$.

Recently D.~Prasad 
[\refDipendra] investigated this question in terms of toric
varieties.

The qualitative density question is solved by the following statement, which
is a consequence of Conjecture \numeroConjectureialogs.

\proclaim Conjecture {\numeroConjectureDensiteTores}.
Let $m$, $n$, $k$ be positive integers  and $a_{ij\kappa}$ rational integers
($1\le i\le n$, $1\le j\le m$, $1\le\kappa\le k$). For
$\ux=(x_1,\ldots,x_k)\in\bigl(\bR_+^\times\bigr)^k$ denote by $\Gamma(\ux)$
the following finitely generated subgroup of $\bigl(\bR_+^\times\bigr)^n$:
$$
\Gamma(\ux)=\left\{
\left(
\prod_{j=1}^m \prod_{\kappa=1}^k x_k^{a_{ij\kappa}s_j}\right)_{1\le i\le n}
\; ;\;
\us=(s_1,\ldots,s_m)\in\bZ^m\right\}.
$$
Assume there exists $\ux \in\bigl(\bR_+^\times\bigr)^k$ such that
$\Gamma(\ux)$ is dense in $\bigl(\bR_+^\times\bigr)^n$. Then for any
$\ugamma=(\gamma_1,\ldots,\gamma_k)$ in 
$\bigl(\bR_+^\times\bigr)^k$ with $\gamma_1,\ldots,\gamma_k$ algebraic and
multiplicatively independent, the subgroup $\Gamma(\ugamma)$ is dense in 
 $\bigl(\bR_+^\times\bigr)^n$.

If there is a $\ux$ in $\bigl(\bR_+^\times\bigr)^k$ such that
$\Gamma(\ux)$ is dense in $\bigl(\bR_+^\times\bigr)^n$, then the
set of such $\ux$ is dense in 
$\bigl(\bR_+^\times\bigr)^k$. Hence again, loosely speaking,   Conjecture
{\numeroConjectureDensiteTores} means that logarithms of algebraic numbers
should behave like almost all numbers (see also
[\refLangEncyclopaedie] Chap.~IX, \S 7 p.~235). 

Conjecture {\numeroConjectureDensiteTores} would provide an effective
solution to the question raised by J-L.~Colliot-Th\'el\`ene and J-J.~Sansuc and
solved by D.~Roy (see [\refRoyInventiones]): 

\smallskip
\item{}{\it  Let $k$ be a number field of degree 
$d=r_1+2r_2$, where $r_1$ is the number of real embeddings and $r_2$ the
number of pairwise non-conjugate embeddings of $k$. Then there exists a
finitely generated subgroup
$\Gamma$ of $k^\times$, with rank $r_1+r_2+1$, whose image in
$\bR^{r_1}\times\bC^{r_2}$ is dense.}

\smallskip\noindent
The existence of $\Gamma$ is known, but the proof of D.~Roy does not yield an
explicit example. 

Density questions are closely related to transcendence questions. For instance
 the multiplicative subgroup of $\bR_+^\times$ generated
by
$e$ and $\pi$ is dense if and only if $\log\pi$ is
irrational (which is an open question). 

The simplest case of Conjecture {\numeroConjectureDensiteTores} is obtained
with $n=2$ and $m=3$. It reads as follows.

\smallskip

\item{(?)}{\it 
Let $\alpha_1,\alpha_2,\alpha_3$, $\beta_1,\beta_2,\beta_3$ be nonzero
positive algebraic numbers. Assume that for any  
$(a,b)\in\bZ^2\setminus\{(0,0)\}$, 
two at least of the three numbers
$$
\alpha_1^a\beta_1^b,\; 
\alpha_2^a\beta_2^b,\; \alpha_3^a\beta_3^b 
$$
are multiplicatively independent. Then the subgroup
$$
\Gamma=\bigl\{
\bigl(\alpha_1^{s_1}\alpha_2^{s_2}\alpha_3^{s_3},
\alpha_1^{s_1}\alpha_2^{s_2}\alpha_3^{s_3}\bigr)
\; ;\;
(s_1,s_2,s_3)\in\bZ^3\bigr\}
$$
of $\bigl(\bR_+^\times\bigr)^2$ is dense.}

\smallskip\noindent
It is easy to deduce this statement from the four exponentials Conjecture
\numeroConjectureQuatreExponentielles.

\medskip
The next question is to consider a quantitative refinement. 
Let $\Gamma$ be a finitely generated subgroup of $(\Qbar\cap\bR_+^\times)^n$
which is dense in  $\bigl(\bR_+^\times\bigr)^n$. Fix a set of generators
$\ugamma_1,\ldots,\ugamma_m$ of $\Gamma$. For $\us=(s_1,\ldots,s_m)\in\bZ^m$
and $1\le i\le n$ define
$$
\gamma_i(\us)=
\prod_{j=1}^m  \gamma_{ij}^{s_j}\in\Qbaretoile.
$$
The density assumption means that for any $\uzeta=(\zeta_1,\ldots,\zeta_n)\in
\bigl(\bR_+^\times\bigr)^n$ and any $\varepsilon>0$, there exists $\us\in\bZ^m$
such that
$$
\max_{1\le i\le n} |\gamma_i(\us)-\zeta_i|\le\varepsilon.
$$
We wish to bound $|\us|=\max_{1\le j\le m}|s_j|$ in terms of
$\varepsilon$.

We fix a compact neighborhood $\calK$ of the origin $(1,\ldots,1)$ in
$\bigl(\bR_+^\times\bigr)^n$. For instance
$$
\calK=\{\uzeta\in \bigl(\bR_+^\times\bigr)^n
\; ;\; 1/2\le |\zeta_i|\le 2 \;
(1\le i\le n)\}
$$
would do.

\proclaim Conjecture {\numeroConjectureDensiteToresQuantitative}. For any
$\varepsilon>0$ there exists $S_0>0$ (depending on $\varepsilon$,
$\gamma_1,\ldots,\gamma_m$ and $\calK$) such that, for any $S\ge S_0$ and any
$\uzeta\in\calK$, there exists $\us\in\bZ^m$ with $|\us|\le S$ and
$$
\max_{1\le i\le n} |\gamma_i(\us)-\zeta_i|\le
S^{-1-(1/n)+\varepsilon}.
$$

These questions raise a new  kind of Diophantine approximation problems.
 
\goodbreak\goodbreak

\section{5. \sectioncinq}

\subsection{5.1. \subsectioncinqpointun}

Among the motivations to study metric problems in Diophantine analysis
(not to mention secular perturbations in astronomy and the statistical
mechanics of a gas -- see [\refHarman]), one expects to guess the behavior
of certain classes of numbers (like algebraic numbers, logarithms of
algebraic numbers, and numbers given as values of classical functions,
suitably normalized [\refLangTNDA] p.~658 and 664). 

\smallskip

A first example is related to Wirsing-Schmidt's
Conjecture: in 1965, V.G.~Sprindzuk
showed that the conjecture {\numeroConjectureWirsingSchmidt}  is true for almost
all
$\theta$ (for Lebesgue measure). 
\smallskip

A second example is the
question of refining  Roth's Theorem: Conjecture {\numeroRothRaffine} is
motivated by Khinchine's Theorem  ([\refSprindzuk] Chap.~I, \S 1, Th.~1 p.~1)
which answers the question of rational Diophantine approximation for almost
all real numbers. In 1926  A.~Khinchine himself extended his result to
simultaneous Diophantine rational approximation
$$
\max_{1\le i\le n}|q\alpha_i -p_i|
$$
([\refSprindzuk] Chap.~I, \S 4, Th.~8 p.~28), and  in 1938 A.V.~Groshev
proved the first very general theorem of Khinchine type for systems of linear
forms
$$
\max_{1\le i\le n}|q_1\alpha_{i1}+\cdots+q_m\alpha_{im} -p_i|
$$
([\refSprindzuk] Chap.~I, \S 5, Th.~12, p.~33). Using the same heuristic, one
may extend  Conjecture {\numeroRothRaffine} in the context of simultaneous
linear combinations of algebraic numbers. 

In Conjecture {\numeroRothRaffine} (as well as in 
Khinchine's result for  almost
all real numbers) the function $q\psi(q)$ is assumed to be non-increasing.
A conjecture of Duffin and Schaeffer (see [\refSprindzuk] Chap.~1, \S 2,
p.~17 and [\refHarman]) would enable one to work without such a
restriction. Denote by $\varphi(n)$ Euler's function
$$
\varphi(n)=\sum_{1\le k\le n\atop \gcd(k,n)=1}1.
$$

\proclaim Conjecture {\numeroConjectureDuffinSchaeffer}
 {\rm (Duffin and Schaeffer)}.
Let
$\psi$ be a  positive real valued function.
Then, for
almost all $\theta\in\bR$, inequality (\numeroinegalitepourRothRaffine)
has  an infinite number of solutions in integers $p$ and $q$ with $q>0$
and $\gcd(p,q)=1$ if and only if the series
$$
\sum_{q=1}^\infty{1\over q} \psi(q)\varphi(q)
$$
diverges.

The Khinchine-Goshev Theorem has been extended to certain manifolds
(see [\refSprindzuk], [\refBernikDodson], as well as more recent papers
by V.~Bernick, M.~Dodson, D~Kleinbock and G.~Margulis). Further, connexions
between the metrical theory of Diophantine approximation on one hand,  
hyperbolic geometry, ergodic theory and dynamics of flows on homogeneous spaces
of Lie groups on the other, have been studied by several mathematicians,
including  D.~Sullivan, S.J.~Dani,
G.~Margulis and
D.~Kleinbock. Also S.~Hersonsky and F.~Paulin [\refHershonskyPaulin] recently
investigated the
Diophantine approximation properties of geodesic lines on the Heisenberg
group, which gives rise to new open questions, for instance to study
$$
\max_{1\le i\le n}|q\alpha_i -p_i|^{\kappa_i}
$$
when $\kappa_1,\ldots,\kappa_n$ are positive real numbers. 

\medskip
The set of real numbers with bounded partial quotients is countable. This is 
the set of real numbers which are badly approximable by rational numbers.  
Y.~Bugeaud asks a similar question for numbers which are badly approximable by
algebraic numbers  of bounded degree.

\proclaim Question {\numeroquestionBugeaud} {\rm (Bugeaud)}.
Let $n \ge 2$. Denote by
$\ZH_n$ the set of real numbers $\xi$ with the following property: there
exists $c_1(\xi) >0$ and 
$c_2(\xi)>0$ such that for algebraic number
$\alpha$ of degree $\le n$,
$$
|\xi-\alpha|\ge c_2(\xi) \rmH(\alpha)^{-n-1},
$$
and such that there are infinitely many algebraic numbers
$\alpha$ of degree $\le n$ with
$$
|\xi-\alpha|\le c_1(\xi) \rmH(\alpha)^{-n-1},
$$
Do the set $\ZH_n$ strictly contain the set of algebraic numbers of
degree $n+1$?
 
In connexion with the algebraic independence problems of \S 3.1, one would
like to understand better the behavior of real (or complex) numbers with
respect to Diophantine approximation by algebraic numbers of large degree
(see Conjecture \numeroPremiereConjectureGraz). A natural question is to
consider this question from a metrical point of view. Roughly speaking, what
is expected is that for almost all real numbers
$\xi$, the quality of approximation by algebraic numbers of degree $\le d$ and
measure $\le t$ is $e^{-dt}$. Here is a precise suggestion of Y.~Bugeaud
[\refBugeaudJNT].

For a real number $\kappa>0$, denote by $\calF_\kappa$  the set of real numbers
$\xi$ with the following property: for any $\kappa'$ with $0<\kappa'<\kappa$
and any
$d_0\ge 1$, there exists a real number $h_0\ge 1$ such that, for any $d\ge d_0$ and
any $t\ge h_0d$, the inequality
$$
|\xi-\gamma|\le e^{-\kappa'dt}
$$
has a solution $\gamma\in\Qbar$ with $[\bQ(\gamma):\bQ]\le d$ and
$\mu(\gamma)\le t$.

Also, denote by $\calF'_\kappa$  the set of real numbers $\xi$
with the following property: for any $\kappa'>\kappa$ there exists
$d_0\ge 1$ and $h_0\ge 1$ such that, for any $d\ge d_0$ and
any $t\ge h_0d$, the inequality
$$
|\xi-\gamma|\le e^{-\kappa'dt}
$$
has no solution $\gamma\in\Qbar$ with $[\bQ(\gamma):\bQ]\le d$ and
$\mu(\gamma)\le t$.

These definition are given more concisely in [\refBugeaudJNT]: for $t\ge d\ge
1$ denote by
$
\Qbar(d,t)
$
the set of real algebraic numbers $\gamma$ of degree $\le d$ and measure $\le t$. Then
$$
 \calF_\kappa = \bigcap_{\kappa' < \kappa} \, \bigcap_{d_0 \ge 1} \,
\bigcup_{h_0 \ge 1}
\,
\bigcap_{d \ge d_0} \, \bigcap_{t \ge h_0 d} \, 
\bigcup_{\gamma \in \Qbar(d, t) \cap \bR}
\, ]\gamma -e^{-\kappa' d t}, \gamma + e^{-\kappa' d t}[,
$$
$$
 \calF'_\kappa = \bigcap_{\kappa' > \kappa} \, \bigcup_{d_0 \ge 1} \,
\bigcup_{h_0 \ge 1} \,
\bigcap_{d \ge d_0} \, \bigcap_{t \ge h_0 d} \, 
\bigcap_{\gamma \in \Qbar(d, t) \cap \bR}
\, ]\gamma -e^{-\kappa' d t}, \gamma + e^{-\kappa' d t}[^c
$$
where $]a,b[^c$ denotes the complement of the intervall $]a,b[$.
According to Theorem 4 of [\refBugeaudJNT], there exists two positive constants
$\kappatilde$ and
$\kappatilde'$ such that, for almost all $\xi\in\bR$,
$$
\max \{\kappa>0 \; ;\; \xi \in  \calF_\kappa\}=\kappatilde
\and
\min \{\kappa>0 \; ;\; \xi \in  \calF'_\kappa\}=\kappatilde'.
$$
Further, 
$$
{1\over 850}\le \kappatilde\le\kappatilde'\le 1.
$$
Bugeaud's conjecture is $\kappatilde=\kappatilde'= 1$.

It is an important open question  to investigate the simultaneous
approximation of almost all tuples in
$\bR^n$  by algebraic tuples $\ugamma=(\gamma_1,\ldots,\gamma_n)$ in terms of
the  degree $[\bQ(\ugamma):\bQ]$ and the measure $\mu(\ugamma)$. Most often
previous authors devoted much attention to the dependence on the height, but
now it is necessary to investigate more thoroughly the behavior of the
approximation for large degree.

Further problems which we considered in the previous sections deserve to be
investigated under the metrical point of view. Our next example is a strong
quantitative form of Schanuel's Conjecture for almost all tuples
([\refmiwGraz] Conjecture 4). 

\proclaim Conjecture {\numeroQuatriemeConjectureGraz}. Let $n$ be a positive
integer. For almost all $n$-tuples $(x_1,\ldots,x_n)$, there are
positive constants
$c$ and $D_0$  (depending on $n$, $x_1,\ldots,x_n$ and $\varepsilon$),
with the following property. For any integer $D\ge D_0$, any real number
$\mu\ge D$ and any $2n$-tuple
$\alpha_1,\ldots,\alpha_n,\beta_1,\ldots,\beta_n$ of algebraic numbers
satisfying
$$
[\bQ( \alpha_1,\ldots,\alpha_n,\beta_1,\ldots,\beta_n):\bQ]\le D 
$$
and
$$
[\bQ( \alpha_1,\ldots,\alpha_n,\beta_1,\ldots,\beta_n):\bQ]
\max\left\{\rmh(\alpha_i),\; \rmh(\beta_i)\; ; \; 1\le i\le n\right\}\le\mu,
$$
 we have
$$
\max\left\{|x_i-\beta_i|,\;  
|e^{x_i}-\alpha_i|\; ; \; {1\le i\le n}\right\}
\ge \exp\{-c D^{1/(2n)}\mu\}.
$$

One may also expect that $c$ does not depend on $x_1,\ldots,x_n$.

\medskip

An open metrical problem of uniform distribution has been suggested
by P.~Erd\H{o}s to R.C.~Baker in 1973 (see [\refHarman] Chap.~5, p.~163).
It is a counterpart to a conjecture of Khinchine which was disproved by
J.M.~Marstrand in 1970. 
\smallskip

\item{(?)}{\it Let $f$ be a bounded measurable function with period $1$.
Is it true that
$$
\lim_{N\rightarrow\infty} {1\over \log N}\sum_{n=1}^N {1\over n} f(n\alpha)
=\int_0^1f(x)dx
$$
for almost all $\alpha\in\bR$?}

\smallskip\noindent

\subsection{5.2. \subsectioncinqpointdeux}

Let $K$ be a field and $\calC=K((T^{-1}))$ the field of Laurent
series on $K$.  The field $\calC$ has similar properties to the real
number field, where
  $\bZ$ is replaced by $K[T]$ and  $\bQ$ by
$K(T)$. An  absolute value on
$\calC$ is defined by selecting  $|T|>1$: we set
$|\alpha|=|T|^k$ if $\alpha=\sum_{n\in\bZ}a_nT^{-n}$ is a nonzero element of
$\calC$, where $k=\deg(\alpha)$ denotes the least index such that
$a_k\not=0$. Hence
$\calC$ is the completion of $K(T)$ for this absolute value.

A theory of Diophantine approximation has been developed on  $\calC$ 
in analogy with the classical one. If $K$ has zero characteristic, the result
are very similar to the classical ones. But if $K$ has finite characteristic,
the situation is completely different (see
[\refDeMathanLasjaunias] and [\refSchmidtSurveyAA]). It is not yet even clear
how to describe the situation from a conjectural point of view: a conjectural
description of the set of algebraic numbers for which a Roth type inequality
is valid is still missing.  Some algebraic elements satisfy a Roth type inequality, while for some others, Liouville's estimate is optimal. However, from a certain point of view, much more is known in the function field case, since the exact approximation
exponent is known for several classes of algebraic numbers. 
 
\medskip

There is also a transcendence theory over function fields. The starting point is
a work of Carlitz in the 40's. He defines functions on $\calC$ which behave like
an analog of the exponential function ({\it  Carlitz module}). A generalization
is due to V.G.~Drinfeld ({\it  Drinfeld modules}), and a number of results on the
transcendence of numbers related to these objects are known, going much further
than their classical (complex) counterpart. For example the {\it number}
 $$ 
 \prod_{p} (1-p^{-1})^{-1}
 $$
(in a suitable extension of a finite field), where $p$ runs over the monic irreducible polynomials over the given finite field, is known to be transcendental (over the field of rational functions on the finite field)  [\refAndersonThakur]; it may be considered as an analog of Euler's constant $\gamma$, since
 $$
\gamma=\lim_{s\rightarrow 1}\left(\zeta(s)-{1\over s-1}\right).
 $$
 However the theory is far from
being complete. An analog of Schanuel's Conjecture for Drinfeld modules is
proposed by W.D.~Brownawell in  [\refBrownawell], together with many further
related problems, including large transcendence degree, Diophantine geometry,
values of Carlitz-Bessel functions and values of gamma functions. 

For the study of Diophantine approximation, an important tool  (which is not
available in the classical number theoretic case) is the derivation $d/dT$.
In the transcendence
theory this gives rise to new questions which started to be investigated by
L.~Denis. Also in the function field  case interesting new questions arise
by considering several characteristics. So   Diophantine analysis for
function fields involve different aspects, some which are reminiscent of the
classical theory, and some which have no counterpart.

\medskip 

For the related transcendence theory involving
automata theory, we refer to the paper by D.~Thakur [\refThakur] 
(especially on p.~389--390) and to [\refAlloucheShallit] for the state of the art concerning the following
open problem:

\proclaim Conjecture {\numeroConjectureLoxtonVdP} { \rm  
(Loxton and van der Poorten)}. Let $(n_i)_{i\ge 0}$ be an
increasing sequence of positive integers. Assume  there is a prime
number $p$ such that the power series 
$$
\sum_{i\ge 0}z^{n_i}\in\bF_p[[z]]
$$
is algebraic over $\bF_p(z)$ and irrational (not in $\bF_p(z)$). Then the real 
number
$$
\sum_{i\ge 0}10^{-n_i} 
$$
is transcendental.
 
%

\goodbreak

\section{References}

\parskip=-3pt plus 1 pt minus .5 pt

\def\refmark#1{\medskip
\noindent
\hskip -3 true cm 
\hbox to 2.65 true cm {\sanserif
 [#1]\hfill } 
}
 
 \leftskip 2 true cm



\refmark{\refAlloucheShallit}
Allouche, J-P.; Shallit, J{.}--
{\it Automatic Sequences:  Theory, Applications, Generalizations.}
Cambridge University Press (2003), 

\refmark{\refAmbros-et-al}
D'Ambros, P.; Eversest, G.; Miles, R.; Ward, T{.} --
Dynamical systems arising from elliptic curves.
Colloq. Math. {\bf 84/85} \Numero{1} (2000), 95--107.

\refmark{\refAmorosoDavidCrelle} 
Amoroso, F{.};  David, S{.} --  
Le probl\`eme de Lehmer en dimension sup\'erieure. 
J{.} reine angew{.} Math.,  {\bf  513}  (1999), 145--179.
 
\refmark{\refAmorosoDavidAA} 
Amoroso, F{.};  David, S{.} --
Minoration de la hauteur normalis\'ee des hypersurfaces.
Acta Arith{.} {\bf  92} \Numero{4} (2000), 
339--366.

\refmark{\refAndersonThakur}
Anderson, G.W{.}; Thakur, {D.} --
Tensor powers of the Carlitz module and zeta values.  
Ann. Math., II. Ser. {\bf 132},  \Numero{1}  (1990), 159--191.

\refmark{\refAndreLivre}
Andr\'e, Y{.} --
{\it   $G$-functions and geometry. }
Aspects of
Mathematics, {\bf E13}. Friedr. Vieweg {\&} Sohn, Braunschweig, 1989. 

\refmark{\refAndre}
Andr\'e, Y{.} --
Quelques conjectures de transcendance issues de la g\'eom\'etrie 
alg\'e\-brique. 
Pr\'epu\-blication de l'Institut de Math\'ematiques de
Jussieu, {\bf  121} (1997), 18~p.
\hfill\break
{\seventt $<$http{$:$}//www.institut.math.jussieu.fr/$
\sim$preprints/index-1997.html$>$
}

\refmark{\refBakerLivre} 
Baker, A{.} --
{\it  Transcendental number theory.}
Cambridge Mathematical Library. Cambridge University Press, 
Cambridge, 1975. Second edition,
1990.

\refmark{\refBakerEger} 
Baker, A{.} -- 
Logarithmic forms and the $abc$-conjecture. 
Gy\H{o}ry, Kalman (ed.) et al., 
{\it  Number theory. Diophantine, computational and
algebraic aspects. }
Proceedings of the international conference, Eger, Hungary, 
July 29--August 2, 1996. Berlin: de Gruyter (1998),  37--44.

\refmark{\refBLSW}
Balasubramanian, R.; Langevin, M.; Shorey, T. N.; Waldschmidt, M{.}--
On the maximal length of two sequences of integers in arithmetic
progressions with the same prime divisors. 
Monatsh. Math. {\bf  121}
\Numero{4} (1996), 
295--307.

\refmark{\refBallRivoal} 
Ball, K.; Rivoal, T{.} -- Irrationalit\'e d'une
infinit\'e de valeurs de la fonction z\^eta aux entiers impairs.
Invent. math. {\bf 146} (2001) 1, 193--207.

\refmark{\refBernikDodson}
Bernik, V{.} I{.}; Dodson, M{.} M{.} --
{\it Metric Diophantine approximation on manifolds}. 
 Cambridge Tracts in Mathematics  {\bf 137}   
Cambridge University Press (1999).

\refmark{\refBertolin}
Bertolin, C{.} --
PŽriodes de 1-motifs et transcendance. 
J. Number Theory  {\bf 97}  (2002),  \Numero{2}, 204--221. 

\refmark{\refBertrandAustralie}
Bertrand, D{.} --
Duality on tori and multiplicative dependence relations.
 J. Austral. Math. Soc. Ser. A {\bf 62} \Numero{2} (1997), 
198--216. 

\refmark{\refBertrandRamanujanJ}
Bertrand, D{.} --
Theta functions and transcendence.
International Symposium on Number Theory (Madras, 1996). Ramanujan J. {\bf 1}
\Numero{4} (1997),  339--350. 

\refmark{\refBilu}
 Bilu, Yu{.} F{.} --
 Catalan's Conjecture (after Mih\u{a}ilescu).
{\it S\'em. Bourbaki, } 55$^{\rm e}$ ann\'ee, 2002-2003,
\Numero{909}, 36~p.
 \hfill\break
\hbox{\seventt $<$http{$:$}//www.math.u-bordeaux.fr/${ 
\sim}$yuri}$>$)

\refmark{\refBombieriLang}
Bombieri, E{.}; Lang, S{.} --
Analytic subgroups of group
varieties. 
Invent. Math. {\bf 11} (1970), 1--14.

\refmark{\refBoydExperimental}
Boyd, D{.} W{.} -- Mahler's measure and special values of
$L$-functions. Experiment. Math. {\bf 7}  (1998), \Numero 1, 37--82. 

\refmark{\refBoydZakopane}
Boyd, D{.} W{.} --
Mahler's measure and
special values of
$L$-functions---some conjectures. 
{\it Number theory in progress,}  Vol. {\bf 1}
(Zakopane-Ko\'scielisko, 1997),  27--34, de Gruyter, Berlin, 1999.

\refmark{\refBriggs}
Briggs, K{.} --
Some explicit badly approximable pairs.
\hfill\break
\hbox{\seventt $<$http{$:$}//arXiv.org/abs/math.NT/0211143$>$}


\refmark{\refBrowkin}
Browkin, J{.} -- 
The $abc$--conjecture.
{\it  Number theory,} 
 R.P. Bambah, V.C. Dumir and R.J. Hans Gill, Eds,
Hindustan Book Agency, New-Delhi,  Indian National Science Academy 
and Birkh\"auser-Verlag,
(1999), 75--105.

\refmark{\refBrownawell}
Brownawell, W.D{.} --
Transcendence in positive characteristic.
{\it Number theory},
Ramanujan Mathematical Society, January 3-6, 1996, 
Tiruchirapalli, India;
V.~Kumar Murty and Michel Waldschmidt, Eds.,
Amer. Math. Soc., Contemporary Math.  {\bf 210} (1998), 317--332.

\refmark{\refBugeaudGraz}
Bugeaud, Y{.} --
On the approximation by algebraic numbers with bounded degree.
{\it  Algebraic Number Theory and Diophantine Analysis,} 
F. Halter-Koch and R. F. Tichy, Eds, W{.} de Gruyter, Berlin, (2000), 47--53.

\refmark{\refBugeaudJNT}
Bugeaud, Y{.} --
Approximation par des nombres alg\'ebriques.
J{.} Number Theory,
{\bf  	84} \Numero{1} (2000),  15--33. 
\hfill\break
\hbox{\seventt 
doi:10:1006.jnth.2000.2517 JNT 84.1}

\refmark{\refBugeaudLivre}
Bugeaud, Y{.} --
{\it Approximation by algebraic numbers.}
Cambridge Tracts in Mathematics, 
Cambridge 
University Press, Cambridge-New York, 2004.

\refmark{\refBugeaudMignotte}
Bugeaud, Y{.}; Mignotte, M{.}    --
Sur l'Žquation diophantienne
$\frac{x^n-1}{x-1}= y^q$. II. 
C. R. Acad. Sci., Paris, SŽr. I,
Math. {\bf 328}, \Numero 9  (1999), 741--744.

\refmark{\refBugeaudShorey}
Bugeaud, Y{.}; Shorey, T{.} N{.}    --
On the number of solutions of 
the generalized Ramanujan--Nagell equation. 
J. Reine angew. Math. {\bf 539 } (2001), 55--74.

\refmark{\refBundschuh}
Bundschuh, P{.} --
Zwei Bemerkungen \"uber transzendente
Zahlen. 
Monatsh. Math. {\bf 88} \Numero 4 (1979),  293--304.

\refmark{\refCartierBourbaki}
 Cartier, P{.} --
Fonctions polylogarithmes, nombres polyz\^eta et groupes pro-unipotents.
{\it S\'em. Bourbaki, } 53$^{\rm e}$ ann\'ee, 2000-2001,
\Numero{885}, 36~p.
AstŽrisque {\bf 282} (2002), 137--173.

\refmark{\refCassels}
Cassels, J. W. S{.}-- 
{\it  An introduction to Diophantine approximation.}
Cambridge Tracts in Mathematics and Mathematical Physics, \Numero{\bf 
45}. Cambridge University Press, New York, 1957. 

\refmark{\refCatalan}
Catalan, E{.} --
Note extraite d'une lettre adress\'ee \`a l'\'editeur. 
J{.} reine
angew{.} Math.,  {\bf  27}  (1844), 192.

\refmark{\refChambertLoir}
Chambert-Loir, A{.} --
Th\'eor\`emes d'alg\'ebricit\'e en g\'eom\'etrie diophantienne d'apr\`es
J-B. Bost, Y.~Andr\'e, D \&\ G. Chudnovsky. 
{\it S\'em. Bourbaki, } 53$^{\rm
e}$ ann\'ee, 2000-2001,
\Numero{886}, 35~p.
AstŽrisque {\bf 282} (2002), 175--209 

\refmark{\refChudnovsky}
Chudnovsky, G. V{.} --
Singular points on complex hypersurfaces
and multidimensional Schwarz lemma. 
{\it Seminar on Number Theory, Paris
1979--80,} 29--69, Progr. Math., {\bf 12}, Birk\-h\"auser, Boston, Mass., 1981.

\refmark{\refPaula} 
Cohen, P.B{.} --
Perspectives de l'approximation Diophantienne et la transcendence. 
The Ramanujan Journal,
{\bf  7} \Numero{1-3} (2003), 367--384.
%

\refmark{\refCohn}
Cohn, H{.} --
 Markoff geodesics in matrix theory.
{\it Number theory with an emphasis on the Markoff
spectrum (Provo, UT, 1991),}  Lecture Notes in Pure and Appl. Math.,
{\bf 147}, Dekker, New York, (1993) 69--82.

\refmark{\refColliot}
Colliot-Th\'el\`ene, J-L.; Skorobogatov, A.N.;    
Swinnerton-Dyer, P{.} --
Double fibres and double covers: paucity of rational points.
Acta Arith. {\bf  79} \Numero{2} (1997), 113--135.

\refmark{\refCusickFlahive}
Cusick, T.W.; Flahive, M.E{.}--  
{\it  The Markoff and Lagrange Spectra}.
Mathematical Surveys and Monographs \Numero{\bf  30} 
 American
Mathematical Society, Providence, RI, 1989.


\refmark{\refCusickPomerance}
Cusick, T.W.; Pomerance, C{.} --
View-obstruction problems.
III. J. Number Theory {\bf 19} \Numero{2} (1984), 131--139. 

\refmark{\refDavenportSchmidt}
Davenport, H.; Schmidt, W.M. -- 
Approximation
to real numbers by algebraic integers. Acta Arith. {\bf 15}
(1969), 393-416.

\refmark{\refDavidSinnouLille}
David, S{.} -- Petits points, points rationnels.
{\it XXII\`emes Journ\'ees Arithm\'etiques de Lille, 2001.}

\refmark{\refDavidSinnou}
David, S{.} -- On the height of subvarieties of group varieties. 
The Ramanujan Journal,
to appear.

\refmark{\refDavidHindry}
David, S.; Hindry, M{.} -- 
Minoration de la hauteur de N\'eron-Tate sur les vari\'et\'es ab\'eliennes 
de type C.M.
J. Reine angew. Math. {\bf  529} (2000),
1--74.

\refmark{\refDavidPhilipppon}
David, S.; Philippon, P{.} -- 
Minorations des hauteurs normalis\'ees des sous-vari\'et\'es des
tores. 
Ann. Scuola Norm. Sup. Pisa Cl. Sci. (4) {\bf 28}  \Numero{3}
(1999),
 489--543.  

\refmark{\refDavisMatiyasevichRobinson}
Davis, M.; Matiyasevich, Y.; Robinson, J{.} --
Diophantine equations: a positive aspect of a negative solution.
{\it  Mathematical developments arising from Hilbert
problems,} (Proc. Sympos. Pure Math., {\bf  28}, Part 2, Northern
Illinois Univ., De Kalb, Ill., 1974). Amer. Math. Soc.,
Providence, R. I., (1976),  323--378. 

\refmark{\refDeMathanLasjaunias}
de Mathan, B{.}; Lasjaunias, A.  --
Differential equations and
Diophantine approximation in positive characteristic. 
Monatsh. Math. {\bf 128} \Numero{1}
(1999),  1--6. 

\refmark{\refDiazJNTBx} 
Diaz, G{.} --
La conjecture des quatre exponentielles et les
conjectures de D{.} Bertrand sur la fonction modulaire. 
 J. ThŽor. Nombres Bordeaux  {\bf  9}  \Numero{1} (1997), 
229--245.

\refmark{\refDiazRAMA} 
Diaz, G{.} --
Transcendance et ind\'ependance alg\'ebrique: liens entre les points de vue
elliptique et modulaire.
The Ramanujan Journal,
{\bf  4} \Numero{2} (2000), 157--199.

\refmark{\refDobrowolskiAA} 
Dobrowolski, E{.} -- 
On a question of Lehmer and the number of irreducible factors 
of a polynomial.  
Acta Arith{.} {\bf  34} \Numero{4} (1979), 
391--401.

\refmark{\refDujella} 
Dujella, A{.} -- 
There are only finitely many Diophantine quintuples.
J. Reine angew. Math., to appear. (28pp).
\hfill\break
\hbox{\seventt $<$http$:$//www.math.hr/${\scriptstyle
\sim}$duje/dvi/finit.dvi$>$}

\refmark{\refErdosMonthly} 
Erd\H{o}s, P{.} --
How many pairs of products of consecutive integers have the same prime
factors?
Amer. Math. Monthly {\bf  87} (1980), 391--392.

\refmark{\refErdosDurham} 
Erd\H{o}s, P{.} --
On the irrationality of certain series: problems
and results. 
{\it  New advances in transcendence theory (Durham, 1986)},
Cambridge Univ. Press, Cambridge-New York, (1988), 102--109.

\refmark{\refFaltings} 
Faltings, G{.} --
Diophantine Equations.
{\it Mathematics Unlimited - 2001 and Beyond.} 
Engquist, B.; Schmid, W., Eds., Springer  (2000), 449--454. 
 
\refmark{\refFeldmanNesterenko}
Fel'dman, N{.} I{.};
Nesterenko, Y{.} V{.} -- {\it  Number theory. IV.
Transcendental Numbers.}
Encyclopaedia of
Mathematical Sciences {\bf 
44}. Springer-Verlag,
Berlin, 1998.

\refmark{\refFeldmanShidlovskii} 
Fel'dman, N{.} I{.}; \v Sidlovski\u\i, A{.} B{.}  --
The development and present
state of the theory of
transcendental numbers.
(Russian) Uspehi Mat{.} Nauk
{\bf 22}  \Numero{3} (135) (1967), 
3--81;  Engl{.} transl{.} in
Russian Math{.} Surveys, {\bf
22}  \Numero{3}  (1967), 1--79.

\refmark{\refFinch}
Finch, S{.} R{.} --
{\it Mathematical Constants.}
 Encyclopedia of Mathematics and its Applications, {\bf 94}. Cambridge University Press, Cambridge,  2003.

\refmark{\refFischler}
Fischler, S{.} -- 
IrrationalitŽ de valeurs de zta [d'aprs ApŽry, Rivoal,\dots].
{\it S\'em. Bourbaki, } 55$^{\rm
e}$ ann\'ee, 2002-2003,
\Numero{910}, 35~p.

\refmark{\refFlattoLagariasPollington}
Flatto, L{.}; Lagarias, J{.} C{.}; Pollington, A{.} D{.} --
On the range of fractional parts 
{\rm $ \{ \xi (p/q)^n \} $.}
{Acta Arith.,  {\bf  70}
\Numero{2} (1995),  125--147.

\refmark{\refGalochkin}
Galochkin, A.I{.} --
Sur des estimations pr\'ecises par rapport \`a la hauteur de certaines
formes lin\'eaires.  
{\it   Approximations  Diophan\-tiennes  et  Nombres Transcendants,} 
(Luminy, 1982), Progr. Math., {\bf  31}, Birkh\"auser, Boston,  Mass.,
(1983), 95--98.

\refmark{\refGelfondCRAS} 
Gel'fond, A{.} O{.} -- Sur quelques r\'esultats nouveaux dans la 
th\'eorie des nombres transcendants.
C.R. Acad. Sc. Paris, S\'{e}r. A, {\bf 199} (1934), 259. 

\refmark{\refGelfondConjecture} 
Gel'fond, A{.} O{.} --
The approximation of algebraic numbers
by algebraic numbers and the theory of transcendental numbers.  
Uspehi Matem. Nauk (N.S.) {\bf  4} \Numero{4} (32) (1949), 19--49.
Engl. Transl.:  Amer. Math. Soc. Translation,  {\bf 65} (1952)  81--124.

\refmark{\refGelfondLivre} 
Gel'fond, A{.} O{.} --
{\it  Transcendental and algebraic numbers.}
Gosudarstv{.} Izdat{.} Tehn.-Teor{.} Lit., Moscow, 1952. 
Engl. transl.,
Dover Publications, Inc., New York 1960.

\refmark{\refGramain}
Gramain, F{.} --
 Sur le th\'eor\`eme de Fukasawa--Gel'fond. 
Invent. Math. {\bf 63} \Numero{3} (1981), 495--506.

\refmark{\refGramainWeber}
Gramain, F{.}; Weber, M{.} --
 Computing an arithmetic constant
related to the ring of Gaussian integers. 
Math. Comp. {\bf 44}  \Numero{169} (1985), 
241--250
. Corrigendum:
Math. Comp.  {\bf 48}  \Numero{178}  (1987), 854.

\refmark{\refGras}
Gras, G{.} --
{\it Class field theory. From theory to practice.}
 Springer Monographs in Mathematics, 2002.

\refmark{\refGrimm}
Grimm, C.~A{.} --
A conjecture on consecutive composite numbers.
Amer. Math. Monthly {\bf   76} (1969), 1126--1128. 
 
\refmark{\refGrinspan}
Grinspan, P{.} --
Measures of simultaneous approximation for quasi-periods of abelian
varieties. 
J{.} Number Theory, {\bf 94}  (2002), \Numero 1, 136--176. 
  
\refmark{\refGuinness}
Guiness, I. G{.} --
A sideways look at Hilbert's twenty-three problems.
Notices Amer. Math. Soc. {\bf 47} \Numero{7} (2000), 752--757.

\refmark{\refGuy}
Guy, R{.} --
{\it   Unsolved problems in number theory.}
Problem Books in Mathematics. Unsolved Problems in Intuitive
Mathematics, {\bf I}. Springer-Verlag,   New York, 1981. Second
edition, 1994.

\refmark{\refGyarmati}
Gyarmati, K{.} --
On a problem of Diophantus.
Acta Arith.,  {\bf  97} \Numero{1} (2001),  53--65.

\refmark{\refHall}
Hall, M.~Jr{.} --
The Diophantine equation
$x\sp{3}-y\sp{2}=k$. 
{\it Computers in number theory,} Proc. Sci. Res. Council
Atlas Sympos. \Numero{2}, Oxford, 1969. Academic Press, London,
(1971), 173--198. 

\refmark{\refHardyWright}
Hardy, G.~ H{.}; Wright, E.~M{.} --
{\it  An introduction to the theory of numbers.}
The Clarendon Press, Oxford University Press, New York,  1938.  Fifth edition, 1979.

\refmark{\refHarman}
Harman, G{.} --
{\it Metric number theory.}
London Mathematical
Society Monographs. New Series, {\bf 18}. The Clarendon Press, Oxford
University Press, New York, 1998.
 
\refmark{\refHershonskyPaulin}
Hersonsky, S{.} and  Paulin, F{.} --
Hausdorff dimension
of Diophantine geodesics in negatively curved manifolds. J. Reine angew.
Math. 539  (2001), 29--43.

\refmark{\refHilbert}
Hilbert, D{.} --
Mathematical Problems.
Bull. Amer. Math. Soc. {\bf  37} \Numero{4} (2000), 407--436. Reprinted from 
Bull. Amer. Math. Soc. {\bf  8}   (1902), 437--479.

\refmark{\refClayMathematicalInstitute}
Jackson, A{.} --
Million-dollar Mathematics Prizes Announced.
Notices Amer. Math. Soc. {\bf 47} \Numero{8} (2000), 877--879.
\hfill\break
\hbox{\seventt $<$http{$:$}//www.claymath.org/$>$}

\refmark{\refKontsevichZagier}
Kontsevich, M{.}; Zagier, D{.} --
Periods.
{\it 
Mathematics Unlimited - 2001 and Beyond.} 
Engquist, B.; Schmid, W., Eds., Springer  (2000), 771--808.

\refmark{\refKraus}
Kraus, A{.} --
On the equation $x^p+y^q=z^r$.		
The Ramanujan Journal,		
{\bf  3} \Numero{3} (1999), 315--333.

\refmark{\refLagarias} 
Lagarias, J.C{.} --
On the Normality of Arithmetical Constants.
Experimental Mathematics {\bf 10},  \Numero{3} (2001), 355--368. 
\hfill\break
\hbox{\seventt $<$http{$:$}//arXiv.org/abs/math.NT/0101055$>$}

\refmark{\refLangITN}
Lang, S{.} --
{\it  Introduction to transcendental numbers.}
Addison-Wesley Publishing Co., Reading,  Mass.-London-Don
Mills, Ont{.} 1966.   {\it Collected Papers,} vol.~I, Springer (2000),
396--506.

\refmark{\refLangTNDA}
Lang, S{.} --
Transcendental numbers and Diophantine
approximations. 
Bull. Amer. Math. Soc. {\bf  77} (1971), 635--677.
{\it Collected Papers,} vol.~II, Springer (2000),
1--43.

\refmark{\refLangHDDP}
Lang, S{.} --
Higher dimensional Diophantine problems. 
Bull.
Amer. Math. Soc. {\bf  80} (1974), 779--787.
{\it Collected Papers,} vol.~II, Springer (2000),
102--110.

\refmark{\refLangDistributions}
Lang, S{.} --
Relations de distributions et exemples
classiques.   
{\it S\'eminaire Delange-Pisot-Poitou, 19e ann\'ee: 1977/78,
Th\'eorie des nombres, } Fasc. 2, Exp. \Numero{40}, 6 p.
{\it Collected Papers,} vol.~III, Springer (2000),
102--110.

\refmark{\refLangECDA}
Lang, S{.} --
{\it  Elliptic curves: Diophantine analysis.}
Grundlehren der Mathematischen 
Wissenschaften {\bf  231}. Springer-Verlag,
Berlin-New York, 1978. 

\refmark{\refLangCyclotomic}
Lang, S{.} --
{\it  Cyclotomic fields.}
Graduate Texts in Mathematics,
{\bf  59}. Springer-Verlag, New York-Heidelberg, 1978.

\refmark{\refLangONCDI}
Lang, S{.} --
Old and new conjectured Diophantine inequalities.
Bull. Amer. Math. Soc. (N.S.) {\bf 23} (1990), \Numero 1, 37--75.
{\it Collected Papers,} vol.~III, Springer (2000),
355--393.

\refmark{\refLangEncyclopaedie}
Lang, S{.} --
{\it  Number theory. III. Diophantine geometry.}
Encyclopaedia of Mathematical Sciences, {\bf  60}. 
Springer-Verlag, Berlin, 1991.
Corrected second printing:
{\it  Survey of Diophantine geometry;} 1997.

\refmark{\refLangAlgebra}
Lang, S{.} --
{\it  Algebra. }
 Addison-Wesley
Publishing Co., Reading, Mass.,  1965. Third edition, 1993. 

\refmark{\refLangGazette}
Lang, S{.} --
La Conjecture de Bateman-Horn.
Gazette des Math\'ematiciens,  {\bf 67} (1996), 214--216. 
{\it Collected Papers,} vol.~IV, Springer (2000),
213--216.

\refmark{\refLangevin}
Langevin, M{.} --
Plus grand facteur premier d'entiers en
progression arithm\'etique. 
{\it  S\'eminaire Delange-Pisot-Poitou, 18e
ann\'ee: 1976/77, Th\'eorie des nombres, } Fasc. 1, Exp. \Numero{3}, 7 p.,
Secr\'etariat Math., Paris, 1977. 

\refmark{\refLangevinLuminy}
Langevin, M{.} --
Partie sans facteur carr\'e d'un produit
d'entiers voisins. 
{\it Approximations diophantiennes et nombres transcendants} (Luminy,
1990), de Gruyter, Berlin (1992), 203--214. 

\refmark{\refLangevinRocky}
Langevin, M{.} --
Sur quelques cons\'equences de la conjecture ($abc$) en arithm\'etique et en logique.   Rocky Mount. J. of Math., {\bf 26} (1996), 1031-1042.

\refmark{\refLangevinJNTBx}
Langevin, M{.} --
\'Equations diophantiennes polynomiales \`a hautes multiplicit\'es.
J. ThŽor. Nombres Bordeaux {\bf
13} (2001),  \Numero 1, 211--226.

\refmark{\refLaurentEger} 
Laurent, M{.} --
New methods in algebraic independence.
Gy\H{o}ry, Kalman (ed.) et al., 
{\it  Number theory. Diophantine, computational and
algebraic aspects. }
Proceedings of the international conference, Eger, Hungary, 
July 29--August 2, 1996. Berlin: de Gruyter, (1998), 311--330.

\refmark{\refLaurentSalonique} 
Laurent, M{.} --
Some remarks on the approximation of complex numbers by 
algebraic numbers. 
Bulletin  Greek Math{.} Soc{.} {\bf  42} (1999), 49--57.

\refmark{\refLaurentTokyo} 
Laurent, M{.} --
Diophantine approximation and algebraic independence.
{\it Maison Franco-Jap\-onaise de Tokyo, November 9-11, 1998, 
France-Japan Conference on Transcendental Number Theory
 (Colloque Franco-Japonais de Th\'eorie des Nombres Transcendants), }
Department of Mathematics of Keio University,
{\it  Seminar on Mathematical Sciences},  {\bf 27} (1999), 75--89.

\refmark{\refLeopoldt}
Leopoldt, H.~W{.} --
Zur Arithmetik in abelschen
Zahlk\"orpern. 
J. Reine angew. Math. {\bf  209} (1962),
54--71.

\refmark{\refMahlerZnumbers} 
Mahler, K{.} --
An unsolved problem on the powers of $3/2$. 
J.
Austral. Math. Soc. {\bf 8} (1968), 313--321. 

\refmark{\refMahlerLivre} 
Mahler, K{.} --
{\it  Lectures on transcendental numbers.} 
Lecture Notes in Mathematics, Vol{.} {\bf  546}. 
Springer-Verlag, Berlin-New York, 1976. 

\refmark{\refMahlerSuggestions} 
Mahler, K{.} --
Some suggestions for further research. 
Bull.
Austral. Math. Soc. {\bf 29}  \Numero{1} (1984), 101--108.

\refmark{\refManin} 
Manin, Yu.~I{.} -- {\it  A course in mathematical logic.}
Graduate Texts in Mathematics, Vol. {\bf 53}. Springer-Verlag, 
1977. 

\refmark{\refMasserabc}
Masser, D.~W{.} --
Note on a conjecture of Szpiro.
{\it S\'eminaire sur
les Pinceaux de Courbes Elliptiques (\og \`a la recherche de Mordell
effectif\fg)}; Paris, 1988. Soc.
Math. France, Ast\'{e}risque, {\bf  183} (1990), 19--23.

\refmark{\refMatiyasevich}
Matiyasevich, Yu{.} --
Le dixi\`eme probl\`eme de Hilbert: que peut-on faire avec les \'equations
diophantiennes? 
 {\it  La Recherche de la V\'erit\'e,} coll. L'\'ecriture des
Math\'ematiques, ACL - Les \'Editions du Kangourou (1999), 281--305.
\hfill\break
\hbox{\seventt $<$http{$:$}//logic.pdmi.ras.ru/Hilbert10$>$}

\refmark{\refMatveev}
Matveev, E{.} M{.} --
Explicit lower estimates for rational homogeneous linear 
forms in logarithms of algebraic numbers.  
II,  
Izv{.} Akad{.} Nauk SSSR{.} Ser{.} Mat{.} 
{\bf 64} \Numero{6}  (2000), 125-180.

\refmark{\refBeal}
Mauldin, R.D{.} --
A generalization of Fermat's last theorem: the Beal conjecture and prize
problem. 
Notices Amer. Math. Soc. {\bf 44} \Numero{11} (1997),  
1436--1437.

\refmark{\refMazurTRP}
Mazur, B{.} -- 
The topology of rational points.
Experiment. Math. {\bf  1} \Numero{1} (1992),  35--45.
 
\refmark{\refMazurQDUNT}
Mazur, B{.} --   
Questions of decidability and undecidability in number theory.
J. Symbolic Logic {\bf  59} \Numero{2} (1994),  353--371.

\refmark{\refMazurSTRP}
Mazur, B{.} --   
Speculations about the topology of rational points: an up-date. 
{\it  Columbia University Number Theory Seminar (New York 1992),}
Ast\'erisque {\bf  128} \Numero{4} (1995),  165--181.

\refmark{\refMazurPowersofNumbers}
Mazur, B{.} --   
Questions about powers of numbers.
Notices Amer. Math. Soc. {\bf 47} \Numero{2} (2000), 195--202.

\refmark{\refMetsankyla}
MetsŠnkylŠ, T{.} --  
 Catalan's Conjecture: Another old Diophantine problem solved.  Bull. Amer. Math. Soc.  {\bf 41}  (2004),  43-57.  

\refmark{\refMihailescu}
Mih\u{a}ilescu, P{.} --
Primary units and a proof of Catalan's conjecture, submitted to  Crelle Journal.  \hfill\break
\hbox{\seventt $<$http${:}$//www-math.uni-paderborn.de/$\sim$preda/papers/zips/catcrelle.ps.gz$>$}

\refmark{\refMullerTisserand} 
Muller, J{.}-M{.}; Tisserand, A{.} --  
Towards exact rounding of the elementary functions. 
Alefeld, Goetz (ed.) et al., 
{\it Scientific computing and validated numerics. }
Proceedings of the international symposium on scientific 
computing, computer arithmetic and validated numerics
SCAN-95, Wuppertal, Germany,  September 26-29, 1995. Berlin:
Akademie Verlag. Math{.} Res{.} 90, 59-71 (1996).
 
\refmark{\refNarkiewiczCPNT}
Narkiewicz, W{.} --
{\it  Classical problems in number theory.}
Polish Scientific Publ. {\bf  62} (1986).

\refmark{\refNarkiewiczEATAN}
Narkiewicz, W{.} -- 
{\it Elementary and analytic theory of algebraic numbers. }
Springer-Verlag, Berlin;
PWN---Polish Scientific Publishers, Warsaw, 1974.
 Second edition, 1990.

\refmark{\refNesterenkoPhilippon}
Nesterenko, Y{.} V{.}; Philippon, P{.}, Eds~--
{\it  Introduction to algebraic independence theory. }
Instructional Conference (CIRM Luminy 1997).  
Lecture Notes in Math., {\bf  1752}, Springer, Berlin-New York,
2001.  
 
 \refmark{\refNitaj}
Nitaj, A{.} -- Web site:
\hfill\break
  {\tt 
$<$http{$:$}//www.math.unicaen.fr/${ 
\sim}$nitaj/abc.html$>$}

\refmark{\refOesterle}
{\OE}sterl\'e, J{.} --
Nouvelles approches du \og th\'eor\`eme\fg\ de Fermat.
{\it  S\'{e}m. Bourbaki},  
1987/88, \Numero{694}; Soc. Math. France, Ast\'{e}risque, {\bf 
161--162} (1988), 165--186.

\refmark{\refOort}
Oort, F{.} --
Canonical liftings and dense sets of CM-points. 
{\it Arithmetic Geometry,}, Cortona 1994, Ist. Naz. Mat. F. Severi,
Cambdridge Univ. Press (1997), 228--234.

\refmark{\refPhilipponLW} 
Philippon, P{.} --  
Ind\'ependance et groupes alg\'ebriques.
{\it Number theory} (Montr\'eal, Qu\'e.,
1985), CMS Conf. Proc., {\bf 7},  Amer. Math. Soc.,
Providence, RI, (1987), 279--284.

\refmark{\refPhilipponJNT} 
Philippon, P{.} --  
Une approche m\'ethodique pour la transcendance et 
l'ind\'ependance  alg\'e\-brique de valeurs de fonctions
analytiques.   
J{.} Number Theory {\bf  64} \Numero{2} (1997), 
291--338.

\refmark{\refPhilipponAustralie} 
Philippon, P{.} --  
Quelques remarques sur des questions d'appro\-ximation 
diophantienne. 
Bull{.} Austral{.} Math{.} Soc., {\bf  59}
(1999),  323--334. Addendum, Ibid., {\bf  61}
(2000),  167--169.

\refmark{\refPhilipponGammaunquart}
Philippon, P{.} --  
Mesures d'approximation de valeurs de fonctions analytiques.
Acta Arith.,  {\bf  88} \Numero{2} (1999),  113--127.

\refmark{\refPillai}
Pillai, S.~S{.} --
On the equation $2^x-3^y=2^X+3^Y$.
Bull Calcutta Math. Soc.  {\bf  37} (1945),
 15--20.

\refmark{\refPollingtonVelani}
Pollington, A{.} D{.};  Velani, S{.} L{.} --
On a problem in simultaneous Diophantine approximation: Littlewood's
conjecture. Acta Math{.} {\bf  185} \Numero{2} (2000),  287--306.

\refmark{\refDipendra}
Prasad, Dipendra  --
An analogue of a conjecture of Mazur: a question in Diophantine approximation
on tori.
to appear in Shalika Festschrift, a supplemental volume of the American J. of Math. 
\hfill\break
\hbox{\seventt
\hfill\break
\hbox{\seventt $<$http{$:$}//www.mri.ernet.in/${\scriptstyle
\sim}$mathweb/PB440-0662G-23-699-710.pdf$>$
}
}

\refmark{\refRamachandraAA} 
Ramachandra, K{.} --
Contributions to the theory of transcendental numbers. 
I, II. Acta Arith{.} {\bf  14} (1967/68), 65-72
and 73--88.
  
\refmark{\refRauzy}
Rauzy, G{.} --
{\it Propri\'et\'es statistiques de suites
arithm\'etiques. }
Le Math\'ematicien, {\bf 15}, Collection SUP. Presses
Universitaires de France, Paris, 1976.

\refmark{\refRemond}
R\'emond, G{.} --
In\'egalit\'e de Vojta en dimension sup\'erieure.
 Ann. Scuola Norm. Sup. Pisa Cl. Sci. (4) {\bf 29} (2000), 101--151.

\refmark{\refRhinSmyth}
Rhin, G.; Smyth, C{.} --
On the absolute Mahler measure of polynomials having all zeros in a sector.
Math. Comp. {\bf 64} \Numero{209} (1995),  295--304.
  

\refmark{\refRibenboimCatalan}
Ribenboim, P{.} --
{\it  Catalan's conjecture. Are $8$ and $9$ the
only consecutive powers?} 
Academic Press, Inc., Boston, MA, 1994.

\refmark{\refRibenboimNumbersFriends}
Ribenboim, P{.} --
{\it  My numbers, my friends.} Popular Lectures on Number Theory.
Springer-Verlag, Berlin-Hei\-del\-berg, 2000.

\refmark{\refRivoalCRAS} 
Rivoal, T{.} --
La fonction z\^eta de Riemann prend une infinit\'e de valeurs irrationnelles aux
entiers impairs.
C{.} R{.} Acad{.} Sci{.} Paris S\'er{.} I Math{.}, {\bf  331} (2000),
267--270.  
\hfill\break
\hbox{\seventt $<$http{$:$}//arXiv.org/abs/math.NT/0008051$>$}

\refmark{\refRoySchanuel} 
Roy, D{.} -- 
Sur la conjecture de Schanuel pour les logarithmes de nombres 
alg\'e\-briques.  {\it  Groupe d'\'Etudes sur les Probl\`emes
Diophantiens} 1988-1989, Publ{.} Math{.} Univ{.} P{.} et M{.}
Curie (Paris VI), {\bf  90} (1989), \Numero{6}, 12 p. 

\refmark{\refRoyMatrices} 
Roy, D{.} -- 
Matrices dont les coefficients sont des formes lin\'eaires.
{\it  S\'eminaire de Th\'eorie des Nombres, Paris 
1987--88,} Progr{.} Math., {\bf  81}, Birkh\"auser
Boston, Boston, MA, (1990), 273--281.

\refmark{\refRoyInventiones} 
Roy, D{.} --
Simultaneous approximation in number fields.
Invent{.} Math{.} {\bf 109} \Numero{3} (1992), 547--556.

\refmark{\refRoyActaMath} 
Roy, D{.} --
Points whose coordinates are logarithms of
algebraic numbers on algebraic varieties. 
Acta Math{.} {\bf  175} \Numero{1} (1995),  49--73.

\refmark{\refRoyACVEF} 
Roy, D{.} --
An arithmetic criterion for the values of the exponential
function.
Acta Arith{.},  {\bf 97} \Numero{2} (2001),  
183-194.

\refmark{\refRoyAASNL} 
Roy, D{.} --
Approximation alg\'ebrique simultan\'ee de nombres de Liouville.
Canad. Math. Bull., {\bf 44} \Numero{1}  (2001), 115--120.

\refmark{\refRoyRMTNT}
Roy, D{.} -- 
Results and methods of transcendental number theory.
Notes de la SMC,
{\bf 33} \Numero{2} (2001), 9--14.

\refmark{\refRoyIFAF} 
Roy, D{.} -- Interpolation formulas and auxiliary functions.
J. Number Theory {\bf 94}  (2002), \Numero 2, 248--285.

\refmark{\refRoyNote} 
Roy, D{.} -- Approximation simultanŽe d'un nombre et de son carrŽ.
 C. R. Math. Acad. Sci. Paris  {\bf 336}  \Numero 1 (2003),   1--6.

\refmark{\refSchinzelBilkent} 
Schinzel, A{.} --
The Mahler measure of polynomials.
{\it  Number theory and its applications}, 
Cem Y.~Y{\i}ld{\i}r{\i}m and Serguei A.~Stepanov Ed., 
Lecture Notes in Pure and Applied Mathematics {\bf  204},
 Marcel Dekker Inc., (1999), 171--183.

\refmark{\refSchinzelZassenhaus} 
Schinzel, A{.}; Zassenhaus, H{.} --
A refinement of two theorems of Kronecker. 
Michigan Math{.} J{.} {\bf  12} (1965), 81--85.

\refmark{\refSchlickeweiViola}
Schlickewei, 
H{.} P{.}; Viola, C{.}--
Generalized Vandermonde determinant.
Acta Arith{.} {\bf  95} \Numero{2} (2000), 123--137.

\refmark{\refSchmidtLNDA} 
Schmidt, W{.} M{.}--
{\it  Diophantine approximation.} 
Lecture Notes in Mathematics, {\bf  785}. Springer, Berlin, 
1980. 

\refmark{\refSchmidtBilkent} 
Schmidt, W{.} M{.}--
Heights of algebraic points.
{\it  Number theory and its applications}, 
Cem Y.~Y{\i}ld{\i}r{\i}m and Serguei A.~Stepanov Ed., 
Lecture Notes in Pure and Applied Mathematics {\bf  204},  Marcel Dekker
Inc., (1999),
185--225.

\refmark{\refSchmidtSurveyAA} 
Schmidt, W{.} M{.}--
On continued fractions and Diophantine approximation in power series fields.
 Acta Arith{.} {\bf  95} \Numero{2} (2000), 139--166.

\refmark{\refSchneiderLivre} 
Schneider, Th{.} --
{\it  Einf\"uhrung in die transzendenten Zahlen.}
Springer-Verlag, Berlin-G\"ot\-tingen-Heidelberg, 1957. 
{\it  Introduction aux nombres transcendants.}
Tra\-duit de l'alle\-mand par 
 P{.} Eymard. Gauthier-Villars, 
Paris 1959. 

\refmark{\refSerre}
Serre, J-P{.}   -- 
{\it Lectures on the Mordell-Weil theorem}.
Aspects of Mathematics, {\bf E15}. Friedr. Vieweg {\&} Sohn,
Braunschweig, 1989.

\refmark{\refShoreySurvey}
Shorey, T{.} N{.}   -- 
Exponential Diophantine equations involving products of
consecutive integers.
{\it  Number Theory,} 
 R.P. Bambah, V.C. Dumir and R.J. Hans Gill, Eds,
Hindustan Book Agency, New-Delhi  and  Indian National Science Academy 
(1999), 463--495.

\refmark{\refShoreyConjectures}
Shorey, T{.} N{.}   -- 
Some conjectures in the theory of exponential Diophantine equations. 
 Publ. Math. Debrecen {\bf 56 } (2000), \Numero 3-4, 631--641.

\refmark{\refShoreyTijdeman} 
Shorey, T{.} N{.}; Tijdeman, R{.}  -- 
{\it  Exponential Diophantine equations.} 
Cambridge Tracts in Mathematics, {\bf  87}. Cambridge 
University Press, Cambridge-New York, 1986. 

\refmark{\refSiegelPreuss} 
Siegel, C{.} L{.} --
\"Uber einige Anwendungen diophantischer Approximationen.
Abh{.} Preuss{.} Akad{.} Wiss., Phys.-Math., {\bf 1}
(1929), 1--70.
{\it Gesammelte Abhandlungen.}   Springer-Verlag,
Berlin-New York 1966 Band {\bf I}, 209--266. 

\refmark{\refSiegelLivre} 
Siegel, C{.} L{.} --
{\it Transcendental numbers.} 
Annals of Mathematics Studies, \Numero{\bf 16}.
Princeton University Press,
Princeton, N{.} J., 1949. 

\refmark{\refSierpinskiA}
Sierpi\'nski, W{.} --
{\it A selection of problems in the theory of numbers.}
Translated from the Polish by A. Sharma. 
Popular lectures in mathematics, {\bf 11}.
A Pergamon Press Book The Macmillan Co., New York 1964 

\refmark{\refSierpinskiB}
Sierpi\'nski, W{.} --
{\it 250 problems in elementary number theory.}
Elsevier, 1970.
Modern Analytic and Computational Methods in Science and Mathematics,
\Numero{\bf 26} American Elsevier Publishing Co., Inc., New York; PWN
Polish Scientific Publishers, Warsaw 1970.
{\it 250 probl\`emes de th\'eorie \'el\'ementaire
des nombres. }
Translated from the English. Reprint of the 1972 French translation.
\'Editions Jacques Gabay, Sceaux, 1992.

\refmark{\refSilverman}
Silverman, J{.} H{.} --
{\it The arithmetic of elliptic curves.}
Graduate Texts in Mathematics, {\bf 106}. 
Springer-Verlag, New York-Berlin, 1986.

\refmark{\refSmale}
Smale, S{.} --
Mathematical problems for the next century.
Math. Intelligencer {\bf 20} \Numero{2}  (1998),  7--15.
{\it Mathematics: frontiers and perspectives},  271--294, Amer. Math. Soc.,
Providence, RI, 2000. 

\refmark{\refSondow}
Sondow, J{.} --
Criteria for Irrationality of Euler's Constant.
Proc. Amer. Math. Soc. {\bf 131} (2003), 3335-3344.
\hfill\break
\hbox{\seventt 
$<$http{$:$}//arXiv.org/abs/math.NT/0209070$>$}

\refmark{\refSprindzuk}
Sprind\v{z}uk, V{.} G{.} --
{\it Metric theory of Diophantine
approximations. }
Scripta Series in
Mathematics. V. H. Winston {\&} Sons, Washington, D.C.; A Halsted Press Book,
John Wiley {\&} Sons, New York-Toronto, Ont.-London, 1979. 

\refmark{\refStewartYuKunruiI}
Stewart, C{.} L{.}; Yu, Kun Rui --
On the $abc$ conjecture. 
Math{.} Ann{.} {\bf  291} \Numero{2} (1991),
 225--230. 

\refmark{\refStewartYuKunruiII}
Stewart, C{.} L{.}; Yu, Kun Rui --
On the $abc$ conjecture. II, Duke Math. J., {\bf 108} (2001), 169--181.

\refmark{\refTerasoma}
Terasoma, T{.} -- 
Mixed Tate motives and multiple zeta values.
Invent{.} Math{.} {\bf 149} (2002), \Numero{2}, 339-369

\refmark{\refThakur}
Thakur, D{.} --
Automata and transcendence.
{\it Number Theory},
Ramanujan Mathematical Society, January 3-6, 1996, 
Tiruchirapalli, India;
V.~Kumar Murty and Michel Waldschmidt, Eds.,
Amer. Math. Soc., Contemporary Math.  {\bf 210} (1998), 387--399.

\refmark{\refTijdemanHilbert} 
Tijdeman, R{.} --
Hilbert's seventh problem: on the Gel'fond-Baker
method and its applications. 
{\it  Mathematical developments arising from Hilbert
problems,} Proc. Sympos. Pure Math., {\bf  28}, Part 1,  Northern
Illinois Univ., De Kalb, Ill., 1974. Amer. Math. Soc.,
Providence, R. I., (1976),  241--268. 

\refmark{\refTijdemanCatalan} 
Tijdeman, R{.} --
On the equation of Catalan. 
Acta Arith{.} {\bf  29} \Numero{2} (1976), 197--209.

\refmark{\refTijdemanSurveyA}
Tijdeman, R{.} -- Diophantine equations and Diophantine approximation. 
{\it Number theory and applications,} Banff, AB, 1988, NATO Adv. Sci.
Inst. Ser. C: Math. Phys. Sci., 265, Kluwer Acad. Publ., Dordrecht, (1989),
215--243.
 
\refmark{\refTijdemanSurveyB} 
Tijdeman, R{.} --
Exponential Diophantine equations 1986--1996. 
{\it  Number theory,}
Eger, 1996, de Gruyter, Berlin, (1998),  523--539. 

\refmark{\refTijdemanSADA} 
Tijdeman, R{.} --
Some applications of Diophantine approximation.
Math. Inst. Leiden, Report MI 2000-27, September 2000, 19p.

\refmark{\refVojta}
Vojta, P{.} -- 
On the $ABC$ Conjecture and Diophantine approximation by
rational points. 
Amer. J. Math. {\bf 122} \Numero{4} (2000), 843--872.

\refmark{\refmiwSemLelong} 
Waldschmidt, M{.} --
Propri\'{e}t\'{e}s arithm\'{e}tiques de fonctions de plusieurs variables (II).
S\'{e}m. P.  Lelong (Ana\-lyse), 16\`{e}  ann\'{e}e, 1975/76. Lecture Notes  
in Math., {\bf567} (1977), 274--292.

\refmark{\refmiwGAGDT} 
Waldschmidt, M{.} --
Groupes alg\'{e}briques et grands degr\'{e}s de transcendance.
Acta Mathematica, {\bf  156} (1986), 253--294. 

\refmark{\refmiwNagoya} 
Waldschmidt, M{.} --
Density measure of rational points on abelian varieties. 
Nagoya Math. J., {\bf  155} (1999), 27--53.

\refmark{\refmiwSurveyIA} 
Waldschmidt, M{.} --
Algebraic independence of transcendental numbers: a survey.
{\it  Number Theory,} 
 R.P. Bambah, V.C. Dumir and R.J. Hans Gill, Eds,
Hindustan Book Agency, New-Delhi  and  Indian National Science Academy 
(1999), 497--527.
\hfill\break
\hbox{\seventt $<$http{$:$}//www.birkhauser.ch/books/math/6259.htm$>$}

\refmark{\refmiwGraz} 
Waldschmidt, M{.} --
Conjectures for Large Transcendence Degree.
{\it  Algebraic Number Theory and Diophantine Analysis,} 
F. Halter-Koch and R. F. Tichy, Eds, W{.} de Gruyter, Berlin, (2000), 497-520.

\refmark{\refmiwGL}
Waldschmidt, M{.} --
{\it  Diophantine Approximation on linear algebraic groups. 
Transcendence Properties of the Exponential Function in Several 
Variables.} 
Grund\-lehren der Mathematischen Wissenschaften
{\bf  326}, Springer-Verlag,
Berlin-Hei\-del\-berg, 2000. 

\refmark{\refmiwCinquante}
Waldschmidt, M{.} --
Un demi-si\`{e}cle de transcendance. 
{\it  Development of Mathematics 1950--2000,}
\sanserif
 Pier, J.-P., ed., Birkh\"{a}user Verlag, (2000), 
1121--1186.
\hfill\break
\hbox{\seventt $<$http{$:$}//www.birkhauser.ch/books/math/6280.htm$>$}
\hfill
\hbox{\seventt $<$http{$:$}//www.math.jussieu.fr/${\scriptstyle
\sim}$miw/articles/50.html$>$}

\refmark{\refmiwMZV}
Waldschmidt, M{.} --
Valeurs z\^eta multiples: une introduction.
 J. ThŽor. Nombres Bordeaux, 
{\bf 12} (2000), 581--595.

\refmark{\refmiwCetraro}
Waldschmidt, M{.} --
Linear Independence Measures for Logarithms of Algebraic
Numbers. 
`` Diophantine Approximation, Cetraro, Italy 2000'', 
Lecture Notes in Mathematics
 {\bf 1819 }  (2003), 249--344.  
\hfill\break
\hbox{\seventt $<$http{$:$}//www.math.jussieu.fr/${\scriptstyle
\sim}$miw/articles/cetraro.html$>$}

\refmark{\refWaring}
Waring, E{.} --
{\it  Meditationes algebraic{\ae}.}
Cambridge, 1770. English Translation, Amer. Math. Soc. 1991.

\refmark{\refWhittakerWatson} 
Whittaker, E.T.; Watson, G.N{.} --
{\it A Course of modern analysis.}
 Cambridge Univ. Press, 1902. Fourth edition, 1927.

\refmark{\refWirsing}
Wirsing, E{.} --
Approximation mit algebraischen Zahlen beschr\"ankten Grades. 
J. Reine angew. Math. {\bf 206} \Numero{1/2} (1961), 
67--77.

\refmark{\refWong}
Wong, Chi Ho --
 An explicit result related to the $abc$-conjecture, MPhil Thesis, Hong Kong University of Science  \&\ Technology, 1999.

\refmark{\refWoods}
Woods, A{.} -- 
{\it Some problems in logic and number theory.}
{\it Thesis, } Manchester, 1981.

\refmark{\refZagier} 
Zagier, D{.} --
Values of zeta functions and their applications. 
{\it Proc. First European Congress of Mathematics,}
Vol. {\bf 2}, Birkh\"auser, Boston (1994), 497--512.

\refmark{\refZudilin} 
Zudilin, W{.} --
Algebraic relations of multiple zeta values.
Uspekhi Mat. Nauk =
Russian Math. Surveys {\bf 58}:1 (2003), 3--32.

}

\vfill

 \vskip 2truecm plus .5truecm minus .5truecm 

\hfill
\vbox{\ninerm
 \hbox{Michel WALDSCHMIDT}
 \hbox{Universit\'e P.~et M.~Curie (Paris VI)}
 \hbox{Institut de Math\'ematiques  de Jussieu CNRS UMR 7586}
 \hbox{Th\'eorie des Nombres\qquad Case 247}
 \hbox{175, rue du Chevaleret }
 \hbox{F--75013 PARIS}
 \hbox{e-mail: {\ninett miw@math.jussieu.fr}}
 \hbox{URL: {\ninett
$<$http{$:$}//www.math.jussieu.fr/${  \sim}$miw/$>$}}
 \hbox{ \qquad {\ninett
$<$http$:$//arXiv.org/abs/math.NT/0312440$>$}}
} 

\bye